\documentclass[10pt,english]{amsart}
\usepackage{etoolbox}
\patchcmd{\thebibliography}{*}{}{}{}
\pretocmd\thebibliography{\csname c@secnumdepth\endcsname=-2 }{}{}

\usepackage[T1]{fontenc}

\usepackage[margin=8em]{geometry}

\usepackage{amssymb,url,xspace,aliascnt,stmaryrd,mathtools}
\usepackage{accents}
\usepackage{enumitem} %resuming enumerate
\setlist[enumerate]{label=(\roman*)} %labelling for enumerate
\usepackage{soul} %highlight

\usepackage{indextools}              % table de notation
\indexsetup{level=\subsection*}
\makeindex[columns=3, title=Index of notation]

\usepackage{tikz-cd,denpastyle}

\newaliascnt{theo}{subsubsection}
\newaliascnt{lemm}{subsubsection}
\newaliascnt{prop}{subsubsection}
\newaliascnt{coro}{subsubsection}
\newaliascnt{rema}{subsubsection}
\newaliascnt{defi}{subsubsection}
\newaliascnt{exam}{subsubsection}
\newaliascnt{conj}{subsubsection}
\newaliascnt{hypo}{subsubsection}
\newaliascnt{quest}{subsubsection}
\theoremstyle{plain}
\newtheorem{theointro}{Theorem}

\newtheorem{theo}[theo]{Theorem}
\newtheorem{lemm}[lemm]{Lemma}
\newtheorem{prop}[prop]{Proposition}
\newtheorem{coro}[coro]{Corollary}
\theoremstyle{remark}
\newtheorem{rema}[rema]{Remark}
\theoremstyle{definition}
\newtheorem{defi}[defi]{Definition}
\newtheorem{exam}[exam]{Example}

\aliascntresetthe{theo}
\aliascntresetthe{lemm}
\aliascntresetthe{prop}
\aliascntresetthe{coro}
\aliascntresetthe{rema}
\aliascntresetthe{defi}
\aliascntresetthe{exam}
\aliascntresetthe{conj}
\aliascntresetthe{hypo}
\aliascntresetthe{quest}

\let\originalleft\left
\let\originalright\right
\renewcommand{\left}{\mathopen{}\mathclose\bgroup\originalleft}
\renewcommand{\right}{\aftergroup\egroup\originalright}
\newcommand{\FS}{\operatorname{FS}}
\newcommand{\RHom}{\operatorname{RHom}}
\newcommand{\Db}{\operatorname{D}^{\mathrm{b}}}

\newcommand{\DG}{\operatorname{DG}}
\newcommand{\Perv}{\mathrm{Perv}}

\newcommand{\nil}{\mathrm{nil}}

\newcommand{\cont}{\mathrm{cont}}
\DeclareMathOperator{\gHom}{\mathrm{gHom}}

\newcommand{\ext}{\mathrm{ext}}

\newcommand{\ubar}[1]{\underline{\smash{#1}}}
\newcommand{\Gz}{{G_{\ubar 0}}}

\newcommand{\Cq}{{\bbC^{\times}_{q}}}
\newcommand{\LY}{\mathrm{LY}}
\newcommand{\Sp}{\mathrm{Sp}}
\newcommand{\KL}{\mathrm{KL}}

\newcommand{\sgn}{\prescript{\epsilon}{}}
\newcommand{\ci}{\accentset{\circ}}
\def\pH{\prescript{p}{}{\mathscr{H}}}
\newcommand{\fl}{\mathrm{fl}}
\usepackage[backend=biber,doi=false,url=false,isbn=false]{biblatex}
\DeclareFieldFormat[article]{volume}{\textbf{#1}}
\DeclareFieldFormat[article]{number}{-#1}
\DeclareFieldFormat[book]{volume}{\textbf{#1}}
\renewbibmacro{in:}{%
  \ifentrytype{article}{}{\printtext{\bibstring{in}\intitlepunct}}}

\renewbibmacro*{volume+number+eid}{%
  \printfield{volume}
  \printfield[parens]{number}%
  \setunit{\addcomma\space}%
  \printfield{eid}}

\DeclareFieldFormat{postnote}{#1}
\DeclareFieldFormat{multipostnote}{#1}
\title{Bi-orbital sheaves and affine Hecke algebras at roots of unity}
\date{\today}
\author{\sc Wille Liu}

\address{{\sc Institute of Mathematics, Academia Sinica}.  \\
	6F, Astronomy-Mathematics Building, 
	No. 1, Sec. 4, Roosevelt Road, 
	Taipei, Taiwan.
}
\email{\href{mailto:wliu@sinica.edu.tw}{wliu@sinica.edu.tw}}

\begin{document}
\def\smfbyname{}
\maketitle

\begin{abstract}
	Continuing the study of perverse sheaves on the nilpotent cone of a $\bbZ/m$-graded Lie algebra initiated by Lusztig--Yun, we study in this work the parabolic induction and introduce the notion of supercuspidal sheaves on the nilpotent cone. One of our main results shows that simple perverse sheaves with nilpotent singular support (called bi-orbital sheaves) are produced by parabolic induction from supercuspidal sheaves. As application, we provide a proof for a theorem announced by I. Grojnowski on the parametrisation of simple modules of affine Hecke algebras at roots of unity via Deligne--Langlands--Lusztig parameters with nilpotent singular support.
\end{abstract}
%\tableofcontents

\section*{Introduction}

\subsection*{Sheaf theory on graded Lie algebras}

Let $G$ be a connected reductive complex algebraic group and $\frakg$ its Lie algebra. Let $m\in \bbZ_{>0}$ and let $\theta:\mu_m\to \Aut(G)$ be a group homomorphism, where $\mu_m\subset \bbC^{\times}$ is the subgroup of $m$-th roots of unity. Then $\theta$ gives rise to a $\bbZ/m$-grading on $\frakg$, written $\frakg = \bigoplus_{\ubar n\in \bbZ/m} \frakg_{\ubar n}$ with $\frakg_{\ubar n} = \left\{ X\in \frakg \;\vert\; \theta(\omega)X = \omega^nX\;\; \forall \omega\in \mu_m \right\}$. Let $\Gz = (G^{\theta})^{\circ}$ be the neutral component of the subgroup of $\theta$-fixed points. Fix an integer $d\in \bbZ\setminus\left\{ 0 \right\}$ and denote $\ubar d = d\mod m\in \bbZ / m$. The adjoint action of $G$ on $\frakg$ restricts to an action of $\Gz$ on $\frakg_{\ubar d}$. The $\Gz$-equivariant derived category of sheaves on the nilpotent cone $\frakg^{\nil}_{\ubar d} = \frakg_{\ubar d}\cap \frakg^{\nil}$ with bounded constructible cohomology, denoted by $\Db_{\Gz}(\frakg_{\ubar d}^{\nil})$ is the object of our study. \par
The invariant theory of the adjoint $\Gz$-action on $\frakg_{\ubar d}$ has been first studied by Vinberg \cite{vinberg76}. This theory, known as Vinberg's $\theta$-groups, reveals a rich connection with complex reflection groups. More recently, Vinberg's theory has drawn considerable attention for its connection with the representation theory of $p$-adic groups, see \cite{RY14,fintzen21} and references therein.  In view of the success of the Springer theory and the theory of character sheaves of Lusztig, it is natural to expect that the sheaf theory on $\frakg_{\ubar d}$ plays a role in the representation theory of $p$-adic groups. Some observations in this sense have been made in \cite[introduction]{VV09}.\par

In~\cite{LYI}, Lusztig and Yun discovered a uniform way to produce all the simple $\Gz$-equivariant perverse sheaves on $\frakg_{\ubar d}^{\nil}$, which generate $\Db_{\Gz}(\frakg^{\nil}_{\ubar d})$ as thick subcategory. They remarked that in contrast to the $\bbZ$-graded case studied in~\cite{lusztig95b}, the parabolic induction is not enough to produce all $\Gz$-equivariant simple perverse sheaves on $\frakg_{\ubar d}^{\nil}$. Instead, they introduced the notion of spiral induction and spiral restriction, which turn out to be the right notion of induction and restriction of sheaves on $\bbZ/m$-graded Lie algebras as they produce all $\Gz$-equivariant simple perverse sheaves on $\frakg^{\nil}_{\ubar d}$ and share many of the features of parabolic induction and restriction on $\bbZ$-graded Lie algebras. \par

According to~\cite{LYI}, there is a decomposition of triangulated category
\begin{eq}
	\Db_{\Gz}\left( \frakg_{\ubar d}^{\nil} \right) = \bigoplus_{\xi\in \ubar \frakT} \Db_{\Gz}\left( \frakg_{\ubar d}^{\nil} \right)_{\xi}.
\end{eq}
Let $\frakT$ be the set of systems $(L, \frakl_*, \scrC)$, called \emph{Lusztig--Yun (LY)-cuspidal support}, where $L$ is a \emph{pseudo-Levi subgroup} of $(G, \sigma)$, where $\sigma$ is the \emph{outer type} of $\theta$,  $\frakl_*$ is a $\bbZ$-grading on the Lie algebra $\frakl = \Lie L$ such that $\frakl_n\subseteq \frakg_{\ubar n}$ for all $n\in \bbN$, and $\scrC$ is a cuspidal local system on the open $L_0$-orbit in $\frakl_d$ in the sense of~\cite{lusztig95b}, where $L_0$ is the connected subgroup of $G$ with Lie algebra $\frakl_0$. The set above $\ubar\frakT$ is the set of $\Gz$-conjugacy classes in $\frakT$. The sum of the spiral induction of the intersection complex $\IC(\scrC)$ with coefficients in $\scrC$ along various spirals forms an (ind-)complex $\bfI^{\LY}$ in the equivariant category $\Db_{\Gz}( \frakg^{\nil}_{\ubar d} )$. The subcategory $\Db_{\Gz}( \frakg_{\ubar d}^{\nil} )_{\xi}$, called the \emph{$\xi$-block}, is the thick subcategory generated by $\bfI^{\LY}$. The simplest example is the \emph{principal block}: it is the block associated with $\xi_0 = \left(T, \frakt_*, \delta_0  \right)$, where $T$ is a maximal torus of $\Gz$ and $\delta_0$ is the skyscraper sheaf supported on $0$ with fibre $\bbC$. \par

%Let $\Cq$ be a one-dimensional torus acting linearly on $\frakg$ by weight $2$. The action of $\Cq$ commutes with the adjoint action of $G$ so that there is a $\Gz\times \Cq$ action on $\frakg^{\nil}_{\ubar d}$. By Jacobson--Morosov theorem, the $\Gz$-orbits in $\frakg^{\nil}_{\ubar d}$ are stable by the action of $\Cq$. The functor of forgetting the $\Cq$-action 
%\begin{eq}
	%\Db_{\Gz\times\Cq}(\frakg^{\nil}_{\ubar d})\to \Db_{\Gz}(\frakg^{\nil}_{\ubar d}) \quad \text{induces}\quad \Irr\Perv_{\Gz}(\frakg_{\ubar d}^{\nil}) \cong\Irr\Perv_{\Gz\times\Cq}(\frakg_{\ubar d}^{\nil})
%\end{eq}
%on the isomorphism classes of equivariant simple perverse sheaves. Similarly, each admissible system $\xi\in \frakT\left( \frakg_{\ubar d} \right)$ acquires automatically a $\Gz\times \Cq$-equivariant structure. The spiral induction, spiral restriction and the decomposition of equivariant category into blocks also has a $\Gz\times\Cq$-equivariant version. \par

Fix an LY-cuspidal support $\xi\in \frakT$ on $\frakg_{\ubar d}$. From a Springer-theoretic point of view, the conjecture of \cite{LYIII} proven in \cite{liu21} states that the simple perverse sheaves in the $\xi$-block are in bijective correspondence with simple modules of a certain block of the category of integrable modules of a certain degenerate double affine Hecke algebras (dDAHA, a.k.a. trigonometric Cherednik algebras) $\bfH$, depending on $\xi$. The completed extension algebra $\calH = \Hom^{\bullet}_{G_{\ubar 0}}\left(\bfI^{\LY}, \bfI^{\LY}  \right)^{\wedge}$ realises a localised version of $\bfH$, called spectral completion (\autoref{subsec:spectral}), which controls a block of the category of integrable modules $\rmO(\bfH)$. The theory of extension algebras studied in Chriss--Ginzburg~\cite{CG} and further developed by Kato \cite{kato17} allows one to analyse the structure of the block of the category of $\calH$-modules via the geometry of $\frakg_{\ubar d}^{\nil}$ and the complex $\bfI^{\LY}$. In particular, the simple $\calH$-modules are in bijective correspondence with simple constituents of $\bfI^{\LY}$ and the dimension of a simple module is equal to the multiplicity of the corresponding simple factor of $\bfI^{\LY}$. Moreover, there is a so-called \emph{standard module} for each simple module of $\calH$, defined in terms of the $!$-fibre of $\bfI^{\LY}$ at a nilpotent $\Gz$-orbit in $\frakg_{\ubar d}^{\nil}$. \par

\subsection*{Deligne--Langlands correspondence at roots of unity}

Let $\bfK$ denote the affine Hecke algebra in the sense of \cite{lusztig89}. In the case where the parameter of $\bfK$ are of infinite order, Kazhdan--Lusztig \cite{KL} established the famous Deligne--Langlands correspondence for affine Hecke algebras: simple $\bfK$-modules are parametrised by certain $G$-conjugacy classes of triplets $(s, z, \chi)$ consisting of semisimple element $s\in G$, nilpotent element $z\in \frakg$ and an irreducible representation $\chi$ of the component group of the simultaneous stabiliser $Z_G(s,z)$ (Deligne--Langlands--Lusztig parameters), based on Ginzburg's realisation of $\bfK$ via the equivariant K-homology on the Steinberg variety over $\frakg$ and equivariant localisation. The proof can be formulated in terms of equivariant perverse sheaves on a subspace of $\frakg$, see \cite{CG,lusztig88}. This theorem can fail when the parameter of $\bfK$ is a root of unity. In this case, Grojnowski has announced a theorem \cite[Theorem 1]{groj94} describing subset of the data $(s, z, \chi)$ which parametrise finite-dimensional $\bfK$-modules. This subset is defined in terms of a nilpotency condition on the singular support of a perverse sheaf attached to the datum $(s, z, \chi)$. The proof of it seems to be unavailable to this day. \par

The set of data $(s, z, \chi)$ with a fixed $s$ of finite-order can be viewed as simple local systems on $\Gz$-orbits of a graded Lie algebra $\frakg_*$, where the grading $\theta$ is given by the adjoint action of $s$. Let $\bfI^{\KL}\in \Db_{\Gz}(\frakg_{\ubar d})$ denote the sum of parabolic induction of the constant sheaf along various $\theta$-stable Borel subalgebras. The semisimple complex $\bfI^{\KL}$ is a special case (when $\xi = \xi_0$) of parabolic induction from a \emph{supercuspidal support}. The completed extension algebra $\calK = \hom_{G_{\ubar 0}}(\bfI^{\KL}, \bfI^{\KL})^{\wedge}$ is isomorphic to a localised version of $\bfK$ (with parameters at roots of unity) which accounts for a block of finite-dimensional $\bfK$-modules. Again, the theory of extension algebras implies that the simple $\calK$-modules are in bijective correspondence with simple constituents of $\bfI^{\KL}$. \par

The question of the relation between $\bfI^{\LY}$ and $\bfI^{\KL}$ arises naturally in this context. On one hand, $\bfI^{\KL}$ lies in the $\xi$-block $\Db_{\Gz}(\frakg^{\nil}_{\ubar d})_{\xi}$, which is generated by $\bfI^{\LY}$; on the other hand, $\bfI^{\KL}$ may not generate $\Db_{\Gz}(\frakg^{\nil}_{\ubar d})_{\xi}$ as thick subcategory in general. Grojnowski's statement \cite[Theorem 1]{groj94} can be reformulated in terms of Lusztig--Yun's theory as follows: the simple constituents of $\bfI^{\KL}$ are precisely the bi-orbital sheaves in the $\xi_0$-block. One of the objectives of the present article is to prove the following generalisation of \cite[Theorem 1]{groj94}: bi-orbital sheaves are simple constituents of parabolic induction of bi-orbital \emph{supercuspidal sheaves}. The case where $\xi = \xi_0$ recovers \cite[Theorem 1]{groj94}. \par

Our approach is based on the representation theory of dDAHA. The Yoneda product provides a $(\calK, \calH)$-bimodule structure on the Ext-space $\Hom^{\bullet}_{\Gz}(\bfI^{\LY}, \bfI^{\KL})^{\wedge}$. We define a functor  
\begin{eq}
	\bfV = \Hom^{\bullet}_{\Gz}(\bfI^{\LY}, \bfI^{\KL})^{\wedge}\otimes_{\calH}\relbar: \calH\Mod\to \calK\Mod.
\end{eq}
It is identified with the algebraic Knizhnik--Zamolodchikov functor introduced in \cite{liu22}. The characterisation of the kernel of $\bfV$ (\autoref{theo:KZ}) and the algebraic characterisation of orbital and anti-orbital sheaves (\autoref{prop:proj-inj}) are the key ingredients of the proof.

\subsection*{Main results}
The first half of this article is dedicated to the study of parabolic induction and \emph{supercuspidal sheaves} in the $\bbZ/m$-graded setting. We introduce in \autoref{subsec:parab-ind} the $\bbZ/m$-graded version of parabolic induction and restriction. Let $\Pi$ denote the set of pairs $(\rmC, \scrL)$ consisting of a nilpotent $\Gz$-orbit $\rmC\subseteq \frakg^{\nil}_{\ubar d}$ and an irreducible $\Gz$-equivariant local system $\scrL$ on $\rmC$, defined up to isomorphism. The notion of a \emph{supercuspidal pair} is introduced in \autoref{subsec:supercuspidal}: a pair $(\rmO, \scrS)\in \Pi$ is supercuspidal if it is annihilated by parabolic restriction along every proper parabolic subalgebra. 
\begin{theointro}[=\autoref{theo:supercuspidal}]\label{theo:A}
	The following conditions for the pair $(\rmO, \scrS)\in \Pi$ are equivalent:
	\begin{enumerate}
		\item
			$(\rmO, \scrS)$ is a supercuspidal pair;
		\item
			the orbit $\rmO$ is distinguished and the local system $\scrS$ is clean (\autoref{subsec:primitive});
		\item
			$(\rmO, \scrS)$ is the primitive pair (\autoref{subsec:primitive}) attached to an LY cuspidal support $\xi = (L, \frakl_*, \scrC)$ and $(Z_G)^{\theta} / Z_L$ is finite.
	\end{enumerate}
\end{theointro}
A \emph{supercuspidal support} on $\frakg_{\ubar d}$, introduced in \autoref{subsec:supp-supercuspidal} is a triplet $(M, \rmO, \scrS)$, where $M\subseteq G$ is a $\theta$-stable Levi factor of a $\theta$-stable Levi subgroup and $(\rmO, \scrS)$ is a supercuspidal pair on $\frakm_{\ubar d}$. \autoref{theo:A} provides a canonical supercuspidal support attached to each LY cuspidal support, which reduces the classification of supercuspidal sheaves to that of cuspidal local systems on a reductive group due to Lusztig \cite{lusztig84}. \par
Given a $\theta$-stable parabolic subgroup $Q\subseteq G$ let $\Ind^{\frakg_{\ubar d}}_{\frakq_{\ubar d}}$ denote the parabolic induction along $\frakq$. 
\begin{theointro}[=\autoref{theo:comparaison-bloc}]\label{theo:B}
	If $\zeta = (M, \rmO, \scrS)$ is a supercuspidal support attached to an LY cuspidal support $\xi\in\frakT$, then for every pair $(\rmC, \scrL)\in \Pi$ whose intersection complex $\IC\left(\scrL\right)$ lies in the $\xi$-block, there exists a $\bbZ/m$-graded parabolic subgroup $Q\subseteq G$ which contains $M$ as Levi factor such that the local system $\scrL$ is a direct summand of the cohomology of the complex $j^!_{\rmC}\Ind^{\frakg_{\ubar d}}_{\frakq_{\ubar d}}\IC\left( \scrS \right)$, where $j_{\rmC}:\rmC\hookrightarrow \frakg^{\nil}_{\ubar d}$ is the inclusion of nilpotent orbit. 
\end{theointro}
In other word, even if in the $\bbZ/m$-setting the parabolic induction from supercuspidal sheaves cannot produce all simple perverse sheaves on the nilpotent cone as simple constituents, the cohomology of ``Springer fibres'' with coefficients in $\scrS$ for various supercuspidal support still contains all local systems on nilpotent orbits as direct summands. \par

The second half of this article is dedicated to the proof of the aforementioned generalisation of \cite[Theorem 1]{groj94}. A simple $\Gz$-equivariant complex on $\frakg_{\ubar d}$ is called \emph{bi-orbital} if its singular support lies in $\frakg^{\nil}_{\ubar d}\times \frakg^{\nil}_{-\ubar d}\subseteq T^*\frakg_{\ubar d}$. The last main result is the following:
\begin{theointro}[=\autoref{theo:biorbital}] \label{theo:C}
	Let $\scrF$ be a $\Gz$-equivariant simple perverse sheaf on $\frakg^{\nil}_{\ubar d}$. Then, $\scrF$ is bi-orbital if and only if there exist a supercuspidal support $(M, \rmO, \scrS)$ and a $\theta$-stable parabolic subgroup $Q\subseteq G$ with Levi factor $M$ such that $\IC(\scrS)$ is bi-orbital and $\scrF$ is a direct summand of the perverse cohomology of $\Ind^{\frakg_{\ubar d}}_{\frakq_{\ubar d}}\IC(\scrS)$. 
\end{theointro}
This theorem provides an approach to a question of G. Lusztig \cite[0.6(a)]{lusztig11}: describe explicitly bi-orbital simple perverse sheaves on $\frakg_{\ubar d}$. \par
We derive two consequences from \autoref{theo:C}. We show in \autoref{subsec:cusp} that bi-orbital \emph{cuspidal} sheaves are necessarily supercuspidal and we prove \cite[Theorem 1]{groj94} in \autoref{subsec:groj}.

\subsection*{Related works}
Perverse sheaves on $\bbZ/2$-graded Lie algebras have been extensively studied from the perspective of character sheaves on symmetric spaces. The reader is referred to \cite{shoji19} for an overview. \par
Apart from the works of Lusztig--Yun, there are recent works by Grinberg, Vilonen and Xue \cite{GVX23,VX21,VX23} about construction of cuspidal character sheaves on $\bbZ/m$-graded Lie algebras. These works can be thought of as \emph{orthogonal} to the results of the present article. \par
Hennecart's result \cite{hennecart22} about the characterisation of Lusztig's sheaf via singular support can be viewed as analogue of \cite[Theorem 1]{groj94} for affine and loop quivers. 

\section*{Acknowledgements}
This article is a continuation of my doctoral thesis under the supervision of \'Eric Vasserot, to whom I owe the most thanks. I am also grateful to Tsai Cheng-Chiang, who told me about the relation between parabolic inductions and spiral inductions, which inspired this article.

\section*{Convention}
In the present article, all algebraic varieties are defined over a fixed algebraically field $\bfk$ of characteristic $0$ or of prime characteristic larger than the Coxeter number of algebraic groups that are in question. \par

Given a (possibly non-unital) ring $A$, we let $A\Mod$ denote the category of left $A$-modules satisfying $AM = M$, $A\mof$ the full subcategory of finitely presented modules, $A\mof_{\fl}$ the full subcategory of modules of finite length and $A\proj$ the full subcategory of finitely generated projective modules.  \par

Given a linear algebraic group $G$, let $G^{\circ}$ denote the neutral component of $G$ and $Z_G$ the centre of $G$. The Lie algebra of $G$ is denoted with the fraktur small case $\frakg$ and the adjoint representation is denoted by $\Ad: G\to \Aut(G),\; g\mapsto \Ad_g$. We will also write $G^{\ad} = \Ad(G)$ for the adjoint group. Let $\rmH^{\bullet}_G = \rmH^{\bullet}_G(\pt, \bbC)$ denote the equivariant cohomology ring of the point. \par

Given a quasi-projective algebraic variety $X$ over $\bbC$ equipped with a linearisable action of an algebraic group $G$, let $\Db_G(X)$ denote the Bernstein--Lunts \cite{BL94} $G$-equivariant derived category of $\bbC$-sheaves on $X$ with bounded constructible cohomology. For $\scrM, \scrN\in \Db_{G}(X)$, we denote $\Hom^{n}_G\left( \scrM, \scrN \right) = \Hom_{\Db_G(X)}\left(\scrM, \scrN[n]  \right)$ and $\Hom^{\bullet}_G\left( \scrM, \scrN \right) = \bigoplus_{n\in \bbZ}\Hom^n_{G}\left(\scrM, \scrN \right)$, viewed as finitely generated graded $\rmH^{\bullet}_G$-module. We will also denote $\Hom^{\bullet+n}_G\left( \scrM, \scrN \right) = \Hom^{\bullet}_G\left( \scrM, \scrN[ n] \right)$ for $n\in \bbZ$. In addition, we will write $\RHom_{G}(\scrM, \scrN)$ for the derived Hom-space, viewed as object of the derived category of $\rmH^{\bullet}_{G}$-dg-modules. Given $\scrM\in \Db_G(X)$, we write $\rmH^{\bullet}_G(X, \scrM) = \Hom^{\bullet}_G(\bbC, \scrM)$ for the equivariant cohomology with coefficients in $\scrM$. \par

We will denote by $\Perv_G\left( X \right)$ the category of $G$-equivariant perverse sheaves on $X$. We normalise the perverse t-structure on $\Db_G(X)$ in such a way that the forgetful functor $\Db_G(X)\to \Db(X)$ is perverse t-exact. We will denote by $\pH^k$ the $k$-th cohomology with respect to this perverse t-structure and $\pH^{\bullet} = \bigoplus_{n\in \bbZ}\pH^{n}$. \par

We let $\bfX_*(G)$ denote the set of cocharacters of $G$. We adopt the notion of fractional cocharacters introduced in \cite{LYI}:
\[
	\bfX_*(G)_{\bbQ} = (\bfX_*(G) \times \bbZ_{>0}) / \{(\lambda, a)\sim (\mu, b)\;\vert\; q\lambda = b\mu\}.
\]
We will denote write $\lambda / a = (\lambda, a)$. Given a rational representation $G\to \GL(V)$ and $r = p/q\in \bbQ$ with $p,q\in \bbZ$ and $q\neq 0$, we will write $\prescript{\lambda/a}{r}V$ for the $r$-weight space of $\lambda/a$, defined as $\prescript{\lambda/a}{r}V = \left\{ v\in V\;\vert\; \lambda(t^q)v = t^{ap}v,\;\forall t\in \bbC^{\times} \right\}$. \par

\section{Reminder on \texorpdfstring{$\bbZ/m$}{Z/m}-graded Lie algebras}\label{sec:graduee} 
We review in this section basic notions and results related to $\bbZ/m$-graded Lie algebras developed in \cite{LYI,LYIII}. 

\subsection{Semisimple automorphisms}\label{subsec:pin}
In this article, $G$ will denote a connected reductive algebraic group unless otherwise stated. Let $\Out(G) = \Aut(G) / G^{\ad}$ be the group of outer automorphisms of $G$. An \emph{outer type} of $G$ of order $\bfe\ge 1$ is an injective group homomorphism $\sigma:\mu_\bfe\hookrightarrow \Out(G)$. Fix an outer type $\sigma$ of $G$ of order $\bfe$. Let $\mu_\bfe^*\subseteq \mu_\bfe$ denote the subset of primitive elements and $\Aut_\sigma(G)$ the pre-image of $\sigma(\mu_\bfe^*)$ under the quotient map $\Aut(G)\to \Out(G)$. We shall recollect some results about the structure of the pair $(G, \sigma)$.  \par

A pinning $E = (B, T, a)$ for $G$ is a triplet consisting of a Borel subgroup $B\subseteq G$, a maximal torus $T\subseteq B$ and a morphism of algebraic groups $a:[B,B]\to \bbG_{\rma}$ satisfying a condition of non-degeneracy; let $\Aut_E(G)\subseteq \Aut(G)$ denote the subgroup of automorphisms fixing the datum $E$. The set of all pinnings of $G$ forms a principal homogeneous $G^{\ad}$-space and every pinning $E$ yields a decomposition $\Aut(G) = G^{\ad}\rtimes \Aut_E(G)$ and thus a canonical isomorphism $\Aut_E(G)\cong \Out(G)$, under which $\sigma$ is identified with a monomorphism $\sigma_E: \mu_\bfe\to \Aut_E(G)$; the latter yields $\Aut_{\sigma}(G) = G^{\ad}\rtimes \sigma_E(\mu_\bfe^*)$.
\begin{prop}\label{prop:in-out}
	Given $G$ as above and a pinning $E = (B, T, a)$, every semisimple element of $\Aut_\sigma(G)$ is $G$-conjugate to an element of the subset
	\[
		S^{\ad}\times \sigma_E(\mu_\bfe^*), \quad \text{where $S^{\ad} = S / (S\cap Z_G)$ and $S = (T^{\sigma_E})^{\circ}$}.
	\]
\end{prop}
\begin{proof}
	The proof of \cite[2.2.2]{LYIII} goes through for arbitrary connected reductive groups, see also \cite[3.2]{reeder10}.
\end{proof}
Consequently, given any semisimple element $\tau$ of $\Aut_{\sigma}(G)$, there exists a pinning $E = (B, T, a)$ such that $\tau\in S^{\ad}\times \sigma_E(\mu_\bfe^*)$.

The following is also proven in \cite[2.2.2]{LYIII}:
\begin{prop}\label{prop:pl}
	Given a semisimple automorphism $\tau\in \Aut(G)$ whose image in $\Out(G)$ lies in $\sigma(\mu_\bfe^*)$, the subgroup $M = (G^{\tau})^{\circ}$ is reductive. Moreover, if $E = (B, T, a)$ is a pinning of $G$ such that $\tau\in S^{\ad}\times \sigma_E(\mu^*_\bfe)$, then $S$ is a maximal torus of $M$. \hfill\qedsymbol
\end{prop}

\begin{defi}
	A subgroup $M$ of $G$ of the form described in \autoref{prop:pl} is called a \emph{pseudo-Levi subgroup} of $(G, \sigma)$; the Lie algebra of $M$ is called a \emph{pseudo-Levi subalgebra} of $(\frakg, \sigma)$.
\end{defi}
\begin{rema}
	Given a pseudo-Levi subgroup $M$ of $(G, \sigma)$, we can always choose a finite-order automorphism $\tau\in S^{\ad}\times \sigma_E(\mu^*_\bfe)$ such that $M = (G^{\tau})^{\circ}$. 
\end{rema}

\subsection{Cyclic gradings on reductive groups}\label{subsec:grading}
Assume now that we are given a connected reductive group $G$, a positive integer $m\ge 1$ and a group homomorphism
\begin{eq}
	\theta : \mu_m\to \Aut(G).
\end{eq}
Then, $\theta$ defines a $\bbZ / m$-grading on the Lie algebra $\frakg = \Lie G$:
\[
	\frakg = \bigoplus_{\ubar i\in \bbZ / m} \frakg_{\ubar i}, \quad \frakg_{\ubar i} = \left\{ x\in \frakg\;\vert\; \theta(\omega)x = \omega^i,\quad \forall \omega\in \mu_m \right\}.
\]
Let $G^{\theta}$ be the fixed-point subgroup. After \cite[8.1]{steinberg68}, it is known that $G^{\theta}$ is connected whenever $G$ is semisimple and simply connected; however, $G^{\theta}$ may be disconnected in general. We define $G_{\ubar 0} = (G^{\theta})^{\circ}$ to be its neutral component.%For the sake of simplification, we make the following assumption throughout this article:
%\begin{hypo}\label{hypo:conn}
	%$G^{\theta}$ is connected.
%\end{hypo}

Let $\bfe\in \bbZ_{>0}$ be the order of the image of $\theta(\mu_m)$ in the quotient $\Out(G) = \Aut(G)/G^{\ad}$, so that the composite of $\theta$ with the quotient factorises through $m/\bfe: \mu_m\to \mu_{\bfe}$ and induces a monomorphism $\sigma : \mu_\bfe\hookrightarrow \Out(G)$. 
\begin{prop}\label{prop:graduation}
	There exists a pinning $E= (B, T, a)$ for $G$ satisfying the following property: set $S = (T^{\sigma_E})^{\circ}$ and $S^{\ad} = S / (Z_G \cap S)\subseteq G^{\ad}$; then there exists a homomorphism $\theta^{\mathrm{in}}:\mu_m\to S^{\ad}$ such that 
	\[
		\theta(\omega) = \theta^{\mathrm{in}}(\omega)\sigma_E(\omega^{m/\bfe})\quad \text{for every $\omega\in \mu_m$}.
	\]
\end{prop}
\begin{proof}
	Choose a primitive $m$-th root of unity $\omega\in \mu^*_m$ and apply \autoref{prop:in-out} to the automorphism $\theta(\omega)\in \Aut_{\sigma}(G)$. 
\end{proof}

\subsection{Graded \texorpdfstring{$\fraksl_2$}{sl2}-triplets}\label{subsec:JM}
Let $\ubar d\in \bbZ / m$. The Jacobson--Morosov theorem in the $\bbZ/m$-graded setting is proven in \cite[2.3]{LYI}: given every $e\in \frakg_{\ubar d}^{\nil}$, we can complete $e$ into a triplet $\phi = \left( e, h, f \right)$ with $h\in \frakg_{\ubar 0}$ and $f\in \frakg_{-\ubar d}$, satisfying the standard $\fraksl_2$-relations $[h,e] = 2e$, $[h,f] = -2f$ and $[e,f] = h$; moreover, such an $\fraksl_2$-triplet gives rise to a cocharacter $\iota\in \bfX_*(G_{\ubar 0})$ which satisfies $\mathrm{d}\iota(1) = h$. \par

Let $\Cq = \bbC^{\times}$ be a one-dimensional torus. We let $\Cq$ act trivially on $G$ and by weight $-2$ on the Lie algebra $\frakg$. Let $\rmO = \Ad_{\Gz}e\subseteq \frakg^{\nil}_{\ubar d}$. It follows that $\Cq e = \Ad_{\iota(\bbC^{\times})} e\subset \rmO$, so that $\rmO$ is stable under the action of $\Cq$. Moreover, the natural inclusion $\Gz\subset \Gz\times \Cq$ induces an isomorphism of component groups of stabilisers of $e$:
\begin{equation}\label{equa:stab_e}
	\pi_0(Z_{\Gz}(e))\cong \pi_0(Z_{\Gz\times \Cq}(e)). 
\end{equation}
%\begin{rema}
	%In~\cite{LYI}, one is restricted to simply connected semisimple groups. In the present setting, one can apply the results of {\it loc. cit.} to the derived subgroup $G^{\mathrm{der}}$, which is semisimple and simply connected by hypothesis, and then lift them onto $G$. The equivariant cohomology can be calculated by means of the Hochschild--Serre spectral sequence 
	%\begin{eq}
		%\rmH^p\left( \Gz /G^{\mathrm{der}}_{\ubar 0}, \Hom^{q}_{G^{\mathrm{der}}_{\ubar 0}}\left(\relbar,\relbar \right) \right)\Longrightarrow \Hom^{p+q}_{\Gz}\left( \relbar , \relbar \right). 
	%\end{eq}
%\end{rema}

\subsection{Spirals and splittings}\label{subsec:spiral}
Let $\lambda\in \bfX_*\left(\Gz\right)_{\bbQ}$ be a fractional cocharacter and $\epsilon\in \left\{ 1, -1 \right\}$ a sign. We associate to them $\bbZ$-graded Lie algebras $\sgn\frakp^{\lambda}_* =\bigoplus_{n\in \bbZ}  \sgn\frakp^\lambda_n, \frakl^{\lambda}_* =\bigoplus_{n\in \bbZ}  \frakl^\lambda_n$ and $\sgn\fraku^{\lambda}_* =\bigoplus_{n\in \bbZ}  \sgn\fraku^\lambda_n$ by setting\footnote{Our notation differs slightly from Lusztig--Yun \cite{LYI}: our $\sgn\frakp^{\lambda}_*$ equals their $\sgn\frakp^{\epsilon\lambda}_*$, etc. }
\[
	\sgn\frakp^{\lambda}_n = \prescript{\epsilon\lambda}{\ge \epsilon n}\frakg_{\ubar n} = \begin{cases}\prescript{\lambda}{\ge n}\frakg_{\ubar n} & \epsilon = 1 \\ \prescript{\lambda}{\le n}\frakg_{\ubar n}& \epsilon = -1\end{cases},\quad \frakl^{\lambda}_n =  \prescript{\lambda}{n}\frakg_{\ubar n},\quad \sgn\fraku^{\lambda}_n = \prescript{\epsilon\lambda}{>\epsilon n}\frakg_{\ubar n}.
\]

%Note that $\frakp^{\lambda}_*$ and $\fraku^{\lambda}_*$ depend on the sign $\epsilon$. Any such $\bbZ$-graded Lie algebra $\frakp_*^{ \lambda }$ is called the ($\epsilon$-)\emph{spiral} of $\frakg$ attached to $\lambda$, the Lie subalgebra $\sgn\fraku^{\lambda}_*$ is called the \emph{radical} of $\sgn\frakp^{\lambda}_*$. There is a spiral version of the Levi decomposition $\sgn\frakp^\lambda_* = \frakl^\lambda_* \oplus \sgn\fraku^\lambda_*$. If $\lambda, \lambda'\in \bfX_*(G_{\ubar 0})_{\bbQ}$ are such that $\sgn\frakp^\lambda_* = \sgn\frakp^{\lambda'}_*$, then $\sgn\fraku^\lambda_* = \sgn\fraku^{\lambda'}_{*}$ holds but $\frakl^{\lambda}_*$ and $\frakl^{\lambda'}_*$ can differ. Any such subalgebra $\frakl^{\lambda}_*\subset \sgn\frakp^{\lambda}_*$ is called a \emph{splitting} of $\sgn\frakp^\lambda_*$. \par

Let $r\in \bbZ_{>0}$ be such that $\lambda_0 = r\lambda\in \bfX_*(\Gz)$. Pick any primitive $(rm)$-th root of unity $\omega\in \mu^*_{rm}$ and consider $\tau = \theta(\omega^r)\Ad_{\lambda_0(\omega)}^{-1}\in \Aut(G)$. Then, the connected fixed-point subgroup $L^{\lambda} := (G^{\tau})^{\circ}$ is a pseudo-Levi subgroup of $(G, \sigma)$ in the sense of \autoref{subsec:grading} and it can be shown as in \cite[2.6(a)]{LYI} that $\Lie L^{\lambda} = \frakl^{\lambda}_*$. Let $L^{\lambda}_0 := Z_{L^{\lambda}}(\lambda)$ be the centraliser of $\lambda$, so that $L^{\lambda}_0\subseteq \Gz$ and $\Lie L^{\lambda}_0 = \frakl^{\lambda}_0$ hold. We set $\sgn U^{\lambda}_0 = \exp(\sgn \fraku^{\lambda}_0)$ and $\sgn P^{\lambda}_0 = L^{\lambda}_0\sgn U^{\lambda}_0$. 
\begin{defi}
	The $\bbZ$-graded Lie algebra $\sgn\frakp^{\lambda}_*$ is called the \emph{($\epsilon$-)spiral} attached to $\lambda$. The datum $\frakl^{\lambda}_*$ is called a \emph{splitting factor} of $\sgn\frakp^{\lambda}_*$ and $\sgn\fraku^{\lambda}_*$ is called the \emph{nil-radical} of $\sgn \frakp^{\lambda}_*$. 
\end{defi}
Let $\frakP^{\epsilon}$ denote the set of $\epsilon$-spirals and let $\frakG$ denote the set of splitting factors of spirals. A subscript $G$ will sometimes be added to indicate the dependence on $G$ when other groups with a $\bbZ/m$-grading are involved. It is clear that every spiral $\frakp_*\in \frakP^{\epsilon}$ determines the datum $(P_0, \frakp_*, U_0, \fraku_*)$ and every splitting factor $\frakl_*\in \frakG$ determines the datum $(L, L_0, \frakl, \frakl_*)$ as above. \par % In particular, there is an surjection $\frakP^{\epsilon}\to \frakG$ which sends a spiral system $(P_0, U_0, L, L_0, \frakp_*, \fraku_*, \frakl_*)$ to its splitting components $(L, L_0, \frakl_*)$.
It is shown in \cite[2.7(a)]{LYI} that for every spiral $\frakp_*\in \frakP^\epsilon$, the set of splitting factors of $\frakp_*$ is a principal homogeneous $U_0$-space under the natural action $u\cdot \frakl_* = \Ad_u\frakl_*$. \par

\subsection{Spiral attached to a nilpotent element}\label{subsec:spiral-sl2}
Let $d\in \bbZ \setminus \left\{ 0 \right\}$ and let $z\in \frakg^{\nil}_{\ubar d}$ be a nilpotent element. According to \autoref{subsec:JM}, we can complete $e$ into a graded $\fraksl_2$-triplet $\phi =\left( e, h, f \right)$ and it gives rise to $\lambda\in \bfX_*(\Gz)$. It gives rise to a $\epsilon$-spiral together with a splitting factor
\[
	\frakp_* = \sgn\frakp^{(d/2)\lambda}_*\in \frakP^{\epsilon},\quad \frakl_* = \frakl^{(d/2)\lambda}_*\in \frakG
\]
with $\epsilon = d / |d|$. The spiral $\frakp_*$ does not depend on the choice of $\phi$. \par
The following results are proven in~\cite[2.9]{LYI}:
\begin{lemm}\label{lemm:spiral-stab}
	The following statements hold:
	\begin{enumerate}
		\item
			The orbit $\Ad_{P_0}e$ is open in $\frakp_{d}$. 
		\item %(L,L_0, \frakl_*)$ n'est pas forcément un pseudo-Levi
			$Z_{L_0}(e)$ is a maximal reductive subgroup of $Z_{\Gz}(e)$. In particular, the inclusion $Z_{L_0}(e)\subseteq Z_{\Gz}(e)$ induces an isomorphism $\pi_0\left( Z_{L_0}(e) \right)\cong \pi_0\left( Z_{\Gz}(e) \right)$.
	\hfill \qed
	\end{enumerate}
\end{lemm}

\subsection{\texorpdfstring{$\theta$}{θ}-isotropic Levi subgroups}\label{subsec:isotrope}

\begin{defi}
 A $\theta$-stable Levi subgroup of $G$ is called \emph{$\theta$-isotropic} if it is a Levi factor of a $\theta$-stable parabolic subgroup of $G$. 
\end{defi}
\begin{prop}
	Given a $\theta$-stable parabolic subgroup $Q$ of $G$, there exists a cocharacter $\lambda\in \bfX_*(\Gz)$ such that the Lie algebra $\frakq = \Lie Q$ is given by $\frakq = {^\lambda_{\ge 0}}\frakg$. 
	%Every $\theta$-stable parabolic subgroup of $G$ admits a $\theta$-isotropic Levi factor. %Conversely, every parabolic subgroup of $G$ admitting a $\theta$-isotropic Levi factor is $\theta$-stable. %c'est faux!
\end{prop}
\begin{proof}
	By Steinberg's theorem \cite[7.5]{steinberg68}, $Q$ admits a $\theta$-stable Borel pair $(B, T)$, which is also a Borel pair for $G$ because $Q$ is a parabolic subgroup. We may pick $\lambda'\in \bfX_*(T)$ such that $\frakq = {^{\lambda'}_{\ge 0}}\frakg$. Then, since $T$ and $Q$ are $\theta$-stable, the cocharacter $\theta(\omega)\circ\lambda'\in \bfX_*(T)$ for $\omega\in \mu_m$ also satisfies $\frakq = {^{\theta(\omega)\circ\lambda'}_{\ge 0}}\frakg$. Therefore, 
	\[
		\lambda = \sum_{\omega\in \mu_m}\theta(\omega)\circ\lambda'\in \bfX_*(T)^{\theta} \subseteq \bfX_*(\Gz)
	\]
	satisfies the requirement.
\end{proof}
\begin{coro}\label{coro:theta-parabolic}
	Every $\theta$-stable parabolic subgroup of $G$ admits a $\theta$-stable Levi factor. \hfill\qedsymbol
\end{coro}

\begin{coro}\label{coro:theta-isotropic}
	Every $\theta$-isotropic Levi subgroup of $G$ is of the form $G^{\lambda} = Z_G(\lambda)$ for some cocharacter $\lambda\in \bfX_*(\Gz)$. Conversely, given a subtorus $A\subseteq G_{\ubar 0}$, the centraliser $Z_G(A)$ is a $\theta$-isotropic Levi subgroup of $G$. \hfill\qedsymbol
\end{coro}

%
%\begin{prop}\label{prop:levi-sc}
%Let $M$ be a $\theta$-isotropic Levi subgroup of $G$. Then, the derived subgroup $M^{\der} = [M, M]$ is simply connected and its $\theta$-fixed-point subgroup $M^{\der}_{\ubar 0} = (M^{\der})^{\theta}$ is connected.
%\end{prop}
%\begin{proof}
	%Let $T$ be any maximal torus of $M$. We may choose a basis $\Delta$ for $R(G, T)$ such that $M$ is the centraliser of a subspace $\fraka\subseteq \Lie T$ of the form $\fraka = \bigcap_{\alpha\in \Delta\setminus\Delta'}\ker(\rmd \alpha)$ for some $\Delta'\subseteq \Delta$. Then, the simple connectivity of $G$ implies that $\bfX_*(T)$ is generated by the simple coroots $\left\{ \alpha^{\vee} \right\}_{\alpha\in \Delta}$. Since the cocharacter lattice of the maximal torus $T'\subseteq M^{\der}$ is a subgroup of $\bfX_*(T)$ of rank $\dim \frakg - \dim \fraka = \#\Delta'$, it must be generated by $\left\{ \alpha^{\vee} \right\}_{\alpha\in \Delta'}$. We observe that $\Delta'$ is a basis for $R(M^{\der}, T')$ and $\left\{ \alpha^{\vee} \right\}_{\alpha\in \Delta'}$ is its dual basis. Therefore, $M^{\der}$ is simply connected. The connectivity of $M^{\der}_{\ubar 0}$ follows from Steinberg's theorem.
%\end{proof}

\section{Sheaves on \texorpdfstring{$\bbZ/m$}{Z/m}-graded Lie algebras}\label{sec:faisc}
We review the spiral induction and restriction introduced in \cite{LYI} and explain the relation between spiral restriction and Braden's theorem of hyperbolic restriction \cite{braden03}. We then introduce the $\bbZ/m$-graded version of parabolic induction and restriction.

\subsection{Setup}\label{subsec:G}
In this section, we fix a connected reductive algebraic group $G$, a map $\theta:\mu_m\to \Aut(G)$ with outer type $\sigma:\mu_{\bfe}\to \Out(G)$, a pinning $E = (B, T, a)$ for $G$ provided by \autoref{prop:graduation}. In addition, we fix $d\in \bbZ \setminus \left\{ 0 \right\}$ and write $\ubar d = d\mod m\in \bbZ/m$. \par
We let $\Pi$ denote the set of pairs $\left( \rmC, \scrL \right)$ where $\rmC\subset \frakg^{\nil}_{\ubar d}$ is a nilpotent $\Gz$-orbit and $\scrL\in \Loc_{G_{\ubar 0}}(\rmC)$ is an irreducible $\Gz$-equivariant local system on $\rmC$, defined up to isomorphism. Given a $\Gz$-orbit $\rmC\subseteq \frakg^{\nil}_{\ubar d}$, irreducible local systems are equivalent to irreducible representations of the equivariant fundamental group $\pi_1^{\Gz}(\rmC, z)$ for any $z\in \rmC$. The latter term is isomorphic to the component group $\pi_0(Z_{\Gz}(z))$ of the stabiliser of $z$, which is a finite group. Moreover, the set of nilpotent orbits $\frakg^{\nil}_{\ubar d} / \Gz$ is finite, so $\Pi$ is a finite set. Given a pair $\pi = (\rmC, \scrL)$, we shall write $\IC_{\pi} = \IC(\scrL)$ for the intersection complex on $\frakg^{\nil}_{\ubar d}$ with coefficients in $\scrL$. From the theory of perverse sheaves, it is known that $\left\{ \IC_{\pi}\right\}_{\pi\in \Pi}$ forms a complete list of representatives for the isomorphism classes of simple objects of $\Perv_{\Gz}(\frakg^{\nil}_{\ubar d})$. The isomorphism \eqref{equa:stab_e} from \autoref{subsec:JM} implies that every element of $\Pi$ admits a unique $(\Gz\times \Cq)$-equivariant enhancement. \par

\subsection{Spiral induction and restriction}\label{subsec:spiral-ind}
Suppose we are given $\epsilon \in \{ \pm 1\}$, an $\epsilon$-spiral $\frakp_*\in \frakP^{\epsilon}$ and a splitting factor $\frakl_*\in\frakG$ as in \autoref{subsec:spiral}. They determine a datum $(P_0, \frakp_*, L, L_0, \frakl, \frakl_*)$.
Fixing $d\in \bbZ \setminus\left\{ 0 \right\}$, we have the following diagrams:
\begin{eq*}\label{equa:ind-res-spirale}
	\frakl_{d}\xleftarrow{\alpha} \Gz\times^{U_{0}}\frakp_{d}\xrightarrow{\beta} \Gz\times^{P_{0}}\frakp_{d}\xrightarrow{\gamma} \frakg_{\ubar d}, \qquad
	\frakg_{\ubar d}\xleftarrow{i} \frakp_{d}\xrightarrow{p} \frakl_{d}.
\end{eq*}
The \emph{$\epsilon$-spiral induction, restriction} and \emph{corestriction} are defined as\footnote{The cohomological shift is added to make the adjunctions cleaner.}
\begin{eq}
	\Ind^{\frakg_{\ubar d}}_{\frakp_{d}} &= \gamma_*(\beta^*)^{-1}\alpha^*[\dim \fraku_0 + \dim \fraku_\eta]: \Db_{L_{0}}\left( \frakl_{d} \right)\to \Db_{\Gz}\left( \frakg_{\ubar d} \right) \\
	\Res^{\frakg_{\ubar d}}_{\frakp_{d}} &= p_!i^*[\dim \fraku_\eta - \dim \fraku_{0}]: \Db_{\Gz}\left( \frakg_{\ubar d} \right)\to \Db_{L_{0}}\left( \frakl_{d} \right) \\
	\cores^{\frakg_{\ubar d}}_{\frakp_{d}} &= p_*i^![\dim \fraku_0 - \dim \fraku_{\eta}]: \Db_{\Gz}\left( \frakg_{\ubar d} \right)\to \Db_{L_{0}}\left( \frakl_{d} \right).
\end{eq}

There are adjunctions $\Res^{\frakg_{\ubar d}}_{\frakp_{d}} \dashv \Ind^{\frakg_{\ubar d}}_{\frakp_{d}}$ and $\Ind^{\frakg_{\ubar d}}_{\frakp_{d}}\dashv \cores^{\frakg_{\ubar d}}_{\frakp_{d}}$. Moreover, when $\epsilon = d / |d|$, the $\epsilon$-spiral induction produces complexes on the nilpotent cone:
\[
	\Ind^{\frakg_{\ubar d}}_{\frakp_{d}}:\Db_{L_{0}}\left( \frakl_{d} \right)\to \Db_{\Gz}\left( \frakg^{\nil}_{\ubar d} \right).
\]

\subsection{LY-cuspidal supports and decomposition}\label{subsec:dco}
Let $\frakl_*\in \frakG$ be a splitting factor of spiral, $\rmO\subseteq \frakl_d$ an $L_0$-adjoint orbit and $\scrC$ an irreducible $L_0$-equivariant local system on $\rmO$. The pair $(\rmO, \scrC)$ on $\frakl_d$ is called \emph{cuspidal} if it is of the form $(\tilde\rmO\cap \frakl_{d}, \tilde\scrC\mid_{\tilde\rmO\cap \frakl_{d}})$ for some cuspidal pair $(\tilde\rmO, \tilde\scrC)$ on $\frakl$ in the sense of \cite{lusztig84} and satisfies $\til\rmO\cap \frakl_d \neq \emptyset$. The following properties of cuspidal pairs are proven in \cite[4.4]{lusztig95b}: 
\begin{prop}\label{prop:cuspidal-lusztig}
		Let $(\rmO, \scrC)$ be a cuspidal pair on $\frakl_d$.  Then, the following statements hold:
		\begin{enumerate}
			\item
				The orbit $\rmO$ is distinguished: for each $z\in \rmO$, the group $(Z_{L_0}(z) / Z_L)^{\circ}$ is unipotent.
			\item
				Given an $\fraksl_2$-triplet $(e, h, f)$ with $e\in \rmO$, $h\in \frakl_0$ and $f\in \frakl_{-d}$ and the associated cocharacter $\iota\in \bfX_*(L_0)$ characterised by $\rmd\iota(1) = h$, we have 
				\[
					\frakl_n = {^\iota_{2n/d}}\frakl,\quad \text{for $n\in \bbZ$}. 
				\]
			\item
				$\rmO$ coincides with the unique open $L_0$-orbit $\ci\frakl_d\subseteq \frakl_d$
			\item
				The cuspidal pair $(\rmO, \scrC)$ is clean. Namely, the natural morphism $j_{\rmO!}\scrC\to j_{\rmO*}\scrC$ is an isomorphism, where $j_{\rmO}:\rmO\hookrightarrow \frakl_d$ denotes the inclusion.
		\end{enumerate}
\end{prop}
\begin{defi}
	An \emph{LY-cuspidal support} on $\frakg_{\ubar d}$ is a datum $\left(L,\frakl_*, \scrC  \right)$, where $\frakl_*\in \frakG$, $L$ is the pseudo-Levi subgroup of $(G, \sigma)$ whose Lie algebra is $\frakl_*$ and $\scrC$ is a cuspidal local system on the open $L_0$-orbit $\ci\frakl_d\subseteq \frakl_d$. 
\end{defi}

An isomorphism of LY-cuspidal supports $(L, \frakl_*, \scrC)\cong (L', \frakl'_*, \scrC')$ is a pair $(g, \eta)$ of element $g\in \Gz$ such that $gLg^{-1} = L'$, $\Ad_g \frakl_n = \frakl'_n$ and an isomorphism $\eta : \Ad_g^* \scrC'\cong \scrC$.  Let $\frakT = \frakT_G$ denote the groupoid of LY-cuspidal supports on $\frakg_{\ubar d}$ and let $\ubar\frakT$ denote the set of isomorphism classes of its objects. \par

We define $\Perv_{\Gz}( \frakg_{\ubar d}^{\nil})_{\xi}$ to be the Serre subcategory of $\Perv_{\Gz}( \frakg_{\ubar d}^{\nil} )$ generated by the constituents of $\Ind^{\frakg_{\ubar d}}_{\frakp_d}\IC( \scrC )$ for various $\epsilon$-spirals $\frakp_*$ with $\epsilon = d / |d|$ which has $\frakl_*\in \frakG$ as splitting factor. We define $\Db_{\Gz}( \frakg_{\ubar d}^{\nil} )_{\xi}$ to be the full subcategory of $\Db_{\Gz}( \frakg_{\ubar d}^{\nil} )$ of objects $\scrK$ such that $\prescript{p}{}\scrH^k\scrK$ lies in $\Perv_{\Gz}( \frakg_{\ubar d}^{\nil} )_{\xi}$ for all $k\in \bbZ$. We let $\Pi_\xi\subset \Pi$ denote the subsets of pairs $\pi$ such that $\IC_\pi\in \Perv_{\Gz}( \frakg^{\nil}_{\ubar d} )_\xi$. These subcategories and subset are called \emph{$\xi$-block}. \par

The main result \cite[0.6]{LYI} provides decompositions into blocks:
\begin{eq*}\label{eq:dco}
	\Pi = \bigsqcup_{\xi\in \ubar\frakT} \Pi_{\xi},\quad \Perv_{\Gz}\left( \frakg^{\nil}_{\ubar d} \right) = \bigoplus_{\xi\in \ubar\frakT} \Perv_{\Gz}\left( \frakg^{\nil}_{\ubar d} \right)_{\xi}, \qquad \Db_{\Gz}\left( \frakg^{\nil}_{\ubar d} \right) = \bigoplus_{\xi\in \ubar\frakT} \Db_{\Gz}\left( \frakg^{\nil}_{\ubar d} \right)_{\xi}.
\end{eq*}
Note that even though the above decompositions have been proven under the assumption that $G$ is semisimple and simply connected, the proof works {\itshape verbatim} in our setting. 
\begin{rema}\label{rema:Z-graded}
	Similar decompositions exist for $\bbZ$-graded reductive Lie algebras, see \cite{lusztig95b} and \cite[\S 1]{LYI}. Henceforth, given a connected reductive algebraic group $L$ with a $\bbZ$-grading $\frakl_*$ on the Lie algebra $\frakl = \Lie L$, we will let $\Pi_L$, $\frakT_L$ and $\Pi_{L, \xi}$ denote the counterparts of $\Pi$, $\frakT$ and $\Pi_{\xi}$ for $(L, \frakl_*)$ in the $\bbZ$-graded setting. 
\end{rema}

\subsection{Extension of pairs}%\hl{Revoir la necessite}

\iffalse
The following is an immediate consequence of the main theorems of~\cite{LYI}.
\begin{lemm}\label{lemm:ind-bloc}
	Let $\frakp_* = \frakp^{\lambda}_*$ and $\frakl_* = \frakl^{\lambda}_*$ be spiral and splitting factor as above. Let $\xi = (M, \frakm_*, \scrC)\in \frakT$ be an LY-cuspidal support on $\frakl_d$. Then the spiral induction and restriction preserve blocks:
	\begin{eq}
		\Ind^{\frakg_{\ubar d}}_{\frakp_d}: \Db_{L_0}\left( \frakl_d \right)_{\xi}\to \Db_{\Gz}\left( \frakg^{\nil}_{\ubar d} \right)_{\xi} \\
		\Res^{\frakg_{\ubar d}}_{\frakp_d}, \cores^{\frakg_{\ubar d}}_{\frakp_d}: \Db_{\Gz}\left( \frakg^{\nil}_{\ubar d} \right)_{\xi}\to \Db_{L_0}\left( \frakl_d \right)_{\xi}.
	\end{eq}
	On the other hand, if $\xi' = \left( M', M'_0, \frakm'_*, \scrC' \right)$ is an LY-cuspidal support on $\frakg_{\ubar d}$ such that no $\Gz$-conjugate of $(M', \frakm'_*)$ is a $\bbZ$-graded Levi of $(L, \frakl_*)$, then $\Res^{\frakg_{\ubar d}}_{\frakp_d}$ annihilates $\Db_{\Gz}\left( \frakg^{\nil}_{\ubar d} \right)_{\xi'}$.  \hfill\qed
\end{lemm}

\fi

\begin{lemm}\label{lemm:ext-block}
	Let $(\rmO,\scrL)\in \Pi_{\xi}$ for $\xi\in \frakT$. Then $j_{\rmO !}\scrL$ and $j_{\rmO *}\scrL$ lie in $\Db_{\Gz}\left( \frakg^{\nil}_{\ubar d} \right)_{\xi}$. 
\end{lemm}
\begin{proof}
	The decomposition of equivariant derived category yields a decomposition of complex
	\begin{eq}
		j_{\rmO *}\scrL \cong \bigoplus_{\xi'\in \ubar\frakT}\left(j_{\rmO *}\scrL\right)_{\xi'},\quad \left(j_{\rmO *}\scrL\right)_{\xi'}\in \Db_{\Gz}(\frakg^{\nil}_{\ubar d})_{\xi'}.
	\end{eq}
	Let $\xi' = \left( M', \frakm'_*, \scrC' \right)\in \frakT$ be any LY-cuspidal support on $\frakg_{\ubar d}$ such that $\xi'\not\cong \xi$, we show that $\left(j_{\rmO *}\scrC\right)_{\xi'} = 0$. Suppose that $\left(j_{\rmO *}\scrL\right)_{\xi'} \neq 0$. Let $\delta\in \bbZ$ be the smallest integer such that $\prescript{p}{}\scrH^\delta\left(j_{\rmO *}\scrL\right)_{\xi'} \neq 0$. Let $\scrF\in \Perv_{\Gz}\left( \frakg_{\ubar d}^{\nil} \right)_{\xi'}$ be a simple subobject of $\prescript{p}{}\scrH^\delta\left(j_{\rmO *}\scrL\right)_{\xi'}\neq 0$ so that 
	\[
		\Hom\left(\scrF , \prescript{p}{}\tau^{\le \delta}\left( j_{\rmO *}\scrL\right)_{\xi'}[\delta] \right)\neq 0
	\]
	holds. The distinguished triangle of perverse truncation
	\begin{eq}
		\prescript{p}{}\scrH^\delta\left(j_{\rmO *}\scrL\right)_{\xi'} \to \left(j_{\rmO *}\scrL\right)_{\xi'} \to \prescript{p}{}\tau^{> \delta}\left(j_{\rmO *}\scrL\right)_{\xi'} \xrightarrow{[1]}
	\end{eq}
	yields an injective map
	\begin{eq}
		\Hom\left(\scrF , \prescript{p}{}\scrH^\delta\left(j_{\rmO *}\scrL\right)_{\xi'}[\delta] \right)\to \Hom\left(\scrF, \left(j_{\rm O *}\scrL\right)_{\xi'}[\delta]  \right) = \Hom\left(j_{\rmO}^*\scrF , \scrL[\delta]  \right),
	\end{eq}
	so the last term is non-zero, which contradicts the fact that $\xi\not\cong\xi'$ by~\cite[13.8(a)]{LYII}. Thus $j_{\rmO*}\scrL = \left(j_{\rmO*}\scrL\right)_{\xi}\in \Db_{\Gz}\left( \frakg_{\ubar d}^{\nil} \right)_\xi$. The statement about $j_{\rmO !}\scrL$ can be proven in a similar way.
\end{proof}

\begin{lemm}\label{lemm:ext-block'}
	Let $\Sigma\subseteq \frakg^{\nil}_{\ubar d}$ be a $\Gz$-stable locally closed subset and let $j: \Sigma\hookrightarrow \frakg^{\nil}_{\ubar d}$ be the inclusion. Then the functors $j_*j^*$, $j_!j^*$, $j_*j^!$ and $j_!j^!$ preserve the subcategory $\Db_{\Gz}\left( \frakg^{\nil}_{\ubar d} \right)_{\xi}$ for each $\xi\in \frakT$. 
\end{lemm}
\begin{proof}
	This follows easily from \autoref{lemm:ext-block} and \cite[13.8(a)]{LYII} by d\'evissage.
\end{proof}

\iffalse
\subsection{Odd vanishing}

\begin{prop}\label{prop:parite}
	Let $(\rmO, \scrC)$ be a cuspidal pair on $\frakl_d$, let $\frakp_*$ be a spiral on $\frakg$ which contains $\frakl_*$ as splitting factor and let $j_{\rmO'}:\rmO'\hookrightarrow \frakg_{\ubar d}$ be a $\Gz$-orbit. Then the complexes $j_{\rmO'}^*\Ind^{\frakg_{\ubar d}}_{\frakl_d\subset \frakp_d}\IC\left( \scrC \right)$ and $j_{\rmO'}^!\Ind^{\frakg_{\ubar d}}_{\frakl_d\subset \frakp_d}\IC\left( \scrC \right)$ have no cohomology in odd degrees. 
\end{prop}
\begin{proof}
	By Verdier duality and the fact that the dual $\scrC^{\vee}$ is again cuspidal, it suffices to prove the statement for $j_{\rmO'}^*$ only. The statement for $j_{\rmO'}^*$ follows from~\cite[14.10]{LYII}.
\end{proof}
\fi

\subsection{Hyperbolic restriction}\label{subsec:hyperbolic}
Let $\epsilon\in \left\{ \pm 1 \right\}$ and let $(\frakp_*, \frakl_*, \frakp'_*)$ be a triplet formed by an $\epsilon$-spiral $\frakp_*$ a $(-\epsilon)$-spiral $\frakp'_*$ and a common splitting factor $\frakl_*$ of $\frakp_*$ and $\frakp'_*$ such that $\frakp_*\cap \frakp'_* = \frakl_*$. 
Let $\Db_{G_{\ubar 0}}(\frakg_{\ubar d})^{\Cq\operatorname{-mon}}$ be the thick triangulated subcategory of $\Db_{G_{\ubar 0}}(\frakg_{\ubar d})$ generated by the image of the forgetful functor $\Db_{G_{\ubar 0}\times \Cq}(\frakg_{\ubar d})\to \Db_{G_{\ubar 0}}(\frakg_{\ubar d})$.
\begin{prop}\label{prop:hyperbolic-spiral}
There is an isomorphism of functors from $\Db_{G_{\ubar 0}}(\frakg_{\ubar d})^{\Cq\operatorname{-mon}}$ to $\Db_{L_0}(\frakl_d)$:
\[
	\Res^{\frakg_{\ubar d}}_{\frakp_d} \cong \cores^{\frakg_{\ubar d}}_{\frakp'_d}[n],
\]
where $n = \dim \fraku_d + \dim \fraku'_d - \dim \fraku_0 - \dim \fraku'_0$.
\end{prop}
\begin{proof}
	In the notation of \autoref{subsec:spiral}, assume that $\frakp_* = {^\epsilon}\frakp^{\lambda}_*$ and $\frakp'_* = {^{-\epsilon}}\frakp^{\lambda}_*$ for $\lambda\in \bfX_*(\Gz)_{\bbQ}$ and let $k\in \bbZ_{>0}$ be such that $k\lambda\in \bfX_*(\Gz)$. Set $\rho = (2k\lambda, d):\bbC^{\times}\to \Gz\times \Cq$, which yields a $\bbC^{\times}$ representation in $\frakg_{\ubar d}$. Let $\bbC^{\times}_{\rho} = \bbC^{\times}$ act via $\rho$ on $\frakg_{\ubar d}$ and let $\bbC_{\pm}$ be the one-dimensional representation of $\bbC^{\times}_{\rho}$ with weight $\pm 1$. Using \cite[1.4.7]{DG14}, we can identify $\frakp_d$ as the $\rho$-attractor in $\frakg_{\ubar d}$ and $\frakp'_{d}$ as the $\rho$-repeller
	\[
		\frakp_{d} = \Hom(\bbC_+, \frakg_{\ubar d})^{\bbC^{\times}_{\rho}},\quad \frakp'_{d}= \Hom(\bbC_-, \frakg_{\ubar d})^{\bbC^{\times}_{\rho}}.
	\]
	Consider the diagram:
	\[
		\begin{tikzcd}
			\frakg_{\ubar d} & \frakp_d \arrow{l}{i}\arrow{d}{p}\\
			\frakp'_d \arrow{u}{i'}\arrow{r}{p'} & \frakl_d.
		\end{tikzcd}
	\]
	 The theorem of Braden \cite{braden03} (see also \cite{DG14}) applies to the action of $\rho$ on $\frakg_{\ubar d}$ and yields an isomorphism $p_!i^*\cong p'_*i'^!$, whence the statement. 
\end{proof}
By \autoref{subsec:JM}, every simple object of $\Perv_{G_{\ubar 0}}(\frakg_{\ubar d}^{\nil})$ has a unique $(G_{\ubar 0}\times \Cq)$-equivariant structure. In particular, they lie in $\Db_{G_{\ubar 0}}(\frakg_{\ubar d})^{\Cq\operatorname{-mon}}$. Since they generate $\Db_{G_{\ubar 0}}(\frakg_{\ubar d}^{\nil})$ as thick subcategory, the following proposition holds:

\begin{prop}\label{prop:orbital-monodromic}
There is an inclusion $\Db_{G_{\ubar 0}}(\frakg_{\ubar d}^{\nil})\subseteq \Db_{G_{\ubar 0}}(\frakg_{\ubar d})^{\Cq\operatorname{-mon}}$. \hfill\qed
\end{prop}
In particular, \autoref{prop:hyperbolic-spiral} is applicable to the category $\Db_{G_{\ubar 0}}(\frakg_{\ubar d}^{\nil})$. \par

\subsection{Parabolic induction and restriction}\label{subsec:parab-ind}
Let $M\subseteq G$ be a $\theta$-isotropic Levi subgroup and let $Q\subseteq G$ be a parabolic subgroup having $M$ as Levi factor so that $Q = MV$ is a Levi decomposition. The Lie subalgebras $\frakm = \Lie M$ and $\frakq = \Lie Q$ and $\frakv = \Lie V$ has a $\bbZ/m$-graded restricted from $\frakg$. We have the following sequences
\begin{eq}
	\frakm_{\ubar d}\xleftarrow{\alpha} \Gz\times^{U_{\ubar 0}}\frakq_{\ubar d}\xrightarrow{\beta} \Gz\times^{Q_{\ubar 0}}\frakq_{\ubar d}\xrightarrow{\gamma} \frakg_{\ubar d}, \qquad
	\frakg_{\ubar d}\xleftarrow{\delta} \frakq_{\ubar d}\xrightarrow{\varepsilon} \frakm_{\ubar d}.
\end{eq}
The parabolic induction, restriction and corestriction are defined as
\begin{eq}
	\Ind^{\frakg_{\ubar d}}_{\frakq_{\ubar d}} &= \gamma_*(\beta^*)^{-1}\alpha^*[\dim \frakv_{\ubar 0} + \dim \frakv_{\ubar d}]: \Db_{M_{\ubar 0}}\left( \frakm_{\ubar d} \right)\to \Db_{\Gz}\left( \frakg_{\ubar d} \right) \\
	\Res^{\frakg_{\ubar d}}_{\frakq_{\ubar d}} &= \varepsilon_!\delta^*[\dim \frakv_{\ubar d} - \dim \frakv_{\ubar 0}]: \Db_{\Gz}\left( \frakg_{\ubar d} \right)\to \Db_{M_{\ubar 0}}\left( \frakm_{\ubar d} \right) \\
	\cores^{\frakg_{\ubar d}}_{\frakq_{\ubar d}} &= \varepsilon_*\delta^![\dim \frakv_{\ubar 0} - \dim \frakv_{\ubar d}]: \Db_{\Gz}\left( \frakg_{\ubar d} \right)\to \Db_{M_{\ubar 0}}\left( \frakm_{\ubar d} \right) \\
\end{eq}

As in~\autoref{subsec:spiral-ind}, there are adjunctions $\Res^{\frakg_{\ubar d}}_{\frakq_{\ubar d}} \dashv   \Ind^{\frakg_{\ubar d}}_{\frakq_{\ubar d}}$ and $\Ind^{\frakg_{\ubar d}}_{\frakq_{\ubar d}} \dashv   \cores^{\frakg_{\ubar d}}_{\frakq_{\ubar d}}$. Moreover, these functors preserve complexes on the nilpotent cones and induce $\Ind^{\frakg_{\ubar d}}_{\frakq_{\ubar d}}:\Db_{M_{\ubar 0}}\left( \frakm^{\nil}_{\ubar d} \right)\to \Db_{\Gz}\left( \frakg^{\nil}_{\ubar d} \right)$ etc. \par

%\subsection{Transitivity of induction}\label{subsec:trans}
We state two properties of transitivity for the parabolic induction functor. In addition to the datum $(Q, V, M, \frakq, \frakv,\frakm)$ as above, let $\frakq'\subseteq \frakm$ be a $\theta$-stable parabolic subalgebra and let $\frakm'\subseteq \frakq'$ be a $\theta$-stable Levi factor. Let $\til\frakq'$ be the inverse image of $\frakq'$ under the quotient map $\frakq\to \frakm$. Then we have an isomorphism of functors
\begin{equation}\label{equa:trans-parab-parab}
	\Ind^{\frakg_{\ubar d}}_{\til\frakq'_{\ubar d}} \cong \Ind^{\frakg_{\ubar d}}_{\frakq_{\ubar d}}\circ \Ind^{\frakm_{\ubar d}}_{\frakq'_{\ubar d}}.
\end{equation}
Let $\epsilon\in \left\{ \pm 1 \right\}$ and let $\frakp^M_*\in \frakP^{\epsilon}_M$ be an $\epsilon$-spiral on $\frakm_*$ with splitting factor $\frakl_*$. For any $\epsilon$-spiral $\frakp_*\in \frakP^{\epsilon}$ on $\frakg_*$ with splitting factor $\frakl_*$ such that $\frakp_i = \frakp^M_i \oplus \frakv_{\ubar i}$ holds for $i\in \left\{ 0, d \right\}$ (such $\epsilon$-spiral system always exists), there is an isomorphism of functors
\begin{equation}\label{equa:trans-parab-spiral}
	\Ind^{\frakg_{\ubar d}}_{\til\frakp_{d}} \cong \Ind^{\frakg_{\ubar d}}_{\frakq_{\ubar d}}\circ \Ind^{\frakm_{\ubar d}}_{\frakp_{d}}.
\end{equation}
The proof of the above isomorphisms are similar to \cite[4.2]{lusztig85}.

\section{Reminder on sheaf-theoretic construction of dDAHA}\label{sec:relaff}
We recollect some basic facts from \cite{LYIII} about the relative affine root systems attached to an LY-cuspidal support as well as the sheaf-theoretic construction of degenerate double affine Hecke algebras from \cite{liu21}, from which we derive the indecomposability of the $\xi$-blocks $\Db_{\Gz}(\frakg^{\nil}_{\ubar d})_{\xi}$. 
\subsection{Twisted affine root systems}

Let $G$, $\sigma$, $E = (B,T, a)$ and $\theta$ be as in~\autoref{subsec:pin}. The action of $\sigma_E:\mu_{\bfe}\to \Aut_E(G)$ yields a $\bbZ/\bfe$-grading on $\frakg$:
\[
	\frakg = \bigoplus_{i\in \bbZ / \bfe}\frakg(i),\quad \frakg(i) = \left\{ z\in \frakg\;\vert\; \sigma_E(\omega) z = \omega^iz,\;\forall \omega\in \mu_\bfe\right\}.
\]
Since $\sigma_E$ centralises the torus $S$, we have an $S$-weight-space decomposition $\frakg(i) = \bigoplus_{\mu\in \bfX^*(S)}\frakg(i)_{\mu}$ for each $i\in \bbZ / \bfe$. 
Define the $\bbZ/\bfe$-graded root system relative to $S = (T^{\sigma})^{\circ}$:
\[
	\tilde R = \left\{ (\alpha, -i)\in \bfX^*(S) \times \bbZ/\bfe\;\vert\; \alpha\neq 0,\; \frakg(i)_{\alpha}\neq 0 \right\}.
\]
The set of twisted affine roots attached to $(G, \sigma, S)$ is defined to be
\[
	R^{\aff} = \left\{ \alpha + (k/\bfe)\delta\in \bfX^*(S)\times \bbQ \delta\;\vert\; (\alpha,k\!\mod \bfe)\in \tilde R \right\}.
\]
When $G$ is almost simple and simply connected, it is known (see~\cite[Prop 18.8]{carter05}) that $\frakg(i)$ is a highest-weight $G^{\sigma}$-representation for each $i\in \bbZ / \bfe$; let $\vartheta$ denote the highest weight; in this case, let $R^+(G^{\sigma}, B^{\sigma}, S)\subset R(G^{\sigma},S)$ be the set of positive roots of the root system of $(G^{\sigma}, S)$ corresponding to the Borel subgroup $B^{\sigma}$ and let $\Delta\subset R^+(G^{\sigma}, B^{\sigma}, S)$ be the basis. We may define an affine basis $\Delta^{\aff}\subset R^{\aff}$ by setting $\Delta^{\aff} = \Delta\cup \left\{ \alpha_0\right\}$, where $\alpha_0 =  \delta/\bfe - \vartheta$. In general, we can always form an affine basis $\Delta^{\aff}$ by adjoining affine roots to a given basis of the root system $R(G^{\sigma}, S)$. \par

Define $\fraka_{\bbQ} = \bfX_*(S)$ and $\fraka = \fraka_{\bbQ}\otimes_{\bbQ}\bbR$. We will view $\fraka$ as affine space instead of vector space. Every element $\tilde\alpha = \alpha + (k/\bfe)\delta\in R^{\aff}$ defines on $\fraka$ an affine function by 
\[
	\tilde\alpha: \fraka \to \bbR,\; \lambda\mapsto \langle \alpha, \lambda\rangle + k/\bfe. 
\]
We denote its zero locus by $H_{\tilde\alpha} = \tilde\alpha^{-1}(0)\subset \fraka$. The family $\left\{ H_{\tilde\alpha} \right\}_{\tilde \alpha\in R^{\aff}}$ defines on $\fraka$ a hyperplane arrangement, the set of (open) facets of which is denoted by $\frakF = \frakF(G, \sigma, S)$. It can be identified with the realisation of the Coxeter complex of the affine Weyl group $(W^{\aff}, \Delta^{\aff})$. Note that a facet is affine-homeomorphic to the cartesian product of an affine space and the interior of a polysimplex. Let $\frakE$ be the set of affine subspaces of $\fraka$ generated by elements of $\frakF$. \par

%We can to every point $x\in \fraka_{\bbQ}$ the following subsets 
%\[
	%R^{\aff}_{x \ge 0} = \left\{ \tilde \alpha\in R^{\aff}\;\vert\; \tilde\alpha(x) \ge 0 \right\},\quad  R^{\aff}_{x \le 0} = \left\{ \tilde \alpha\in R^{\aff}\;\vert\; \tilde\alpha(x) \le 0 \right\},\quad R^{\aff}_{x = 0} = R^{\aff}_{x \ge 0}\cap R^{\aff}_{x \le 0}.
%\]
%Given points $x, y\in \fraka_{\bbQ}$, the subset $R^{\aff}_{x \le 0}$ (resp. $R^{\aff}_{x \ge 0}$) coincides with $R^{\aff}_{y \le 0}$ (resp. $R^{\aff}_{y \ge 0}$) if and only if $x$ and $y$ lie in the same facet; $R^{\aff}_{x = 0}$ coincides with $R^{\aff}_{y = 0}$ if and only if the facet which contains $x$ and the one which contains $y$ span the same affine subspaces of $\fraka$. \par

%\subsection{Affine Weyl group}\label{subsec:weyl}
For each $\til \alpha = \alpha + (k/\bfe)\delta\in R^{\aff}$, let $s_{\til\alpha}:\fraka\to \fraka$ be the reflection defined by 
\[
	s_{\til\alpha}(\lambda) = \lambda - \til\alpha(\lambda)\alpha^{\vee}. 
\]
The affine Weyl group of the affine root system $(\fraka, R^{\aff})$, denoted by $W^{\aff}$, is the subgroup affine transforms on $\fraka$ generated by $\left\{ s_{\til \alpha} \right\}_{\til \alpha\in R^{\aff}}$. %Upon choosing a {\it special vertex} $v_0$ of $(\fraka, \frakF)$, there is a decomposition $W^{\aff} = W_0\ltimes \Lambda$, where $W_0 = \Stab_{W^{\aff}}(v_0)$ is a parabolic subgroup of $W^{\aff}$ and $\Lambda$ is the free abelian subgroup of $W^{\aff}$ formed by elements which act by translation on $\fraka$.

\subsection{Facets and spirals}
Assume that the pinning $E$ satisfies \autoref{prop:graduation} and $\theta^{\mathrm{in}}:\mu_m\to G^{\ad}$ is the inner component of $\theta$. We fix the choice of a cocharacter $\iota_0\in \bfX_*(S^{\ad})$ which extends $\theta^{\mathrm{in}}$ in the sense that $\iota_0\mid_{\mu_m} = \theta^{\mathrm{in}}$. Let $\frakP^{\epsilon}_S\subseteq \frakP^{\epsilon}$ denote the subset of $\epsilon$-spirals $\frakp_*$ satisfying $S\subseteq P_0$ and $\frakG_S\subseteq \frakG$ the subset of splitting factors $\frakl_*$ satisfying $S\subseteq L_0$. 

\begin{prop}\label{prop:facet-spiral}
	The map $\frakF\to \frakP^{\epsilon}_S$ which sends a facet $F\in \frakF$ to ${^\epsilon}\frakp^{\lambda_y}_*$ with $\lambda_y = (my - \iota_0)/m$ for any $y\in F\cap \fraka_{\bbQ}$ is a well-defined bijection. Similarly, the map $\frakE\to \frakG_S$ which sends $E$ to $\frakl^E_* := \frakl^{\lambda_y}_*$ for any $y\in \fraka_{\bbQ}$ satisfying $E = \bigcap_{y\in E'\in \frakE} E'$ is a well-defined bijection.
\end{prop}
Note that the bijections depend on the choice of $\iota_0$.

\subsection{Relative affine root system}\label{subsec:rars}
Fix an LY-cuspidal support $\xi = (L, \frakl_*, \scrC)\in \frakT$. Up to taking $\Gz$-conjugation, we may assume that $S\subseteq L_0$. There is an affine subspace $\fraka_{\xi}\in \frakE$ such that $\frakl_* = \frakl_*^{\fraka_{\xi}}$ under the correspondence of \autoref{prop:facet-spiral}. The Weyl group $W_L = W(L, S)$ can be identified with the pointwise stabiliser $\Stab_{W^{\aff}}(\fraka_{\xi})$ in $\fraka$. The facets of $\fraka$ which span $\fraka_{\xi}$ as affine subspace are called $\fraka_{\xi}$-alcoves. Let $\Xi$ denote the set of $\fraka_{\xi}$-alcoves. 
\begin{prop}[{\cite[\S 2]{lusztig95c}}]\label{prop:rel-weyl}
	The restriction
	\[
		R^{\aff}\mid_{\fraka_{\xi}} = \left\{ \alpha\!\mid_{\fraka_{\xi}}\;\vert\; \alpha\in R^{\aff},\; \alpha\!\mid_{\fraka_{\xi}} \text{ non-constant} \right\}
	\]
	is a (possibly non-reduced) affine root system on $\fraka_{\xi}$. Moreover, %let $N_{W^{\aff}}(\fraka_{\xi})$ be the subgroup of $W^{\aff}$ consisting of elements which preserve the subspace $\fraka_{\xi}\subseteq \fraka$; then 
	the affine Weyl group of $(\fraka_{\xi}, R^{\aff}\mid_{\fraka_{\xi}})$, denoted by $W^{\aff}_{\xi}$, can be identified with $N_{W^{\aff}}(W_L) / W_L$ %\{w\in W^{\aff}\;\vert\; w\fraka_{\xi} = \fraka_{\xi}\} / \{w\in W^{\aff}\;\vert\; wy = y,\;\forall y\in \fraka_{\xi}\}
	and it acts simply transitively on $\Xi$. \hfill\qedsymbol
\end{prop}
We let $R^{\aff}_{\xi}\subseteq R^{\aff}\mid_{\fraka_{\xi}}$ be the subset of indivisible roots, so that $(\fraka_{\xi}, R^{\aff}_{\xi})$ is a reduced affine root system, called \emph{relative affine root system}. 
The group $W^{\aff}_{\xi}$ defined above is called the \emph{relative affine Weyl group} for the LY-cuspidal support $\xi$. Recall the cocharacter $\iota_0\in \bfX_*(S^{\ad})$ from \autoref{subsec:grading}. Set $\bfx_0 = (\iota_0 / m, 1)\in \fraka$ and let $\bfx\in \fraka_{\xi}$ be the image of $\bfx_0$ under the orthogonal projection $\fraka\to \fraka_{\xi}$ and let $W^{\aff}_{\xi, \bfx} = \Stab_{W^{\aff}_{\xi}}(\bfx)$ denote the stabiliser of $\bfx$ in $W^{\aff}_{\xi}$.  \par

Let $\frakP^{\epsilon}_{\xi}\subseteq \frakP^{\epsilon}_S$ be the set of $\epsilon$-spirals of $\frakg_*$ which have $\frakl_*$ as splitting factor. We define $\ubar\frakP^{\epsilon}_{\xi}$ to be the set of $\Gz$-conjugacy classes in $\frakP^{\epsilon}_{\xi}$.
\begin{prop}\label{prop:FP}
	The bijection $\frakF\xrightarrow{\sim} \frakP^{\epsilon}_S$ provided by \autoref{prop:facet-spiral} induces further bijections:
	\[
		\Xi\xrightarrow{\sim} \frakP^{\epsilon}_{\xi}\quad \text{and}\quad \Xi / W^{\aff}_{\xi, \bfx} \xrightarrow{\sim} \ubar\frakP^{\epsilon}_{\xi}.
	\]
\end{prop}
\begin{proof}
	The first map is clearly bijective. The second follows from \autoref{lemm:weyl-stab} below, which is stated without proof in \cite{liu21}.
\end{proof}
\begin{lemm}\label{lemm:weyl-stab}
	With the identification $\fraka_{\xi} = \bfx + \bfX_*(Z_L^{\circ})\otimes \bbR$, the conjugacy action of $N_{\Gz}(Z_L^{\circ})$ on $\bfX_*(Z_L^{\circ})\otimes \bbR$ induces an isomorphism
	\[
		N_{\Gz}(Z_L^{\circ}) / L_0 \xrightarrow{\sim} W^{\aff}_{\xi, \bfx}.
	\]
\end{lemm}
\begin{proof}
	Given $g\in N_{\Gz}(Z_L^{\circ})$, there exists $h\in L_0$ such that $hgS(hg)^{-1} = S$ since $S$ and $gSg^{-1}$ are both maximal tori of $L_0$. In other words, the natural map
	\[
		N_{\Gz}(Z_L^{\circ}) \cap N_{\Gz}(S)\to N_{\Gz}(Z_L^{\circ}) / L_0
	\]
	is surjective. The identification $\fraka = \bfx_0 + \bfX_*(S)\otimes \bbR$ and the conjugacy action of $N_{\Gz}(S)$ on $\bfX_*(S)\otimes \bbR$ induces an isomorphism $N_{\Gz}(S)/S \cong W^{\aff}_{\bfx_0}$ by \cite[Corollary 4.2]{liu21}; therefore, we have 
	\[
		(N_{\Gz}(Z_L^{\circ}) \cap N_{\Gz}(S) )/ S = N_{W^{\aff}}(W_L)\cap W^{\aff}_{\bfx_0}. 
	\]
	It remains to show that the isomorphism $N_{W^{\aff}}(W_L)/W_L \xrightarrow{\sim} W^{\aff}_{\xi}$ from \autoref{prop:rel-weyl} restricts to an isomorphism:
	\begin{equation}\label{equa:stab-W}
		(N_{W^{\aff}}(W_L)\cap W^{\aff}_{\bfx_0}) / (W_{L}\cap W^{\aff}_{\bfx_0}) \xrightarrow{\sim} W^{\aff}_{\xi, \bfx}.
	\end{equation}
	The map in question is clearly injective. It suffices to show the surjectivity.
	Pick an $\fraksl_2$-triplet $(e, h, f)$ such that $e\in \ci\frakl_d$, $h\in \fraks = \Lie S$ and $f\in \frakl_{-d}$, which yields a cocharacter $\iota\in \bfX_*(S)$ determined by $\rmd \iota(1) = h$. \autoref{prop:cuspidal-lusztig} implies that $\frakl_n = {^{\iota}_{2n/d}}\frakl$ holds for $n\in \bbZ$ by. Set
	\[
		\frakh_* = \bigoplus_{n\in \bbZ} \frakh_n,\quad \frakh_n = {^\iota_{2n/d}}\frakg_{\ubar n}
	\]
	and let $H$ be the pseudo-Levi subgroup of $(G, \sigma)$ with $\Lie H = \frakh$. Under \autoref{prop:facet-spiral}, $\frakh_*\in \frakG_S$ corresponds to a subspace $\fraka_H\in \frakE$ of $\fraka$, which is the common zero locus of the following subset of roots:
	\[
		R_H := \left\{ \alpha\in R^{\aff}\;\vert\; \alpha(\bfx_0 - (d/2)\iota) = 0 \right\}.
	\]
	It is not hard to show that $\bfx = \bfx_0 - (d/2)\iota$.
	Let $W_H = N_H(S)/S$ and $W_L = N_L(S)/S$ be the Weyl group of $H$ and $L$, respectively. We then have $W_H = W^{\aff}_{\bfx}$ and thus $N_{W_H}(W_L) / W_L \cong W^{\aff}_{\xi, \bfx}$. 

	Choose a $\bbZ$-graded Borel subalgebra $\frakb_*\subseteq \frakh_*$ such that $\frakb_*\cap \frakl_*$ is a Borel subalgebra of $\frakl_*$ and $\frakb \cap \frakl\subseteq {^\iota_{\ge 0}}\frakl$ holds. Given $w\in W^{\aff}_{\xi, \bfx}$, we let $\dot w\in N_{W_H}(W_L)$ be the unique lifting of $w$ which is the shortest representative of its $W_L$-coset with respect to the simple reflections of $W_H$ corresponding to the basis $\Delta\subseteq R_H$ attached to the Borel subalgebra $\frakb$, so that $\dot w (\Delta \cap R_L) = \Delta \cap R_L$ holds. We shall show that $\dot w \bfx_0 = \bfx_0$, which will imply the surjectivity of \eqref{equa:stab-W} and complete the proof. The cocharacters $\iota$ and $\dot w\iota$ are both dominant with respect to the basis $\Delta\cap R_L$ --- indeed $\frakb \cap \frakl\subseteq {^\iota_{\ge 0}}\frakl$ implies that $\iota$ is dominant and 
	\[
		\langle \dot w\iota, \alpha\rangle = \langle \iota, w^{-1}\alpha\rangle \ge 0,\quad \forall \alpha\in \Delta\cap R_{L}
	\]
	implies that $\dot w\iota$ is dominant. Let $\ddot w\in N_H(S)\cap N_H(Z^{\circ}_L)$ be a lifting of $\dot w$. The cuspidal local system $\scrC$ is isomorphic to $\scrC^{\circ}\mid_{\ci\frakl_d}$ for some $L$-equivariant local system $\scrC^{\circ}$ on an $L$-orbit $\rmO^{\circ}\subseteq \frakl$. The cuspidality of $(\rmO^{\circ}, \scrC^{\circ})$ implies that $\Ad_{\ddot w}$ preserves the orbit $\rmO^{\circ}$, so that $\Ad_{\ddot w}e\in \Ad_{L}e$ holds. An application of \cite[5.6.8]{carter05} to the $\fraksl_2$-triplets $(e, h, f)$ and $(\Ad_{\ddot w} e, \Ad_{\ddot w}h, \Ad_{\ddot w} f)$ implies that $\langle \dot w \iota, \alpha\rangle = \langle\iota, \alpha\rangle$ holds for every $\alpha\in \Delta\cap R_L$, which implies that $\dot w\iota = \iota$ because $\iota$ and $\dot w\iota$ both lie in $\bfX_*([L,L]\cap S)$. Now, $\dot w\iota = \iota$, $\dot w \bfx = \bfx$ and $\bfx_0 = \bfx + (d/2m) \iota$ together imply that $\dot w \bfx_0 = \bfx_0$ holds. 
\end{proof}

\subsection{Sheaf-theoretic realisation of double affine Hecke algebras}\label{subsec:springer}
For the moment, we assume additionally that $G$ is \emph{almost simple} and \emph{simply connected}. Given an LY-cuspidal support $\xi = (L, \frakl_*, \scrC)$, consider the completed extension algebra
\[
	\calA = \bigoplus_{\frakp_*, \frakq_*\in \ubar\frakP^{\epsilon}_{\xi}} \Hom^{\bullet}_{\Gz}(\Ind^{\frakg_{\ubar d}}_{\frakp_d}\IC(\scrC), \Ind^{\frakg_{\ubar d}}_{\frakq_d}\IC(\scrC))^{\wedge}_0.
\]
It is equipped with a structure of non-unital ring via the Yoneda product. Via the adjunction 
\[
	\Hom^{\bullet}_{\Gz}(\Ind^{\frakg_{\ubar d}}_{\frakp_d}\IC(\scrC), \Ind^{\frakg_{\ubar d}}_{\frakq_d}\IC(\scrC))\cong \Hom^{\bullet}_{L_0}(\Res^{\frakg_{\ubar d}}_{\frakq_d}\Ind^{\frakg_{\ubar d}}_{\frakp_d}\IC(\scrC), \IC(\scrC)),
\]
the $\rmH^{\bullet}_{\Gz}$-module structure on $\calA$ factorises through the natural maps
\[
	\rmH_{\Gz}^{\bullet}\to \rmH^{\bullet}_{L_0}\to \Hom^{\bullet}_{L_0}(\IC(\scrC), \IC(\scrC))\cong \rmH^{\bullet}_{Z^{\circ}_{L}}. 
\]
\autoref{lemm:weyl-stab} implies that the image of the map $\rmH_{\Gz}^{\bullet}\to \rmH^{\bullet}_{Z^{\circ}_{L}}$ coincides with the subalgebra of invariants $(\rmH^{\bullet}_{Z^{\circ}_{L}})^{W^{\aff}_{\xi, \bfx}}$. Consequently, $\calA$ has a natural $(\rmH^{\bullet\wedge}_{Z^{\circ}_{L}})^{W^{\aff}_{\xi, \bfx}}$-module structure. \par

On the other hand, we have a relative affine root system $(\fraka_{\xi}, R^{\aff}_{\xi})$ defined in \autoref{subsec:rars} associated with $\xi$. Fix a basis $\Delta^{\aff}\subset R^{\aff}_{\xi}$, which corresponds to an $\fraka_{\xi}$-alcove $\nu$ and thus to the choice of spiral $\frakp_*\in \frakP^{\epsilon}_{\xi}$ corresponding to $\nu$ (see \autoref{prop:FP}). In \cite{LYIII}, a degenerate double affine Hecke algebra $\bfH_{\xi}$ for $(\fraka_{\xi}, R^{\aff}_{\xi}, \Delta^{\aff})$ is attached to $\xi$ with a suitable set of parameters $(c_{\alpha})_{\alpha\in R^{\aff}_{\xi}}$, see also \autoref{subsec:dDAHA}. The spectral completion (\autoref{subsec:spectral}) of $\bfH_{\xi}$ at the point $\bfx\in \fraka_{\xi}$ (\autoref{subsec:rars}) yields a non-unital ring $\calH_{\xi} = \calC^{\bfx}\bfH_{\xi}$ with categorical centre $\calZ_{\xi} = \End(\id_{\calH_{\xi}\Mod})$. The ring $\calH_{\xi}$ has a complete orthogonal family of idempotents $(1_x)_{x\in W^{\aff}\bfx}$.  For each $w\in W^{\aff}_{\xi}$, let ${^w}\frakp_*\in \frakP^{\epsilon}_{\xi}$ be the spiral corresponding to the $\fraka_{\xi}$-alcove $w\nu$. The main construction of \cite{liu21} can be paraphrased as follows:

\begin{theo}\label{theo:springer}
	Assume that $G$ is simply connected and almost simple. There is an isomorphism of (non-unital) rings $\calH_{\xi} \cong \calA$ satisfying the following properties:
	\begin{enumerate}
		\item
			It intertwines the module structures over $\calZ_{\xi} \cong (\rmH^{\bullet \wedge}_{Z_L^{\circ}})^{W_{\xi, \bfx}^{\aff}}$.
		\item
			For each $w, y\in W^{\aff}_{\xi}$, it restricts to an isomorphism of subspaces:
			\[
				1_{y^{-1} \bfx}\calH_{\xi} 1_{w^{-1} \bfx} \xrightarrow{\sim}\Hom^{\bullet}_{\Gz}(\Ind^{\frakg_{\ubar d}}_{{^{w}}\frakp_d}\IC(\scrC), \Ind^{\frakg_{\ubar d}}_{{^{y}}\frakp_d}\IC(\scrC))^{\wedge}_0.
			\]
	\end{enumerate}\hfill\qed
\end{theo}
The following is proven in \cite{liu21} as a consequence of \autoref{theo:springer}:
\begin{theo}\label{coro:springer}
	Let $\epsilon = d / |d|$. There is a canonical bijective correspondence between $\Pi_{\xi}$ and $\Irr \calH_{\xi}\mof$ which attaches to $\pi\in \Pi_{\xi}$ the following $\calH_{\xi}$-module:
	\[
		\calL_\pi := \operatorname{cosoc}_{\calH_{\xi}}\left(\bigoplus_{\frakp_*\in \ubar\frakP^{\epsilon}}\Hom^{\bullet}_{\Gz}(\IC_{\pi}, \Ind^{\frakg_{\ubar d}}_{\frakp_d}\IC(\scrC))^{\wedge}_0\right),
	\] 
	\hfill\qedsymbol
\end{theo}

\subsection{Indecomposability of blocks}
We drop now the assumption that $G$ is simply connected and almost simple.
\begin{prop}\label{lemm:HomFG}
	Given $\xi = (L, \frakl_*, \scrC)\in \frakT$ and $\pi, \pi'\in \Pi_{\xi}$, we have $\Hom^{\bullet}_{\Gz}(\IC_{\pi}, \IC_{\pi'})\neq 0$.
\end{prop}
\begin{proof}
	Let $E = (B, T, a)$ be a pinning for $G$ satisfying \autoref{prop:graduation}. Moreover, we assume that $\frakl_*$ is $S$-stable. We shall reduce the question to the case where $G$ is a simply connected almost simple group. 
	\begin{itemize}
		\item[Step 0] (reduction to $d = 1$)
			Let $c = \gcd(m, d)$ and let $a, b\in \bbZ$ be such that $ad + bm = c$. Set $H = G^{\theta^{m/c}}$. Then, the morphism 
			\[
				\theta_H : \mu_{m/c}\to \Aut(H),\quad \theta_H(\omega^{c})h = \theta(\omega^a) h\quad \text{for $\omega\in \mu_m$, $h\in H$}
			\]
			induces on a $\bbZ/(m/c)\bbZ$-grading on the Lie algebra $\frakh$. Moreover, we have $G_{\ubar 0} = H_{\ubar 0}$ and $\frakg_{\pm \ubar d} = \frakh_{\pm\ubar 1}$. Let $\xi_H = (L\cap H, \frakl_*\cap \frakh_*, \scrC)$. It is an LY-cuspidal support on $\frakh_{\ubar 1}$. We may regard $\pi$ and $\pi'$ as in $\Pi_{\xi_H}$ and we have
	\[
		\Hom_{\Gz}^{\bullet}(\IC_{\pi}, \IC_{\pi'}) \cong \Hom_{H_{\ubar 0}}^{\bullet}(\IC_{\pi}, \IC_{\pi'}).
	\]
	Therefore, we may assume that $d = 1$.
\item[Step 1] (reduction to simply connected semisimple group) Let $Z = Z_G$ and $G' = [G, G]$. Let $\til G$ be the universal cover of $G'$ and $\til\frakg = \Lie \til G$ its Lie algebra. We have a canonical decomposition $\frakg = \til\frakg \oplus \frakz$. Moreover, the map $\theta:\mu_m\to \Aut(G)$ admits a unique lifting $\til\theta:\mu_m\to \Aut(\til G)$ and yields a $\bbZ/m$-grading on $\til\frakg$. It induces an isomorphism $\til\frakg^{\nil}_{\ubar 1}\cong\frakg^{\nil}_{\ubar 1}$. Let $\til\xi = (L\times_G\til G, \frakl_*\cap \til\frakg_*, \scrC)\in \frakT_{\til G}$. We may regard $\pi$ and $\pi'$ as in $\Pi_{\til\xi}$ and we have 
	\[
		\Hom_{\Gz}^{\bullet}(\IC_{\pi},\IC_{\pi'}) \cong \Hom_{\til G_{\ubar 0}}^{\bullet}(\IC_{\pi},\IC_{\pi'})\otimes\rmH^\bullet_{Z_{\ubar 0}}.
	\]
	Therefore, we may assume that $G$ is semisimple and simply connected.
		\item[Step 2] (reduction to almost simple group) As $G$ is assumed to be simply connected and semisimple, it admits a decomposition
			\[
				G = G^{(1)}\times \cdots  \times G^{(n)},\quad \frakg = \frakg^{(1)}\times \cdots  \times \frakg^{(n)},
			\]
			where $G^{(i)}$ is an almost simple factor of $G$ and $\frakg^{(i)} = \Lie G^{(i)}$. The pinned automorphism $\sigma$ induces a permutation $s\in \mathfrak{S}_n$ such that $\sigma(G^{(i)}) = G^{(s(i))}$ for each $i\in [1, n]$. We may order these factors so that 
			\[
				s = (a_1 \cdots a_1 - 1) (a_1 \cdots a_2 - 1)\cdots (a_{r-1} \cdots a_{r+1} - 1)
			\]
			for some numbers $1 = a_1 < \cdots < a_{r+1} = n + 1$. For $i \in [1, r]$, let 
			\[
				\frakg^{[i]} = \frakg^{(a_{i})}\times \cdots\times \frakg^{(a_{i+1} - 1)},\quad d_i = m / (a_{i+1} - a_{i}).
			\]
			Equip $\frakg^{[i]}$ with the $\bbZ/m$-grading induced by $\theta\mid_{\frakg^{[i]}}$ and $\frakg^{(a_{i})}$ the $\bbZ / d_i$-grading induced by $\theta^{d_i}\mid_{\frakg^{(a_{i})}}$. It follows that there are isomorphisms:
			\[
				G^{[i]}_{\ubar 0} \cong G^{(a_{i})}_{\ubar 0},\quad \frakg^{[i]}_{\ubar 1} \cong \frakg^{(a_{i})}_{\ubar 1}
			\]
			which intertwine the adjoint actions. \par
			Moreover, the LY-cuspidal support $\xi$ can be decomposed $\xi \cong \xi^{[1]}\boxtimes \cdots \boxtimes\xi^{[r]}$, where $\xi^{[i]} = (L^{[i]}, \frakl^{[i]}_*, \scrC^{[i]})\in \frakT_{G^{(a_i)}}$ with $L^{[i]} = L\cap G^{(a_i)}$ and $l^{[i]}_* = \frakl_*\cap \frakg^{(a_i)}_*$.

			The pairs $\pi$ and $\pi'$ can also be decomposed:
			\[
				\pi \cong \pi^{[1]}\boxtimes \cdots\boxtimes\pi^{[r]},\quad \pi' \cong \pi'^{[1]}\boxtimes \cdots\boxtimes\pi'^{[r]}, 
			\]
			where $\pi^{[i]}, \pi'^{[i]}\in \Pi_{\xi^{[i]}}$ are pairs on $\frakg^{(a_i)}_{\ubar 1}$. The K\"unneth formula yields
			\[
				\Hom_{\Gz}^{\bullet}(\IC_{\pi},\IC_{\pi'}) \cong \Hom_{G^{(a_1)}_{\ubar 0}}^{\bullet}(\IC_{\pi^{[1]}},\IC_{\pi'^{[1]}})\otimes \cdots\otimes\Hom_{G^{(a_r)}_{\ubar 0}}^{\bullet}(\IC_{\pi^{[r]}},\IC_{\pi'^{[r]}}).
			\]
			Therefore, we may assume furthermore that $G$ is almost simple and simply connected.
		\item[Step 3]
			Set $\epsilon = d / |d| = 1$. Choose a spiral $\frakp_*\in \frakP^{\epsilon}_{\xi}$ (resp. $\frakq_*\in \frakP^{\epsilon}_{\xi}$) of $\frakg_*$ such that $\IC_{\pi}$ (resp. $\IC_{\pi'}$) is a direct summand of $\pH^\bullet\Ind^{\frakg_{\ubar d}}_{\frakp_d}\IC(\scrC)$ (resp. of $\pH^\bullet\Ind^{\frakg_{\ubar d}}_{\frakq_d}\IC(\scrC)$). By \autoref{theo:springer}, the completed Ext-space
			\[
				\Hom^{\bullet}_{\Gz}(\Ind^{\frakg_{\ubar d}}_{\frakp_d} \IC(\scrC), \Ind^{\frakg_{\ubar d}}_{\frakq_d}\IC(\scrC))^{\wedge}_0
			\]
			is isomorphic to $1_y \calH_{\xi} 1_x$ for some $x, y\in W^{\aff}_{\xi}\bfx$. Similarly, the completed Ext-algebra
			\[
				\calE := \End^{\bullet}_{\Gz}(\Ind^{\frakg_{\ubar d}}_{\frakp_d}\IC(\scrC)\oplus \Ind^{\frakg_{\ubar d}}_{\frakq_d}\IC(\scrC))^{\wedge}_0.
			\]
			is isomorphic to 
			\[
				(1_x\oplus 1_y)\calH_{\xi} (1_x \oplus 1_y) = 1_x\calH_{\xi} 1_x\oplus 1_x\calH_{\xi} 1_y\oplus 1_y\calH_{\xi} 1_x\oplus 1_y\calH_{\xi} 1_y,
			\]
			which is free over $\calZ_{\xi}$. By \autoref{prop:HK}, the idempotent subalgebra 

			\[
				\calE_K := K\otimes_{\calZ_{\xi}} \calE \cong (1_x \oplus 1_y) (K\otimes_{\calZ_{\xi}}\calH_{\xi}) (1_x \oplus 1_y)
			\]
			is isomorphic to a matrix algebra over $K$. The spaces 
			\begin{align*}
				\calM_{\pi} := \Hom^{\bullet}_{\Gz}(\IC_\pi, \Ind^{\frakg_{\ubar d}}_{\frakp_d} \IC(\scrC) \oplus \Ind^{\frakg_{\ubar d}}_{\frakq_d} \IC(\scrC))^{\wedge}_0, \\
				\calM_{\pi'} := \Hom^{\bullet}_{\Gz}(\IC_{\pi'}, \Ind^{\frakg_{\ubar d}}_{\frakp_d} \IC(\scrC) \oplus \Ind^{\frakg_{\ubar d}}_{\frakq_d} \IC(\scrC))^{\wedge}_0
			\end{align*}
			are non-zero projective left $\calE$-modules by the choice of $\frakp_*$ and $\frakq_*$; in particular, they are free over $\calZ_{\xi}$. Consequently, the localised modules $\calM_{\pi,K}$ and $\calM_{\pi',K}$ are non-zero. As $\calE_K$ is Morita-equivalent to the base field $K$, we have 
			\[
				K\otimes_{\calZ_{\xi}}\Hom^{\bullet}_{\Gz}(\IC_{\pi}, \IC_{\pi'})^{\wedge}_0\cong \Hom_{\calE_K}(\calM_{\pi',K}, \calM_{\pi,K}) \neq 0.
			\]
			Therefore, $\Hom^{\bullet}_{\Gz}(\scrF, \scrG)\neq 0$ holds.
	\end{itemize}

\end{proof}

\begin{coro}\label{coro:indecomp}
	For each LY-cuspidal support $\xi\in \frakT$, the $\xi$-block $\Db_{\Gz}(\frakg_{\ubar d}^{\nil})_{\xi}$ is an indecomposable triangulated category. 
\end{coro}
\begin{proof}
	Since $\Db_{\Gz}(\frakg_{\ubar d}^{\nil})_{\xi}$ is generated as thick subcategory by the simple perverse sheaves $\left\{ \IC_{\pi} \right\}_{\pi\in \Pi_{\xi}}$, the indecomposability follows from \autoref{lemm:HomFG}.
\end{proof}

\iffalse
\begin{lemm} \hl{Modification n\'ecessaire pour le cas tordu!}
	The following statements hold:
	\begin{enumerate}
		\item
			We have $R_{\zeta} = R(G, Z_M^{\circ})$ and $(\frake_{\zeta}, R_{\zeta})$ is a (possibly non-reduced) root system.
		\item
			The Weyl group $W_{\xi}$ can be identified with $N_{G}(Z^{\circ}_M) / M$ via their actions on $\frake_{\xi}$.
		\item
			The action of $W^{\aff}_{\xi}$ on the tangent space $\frake_{\xi}$ induces an isomorphism $W^{\aff}_{\xi} / \bfX_*(Z^{\circ}_M)\xrightarrow{\sim}W_{\xi}$.
	\end{enumerate}
\end{lemm}

\fi

\section{Study of primitive pairs}\label{sec:primitive}
This section serves as a preparation for \autoref{sec:supercuspidal}. We study the primitive pairs introduced in \cite[3.8]{LYI}. The key technical result is \autoref{prop:primitive-induction}, which asserts that $\Pi$ is generated by stalks of parabolically induced complexes from primitive pairs.

\subsection{A geometric lemma}
We will make use of the following geometric lemma:
\begin{lemm}\label{lemm:U-conj}
	Let $\lambda\in \bfX_*(\Gz)_{\bbQ}$ be a fractional cocharacter and set $U_\lambda = \exp({^\lambda_{> 0}}\frakg_{\ubar 0})$. Then, for every integer $e\in \bbZ$, and every $x\in {^\lambda_e}\frakg_{\ubar d}$, we have 
	\[
		\Ad_{\Gz}x\cap (x + {^\lambda_{> e}}\frakg_{\ubar d}) = \Ad_{U_\lambda}x.
	\]
\end{lemm}
\begin{proof}
	Up to replacing $\lambda$ and $e$ with a positive multiple, we may assume that $\lambda\in \bfX_*(\Gz)$. The inclusion $\Ad_{U_\lambda}x\subseteq\Ad_{\Gz}x\cap (x + {^\lambda_{> e}}\frakg_{\ubar d})$ is clear. We claim that for every $f \in \bbZ_{> e}$ and every $x'\in  \Ad_{\Gz}x\cap (x + {^\lambda_{\ge f}}\frakg_{\ubar d})$, there exists $u\in U_\lambda$ such that $\Ad_u x'\in x + {^\lambda_{> f}}\frakg_{\ubar d}$; this will prove the lemma by induction on $f$. Consider the morphism
	\[
		\varphi:\bbC\to \Ad_{\Gz}x\cap (x + {^\lambda_{\ge f}}\frakg_{\ubar d}),\quad t\mapsto \begin{cases}t^{-e}\Ad_{\lambda(t)}x' & t\in \bbC^{\times} \\ \lim_{t\to 0}t^{-e}\Ad_{\lambda(t)}x' = x & t = 0\end{cases}.
	\]
	The power series expansion of $\varphi$ at $0$ yields
	\[
	\varphi(t)\sim x + t^{f - e}x'_f + o(t^{f - e}),\quad \text{where } x'_f =\lim_{s\to 0}s^{-f}\Ad_{\lambda(s)}(x' - x)\in {^\lambda_f}\frakg_{\ubar d}.
	\]
	Consequently, the element $x'_f$ lies in the tangent space of $\Ad_{\Gz}x$ at $x$, which equals $[\frakg_{\ubar 0}, x]$. As $x'_f\in [\frakg_{\ubar 0}, x]\cap {^\lambda_f}\frakg_{\ubar d} = [{^\lambda_{f-e}}\frakg_{\ubar 0}, x]$, we may find $X\in{^\lambda_{f-e}}\frakg_{\ubar 0}$ such that $x'_f = [X, x]$. It follows that
	\[
		\Ad_{\exp(-X)}x' = \exp(-\ad_{X})x' \in x' - x'_f + {^\lambda_{>f}}\frakg_{\ubar d} = x  + {^\lambda_{>f}}\frakg_{\ubar d}
	\]
	and $\exp(-X)\in U_{\lambda}$ as claimed.
\end{proof}
\subsection{Primitive pairs and their cleanness}\label{subsec:primitive} 
We keep the setup of \autoref{subsec:G}.
Let $\xi = \left( L, \frakl_*, \scrC \right)\in \frakT$ be an LY-cuspidal support on $\frakg_{\ubar d}$. Let $\rmO = \Ad_{\Gz}\ci\frakl_d\subseteq \frakg_{\ubar d}^{\nil}$ be the $\Gz$-orbit containing $\ci\frakl_d$ and let $z\in \ci\frakl_d$. By~\cite[3.8(a)]{LYI}, it is known that the inclusion $Z_{L_{0}}(z)\subseteq Z_{\Gz}(z)$ induces an isomorphism on groups of connected components. Therefore, there is a unique extension of $\scrC$ into a $\Gz$-local system on $\rmO$, denoted by $\scrS$. The pair $(\rmO, \scrS)$,  called \emph{primitive pair} attached to $\xi$, lies in the $\xi$-series $\Pi_\xi$. The primitive pair is characterised by the minimality of its support:
\begin{lemm}\label{lemm:supp}
	Let $\pi = (\rmC, \scrL)\in \Pi_{\xi}$ be a pair in the $\xi$-block. Then, $\rmO\subseteq \ba\rmC$ holds. If $\rmC = \rmO$ holds, then $\pi$ is isomorphic to $(\rmO, \scrS)$. 
\end{lemm}
\begin{proof}
	Suppose first that $\rmO\not\subseteq \ba\rmC$. We have $\rmO\cap \ba\rmC = \emptyset$. Let $j_{\rmC}: \rmC\hookrightarrow \frakg_{\ubar d}$ denote the inclusion of orbit. Set $\epsilon = d / |d|$ and let $\frakp_*\in \frakP^{\epsilon}$ be an $\epsilon$-spiral which has $\frakl_*$ as splitting factor. We may choose $\lambda\in \bfX_*(\Gz)$ and $k \in \bbZ_{>0}$ such that $\frakp_n = {^{\epsilon(\lambda/k)}_{\ge \epsilon n}}\frakg_{\ubar n}$ holds for $n\in \bbZ$ (\autoref{subsec:spiral}). The support of the complex $\Res_{\frakp_d}^{\frakg_{\ubar d}} j_{\rmC!}\scrL$ is contained in $\frakl_d \setminus \ci\frakl_d$ --- indeed, the image of $\rmC \cap \frakp_d$ under the projection $\frakp_d \to \frakl_d$ is contained in $\overline\rmC\cap \frakl_d\subseteq \frakl_d \setminus \ci\frakl_d$ because the projection $\frakp_d\to \frakl_d$ is given by 
	\[
		x\mapsto \lim_{t\to 0}t^{-kd\epsilon}\Ad_{\lambda(t)^{\epsilon}}x
	\]
	and the $\Gz$-orbit $\rmC$ is $\bbC^{\times}$-stable (see \autoref{subsec:JM}). Thus, the cleanness of $(\ci\frakl_d, \scrC)$ provided by \autoref{prop:cuspidal-lusztig} implies that 
	\begin{eq}
		\Hom^{\bullet}_{\Gz}\left(j_{\rmC!}\scrL, \Ind^{\frakg_{\ubar d}}_{\frakp_{d}} \IC(\scrC ) \right) \cong \Hom^{\bullet}_{L_{0}}\left(\Res_{\frakp_{d}}^{\frakg_{\ubar d}} j_{\rmC!}\scrL, \IC\left( \scrC \right)\right)\cong \Hom^{\bullet}_{L_{0}}\left(\left( \Res_{\frakp_{d}}^{\frakg_{\ubar d}} j_{\rmC!}\scrL\right)\mid_{\ci\frakl_d}, \scrC[\dim \frakl_d] \right) = 0.
	\end{eq}
	Since $j_{\rmC!}\scrL\in \Db_{\Gz}(\frakg_{\ubar d})_{\xi}$ by \autoref{lemm:ext-block} and $\{\Ind^{\frakg_{\ubar d}}_{\frakp_{d}} \IC(\scrC )\}_{\frakp_*\in \frakP^{\epsilon}}$ generates $\Db_{\Gz}(\frakg^{\nil}_{\ubar d})_{\xi}$ as thick subcategory, it follows that $j_{\rmC!}\scrL = 0$ and hence $\scrL = j_{\rmC}^*j_{\rmC!}\scrL = 0$, contradiction. It follows that $\rmO\subseteq \ba\rmC$ holds. \par
	Assume furthermore that $\rmC = \rmO$. Let $\frakp_*$ be as above. The complex $\left(\Res^{\frakg_{\ubar d}}_{\frakp_d}j_{\rmO!}\scrL\right)\mid_{\ci\frakl_d}$ is isomorphic to $\scrL\mid_{\ci\frakl_d}[-2D]$ for some $D\in \bbN$ --- indeed, $(\ci\frakl_d + \fraku_d)\cap \rmO$ is an affine-space bundle over $\ci\frakl_{d}$ by \autoref{lemm:U-conj} (with $\epsilon \lambda$ and $e = |d|$). 
	If $\scrL\not\cong \scrS$, then $\scrL\mid_{\ci\frakl_d}\not\cong \scrC$. The adjunction and the cleanness of $\scrC$ yields
	\[
		\Hom^{\bullet}_{\Gz}\left( j_{\rmO!}\scrL, \Ind^{\frakg_{\ubar d}}_{\frakp_d} \IC(\scrC) \right)\cong \Hom^{\bullet}_{L_{0}}\left( \Res^{\frakg_{\ubar d}}_{\frakp_d}j_{\rmO!}\scrL, \IC(\scrC) \right)\cong \Hom^{\bullet}_{L_{0}}\left( \scrL\mid_{\ci\frakl_d}, \scrC[\dim \frakl_d+2D] \right) = 0.
	\]
	This leads to $\scrL = 0$, contradiction. It follows that $\scrL\cong \scrS$ holds. 
\end{proof}

We say that a pair $(\rmO, \scrS)\in \Pi$ is \emph{clean} if the natural morphism $j_!\scrS\to j_*\scrS$ is an isomorphism, where $j:\rmO\hookrightarrow \frakg_{\ubar d}$ is the inclusion of nilpotent orbit.
\begin{prop}\label{prop:cleanness}
	If $(\rmO, \scrS)$ is the primitive pair attached to an LY-cuspidal support $\xi$, then $(\rmO, \scrS)$ is clean.
\end{prop}
\begin{proof}
	Let $\scrK = i_*i^!j_!\scrS$. By \autoref{lemm:ext-block}, $\scrK$ lies in $\Db_{\Gz}( \frakg_{\ubar d}^{\nil})_{\xi}$. Applying~\autoref{lemm:supp} to the perverse cohomology of $\scrK$, we see that $\scrK = 0$ and thus $i^!j_!\scrS = 0$. This implies the cleanness.
\end{proof}

\subsection{Distinguished primitive pairs}\label{subsec:dist-prim}
Of special interest are the primitive pairs $(\rmO, \scrS)$ such that $\rmO$ is \emph{distinguished}. Set $Z_{\ubar 0} = Z_G\cap \Gz$.
\begin{defi}\label{defi:dist}
An element $z\in \frakg_{\ubar d}$ is said to be \emph{distinguished} if $(Z_{\Gz}(z) / Z_{\ubar 0})^{\circ}$ is a unipotent group. A $\Gz$-orbit in $\frakg_{\ubar d}$ is called distinguished if every element of it is distinguished.
\end{defi}
Equivalently, an element $z\in \frakg_{\ubar d}$ is distinguished if the inclusion $Z_{\ubar 0}\subseteq Z_{\Gz}(z)$ induces a bijection $\bfX_*(Z_{\ubar 0})\xrightarrow{\sim}\bfX_*(Z_{\Gz}(z))$.

\begin{prop}\label{prop:rang} 
	Let $z\in \frakg_{\ubar d}^{\nil}$ be a distinguished nilpotent element. Then, given any $\frakl_*\in \frakG$ such that $z\in \frakl_d$, the spiral $\frakp_*$ attached to $z$ (see \autoref{subsec:spiral-sl2}) is the unique $\epsilon$-spiral containing $\frakl_*$ as splitting factor with $\epsilon = d / |d|$.
\end{prop}
\begin{proof}
	Suppose that $\frakl_*$ is a splitting factor of spiral of $\frakg_*$ such that $z\in \frakl_d$. Then, we may complete $z$ into an $\fraksl_2$-triplet $(z, h, f)$ such that $h\in \frakl_0$ and $f\in \frakl_{-d}$. Let $\iota\in \bfX_*(L_0)$ be the associated cocharacter characterised by $\rmd\iota(1) = h$. It follows that
	\[
		\frakl'_* = \bigoplus_{n\in \bbZ} \frakl'_n,\quad \frakl'_n = {^{\iota}_{2n/d}}\frakg_{\ubar n}
	\]
	is a splitting factor of the spiral $\frakp_*$ and it contains $\frakl_*$ due to \autoref{prop:cuspidal-lusztig}. In particular, $L$ is a Levi subgroup of $L'$. On the other hand, the distinguishedness of $e$ implies that $Z_L / Z_{\ubar 0}$ is finite; hence $L = L'$ holds. Consequently, $\frakl_*$ is a splitting factor of $\frakp_*$. Fix a maximal torus $S\subseteq L_0$ such that $\image \iota\subseteq S$. Let $\frakp'_*$ be any $\epsilon$-spiral having $\frakl_*$ as splitting factor with $\epsilon = d / |d|$. Since $S\subseteq L_0\subseteq P'_0$, there exists $\iota'\in \bfX_*(S)$ and $q\in \bfZ_{>0}$ such that 
	\[
		\frakp'_n = {^{\epsilon\iota'}_{\ge \epsilon qn}}\frakg_{\ubar n}\quad \text{and}\quad \frakl_n = {^{\iota'}_{qn}}\frakg_{\ubar n}\quad \text{hold for $n\in \bbZ$}. 
	\]
	In particular, the cocharacter $\lambda\in \bfX_*(S)$ defined by $\lambda(t) = \iota'(t)^2\iota(t)^{-qd}$ centralises $\frakl$ and thus lies in $\bfX_*(Z_L)$. The finiteness of $Z_L / Z_{\ubar 0}$ implies that $\lambda\in \bfX_*(Z_{\ubar 0})$. It follows that $\frakp'_* = \frakp_*$ holds, whence the unicity.
\end{proof}

\begin{prop}\label{prop:dist-induction}
Suppose that $\xi = (L, \frakl_*, \scrC)\in \frakT$ is an LY-cuspidal support and that $\rmO = \Ad_{\Gz}\ci\frakl_d$ is distinguished. Then, there is a canonical isomorphism $\Ind^{\frakg_{\ubar d}}_{\frakp_d}\IC(\scrC) \cong \IC(\scrS)$, where $\frakp_*$ is the spiral attached to any element $e\in \ci\frakl_{d}$. Moreover, the $\xi$-series is a singleton: $\Pi_{\xi} = \{(\rmO, \scrS)\}$.
\end{prop}
\begin{proof}
	\autoref{lemm:spiral-stab} implies that $\Ind^{\frakg_{\ubar d}}_{\frakp_d}\IC(\scrC) \cong \IC(\scrS)\oplus \scrK$ for some complex $\scrK\in \Db_{\Gz}(\frakg^{\nil}_{\ubar d})_{\xi}$ concentrated in $\ba\rmO \setminus \rmO$. \autoref{lemm:supp} implies that $\scrK = 0$. The last statement follows from \autoref{prop:rang}.
\end{proof}

\subsection{Parabolic induction and spiral induction}\label{subsec:generic-induction}
Let $\xi = (L, \frakl_*, \scrC)\in \frakT$ be an LY-cuspidal support on $\frakg_{\ubar d}$. Put $M = Z_G(Z^{\circ}_L)$. Then, $M$ is a $\theta$-isotropic Levi subgroup by \autoref{coro:theta-isotropic}. Let $\Pi_M$ denote the counterpart of the set of pairs $\Pi$ on $\frakm_{\ubar d}$. We may regard $\xi$ as an LY-cuspidal support on $\frakm_{\ubar d}$. Let $(\rmO, \scrS)\in \Pi_{M, \xi}$ be the primitive pair on $\frakm_{\ubar d}$ attached to $\xi$. 
We shall relate the parabolic induction of $\IC(\scrS)$ and special cases of spiral induction of $\IC(\scrC)$. \par
\begin{lemm}
	The $M_{\ubar 0}$-orbit $\rmO\subseteq \frakm^{\nil}_{\ubar d}$ is distinguished.
\end{lemm}
\begin{proof}
	Let $z\in \ci\frakl_d$. By \autoref{prop:cuspidal-lusztig}, $(Z_{L_0}(z)/ Z_L)^{\circ}$ is unipotent. On the other hand, $Z_L / Z_M$ is finite and $Z_{L_0}(z)$ is a maximal reductive subgroup of $Z_{M_{\ubar 0}}(z)$ by \autoref{lemm:spiral-stab}. It follows that $(Z_{M_{\ubar 0}}(z) / Z_M \cap M_{\ubar 0})^{\circ}$ is unipotent, so $\rmO$ is distinguished.
\end{proof}
Set $\Lambda_{L} = \bfX_*(Z_L^{\circ})$. An element $\mu\in \Lambda_{L}$ is called \emph{regular} if $Z_G(\mu) = M$ holds.
\begin{defi}\label{defi:positive}
	Given a regular cocharacter $\mu\in \Lambda_{L}$, an $\epsilon$-spiral $\frakp_*\in \frakP^{\epsilon}$ with splitting factor $\frakl_*\in \frakG$ is called \emph{$\mu$-positive} if ${^\mu_{>0}}\frakg_{\ubar i}\subseteq \frakp_i$ and ${^\mu_{<0}}\frakg_{\ubar i}\cap \frakp_i = 0$ hold for $i\in \left\{ 0, d \right\}$. 
\end{defi}
Set $\epsilon = d / |d|$. Given any regular cocharacter $\mu\in \Lambda_{L}$, we may choose a $\mu$-positive $\epsilon$-spiral $\frakp_*$ of $\frakg_*$ with splitting factor $\frakl_*$ satisfying ${^\mu_{>0}}\frakg_{\ubar i}\subseteq \frakp_i$ for $i\in \left\{ 0, d \right\}$. Indeed, given $\lambda\in \bfX_*(\Gz)_{\bbQ}$ such that $\frakl_* = \frakl^{\lambda}_*$ (in the notation of \autoref{subsec:spiral}), then ${^\epsilon}\frakp^{\lambda_1}_{*}$ is $\mu$-positive for $\lambda_1 = \lambda + r\epsilon\mu$ whenever $r \in \bbZ$ and $r \gg 0$. 

\begin{prop}\label{prop:parab-spiral}
	Let $\mu\in \Lambda_{L}$ be a regular cocharacter and $\frakp_*\in \frakP^{\epsilon}$ a $\mu$-positive $\epsilon$-spiral. Set $\frakq = {^{\mu}_{\ge 0}}\frakg$ and $Q = \exp(\frakq)\subseteq G$. Then, $M$ is a Levi factor of $Q$ and there is a natural isomorphism $\Ind^{\frakg_{\ubar d}}_{\frakq_{\ubar d}}\IC(\scrS) \cong \Ind^{\frakg_{\ubar d}}_{\frakp_{d}}\IC(\scrC)$.
\end{prop}
\begin{proof}
	The regularity of $\mu$ implies that $M = Z_G(\mu)$ and thus $M$ is a Levi factor of $Q$. Let $\frakv = {^{\mu}_{> 0}}\frakg$ denote the nil-radical of $\frakq$. Let $\frakp^M_*$ be the spiral of $\frakm_*$ with splitting factor $\frakl_*$ defined by 
\[
	\frakp^M_n = \frakp_n\cap \frakm_{\ubar n}\quad\text{for $n\in \bbZ$}. 
\]
The $\mu$-positivity for $\frakp_*$ can be rephrased as $P_0 = P^M_0\cdot V_{\ubar 0}$ and $\frakp_d = \frakp^M_d \oplus \frakv_{\ubar d}$. The transitivity of induction \eqref{equa:trans-parab-spiral} from \autoref{subsec:parab-ind} reads:
\begin{equation}\label{equa:trans-spiral-parab}
	\Ind^{\frakg_{\ubar d}}_{\frakp_{d}}\cong\Ind^{\frakg_{\ubar d}}_{\frakq_{\ubar d}}\circ \Ind^{\frakm_{\ubar d}}_{\frakp^M_{d}}.
\end{equation}
By \autoref{prop:rang}, $\frakp^M_*$ coincides with the spiral of $\frakm_*$ attached to some element of $\ci\frakl_d\subseteq\rmO$, so \autoref{subsec:spiral-sl2} implies $\IC(\scrS) \cong \Ind^{\frakm_{\ubar d}}_{\frakp^M_{d}}\IC(\scrC)$, which together with \eqref{equa:trans-spiral-parab} implies the statement.
\end{proof}

\subsection{Parabolic induction of distinguished primitive pairs}\label{subsec:parab-prim}
We keep the setup of \autoref{subsec:generic-induction}. Given any $\theta$-stable parabolic subgroup $Q\subseteq G$ with Levi decomposition $Q = MV$, denote by $q:\frakq_{\ubar d}\to  \frakm_{\ubar d}$ the projection; then, the fibre of the induced complex $\Ind^{\frakg_{\ubar d}}_{\frakq_{\ubar d}} \IC(\scrS)$ at $z\in \frakg_{\ubar d}^{\nil}$ is given by
\begin{equation}\label{eq:fib-ind}
	(\Ind^{\frakg_{\ubar d}}_{\frakq_{\ubar d}} \IC(\scrS))_z = (\Ind^{\frakg_{\ubar d}}_{\frakq_{\ubar d}} j_{\rmO!}(\scrS))_z[\dim \rmO] \cong \rmR\Gamma_c(\Sp^Q_z, \dot\scrS)[\dim \rmO + \dim \frakv_{\ubar d} + \dim \frakv_{\ubar 0}],
\end{equation}
where $\frakv = \Lie V$ is the Lie algebra,
\[
	\Sp^Q_z = \left\{ gQ_{\ubar 0}\in G_{\ubar 0} / Q_{\ubar 0}\;\vert\;\Ad_{g^{-1}}z\in q^{-1}(\rmO)\right\}\xrightarrow{p}\rmO, \quad gQ_{\ubar 0}\mapsto q(\Ad_{g^{-1}}z)
\]
and $\dot\scrS = p^*\scrS$ is the inverse image of $\scrS$, which is a local system on $\Sp^Q_z$. We view $\rmR\Gamma_c(\Sp^Q_z, \dot\scrS)$ as an object of the equivariant category $\Db_{Z_{\Gz}(z)}(\pt)$, so that its cohomology is a module over the component group $\pi_0(Z_{\Gz}(z))$ in each degree. This can be view as an analogue of the cohomology of the Springer fibre with coefficients in a cuspidal local system, {\itshape cf.} \cite[\S 6]{lusztig84} and \cite[\S 8]{lusztig88}. The following proposition generalises \cite[\S 5.2]{LYIV}:
\begin{prop}\label{prop:primitive-induction}
	Given a pair $\pi = (\rmC, \scrL)\in \Pi$, the following conditions are equivalent:
	\begin{enumerate}
		\item\label{prop:primitive-induction-i}
			$\pi$ lies in the $\xi$-block $\Pi_{\xi}$;
		\item\label{prop:primitive-induction-ii}
			there exists a $\theta$-stable parabolic subgroup $Q\subseteq G$ containing $M$ as Levi factor such that
			\[
				\Hom_{\pi_0(Z_{\Gz}(z))}(\scrL_z, \rmH^{\bullet}_{c}(\Sp^Q_z, \dot\scrS))\neq 0\quad \text{for $z\in \rmC$};
			\]
		\item\label{prop:primitive-induction-iii}
			there exists a $\theta$-stable parabolic subgroup $Q\subseteq G$ containing $M$ as Levi factor such that
			\[
				\Hom_{\Gz}^{\bullet}(j_{\rmC!}\scrL, \Ind^{\frakg_{\ubar d}}_{\frakq_{\ubar d}}\IC(\scrS))\neq 0.
			\]
	\end{enumerate}
\end{prop}
\begin{proof}
	We prove \ref{prop:primitive-induction-ii} $\Rightarrow$ \ref{prop:primitive-induction-iii}. Making use of \eqref{eq:fib-ind}, we deduce
	\[
		\Hom_{\Gz}^{\bullet}(j_{\rmC!}\scrL, \Ind^{\frakg_{\ubar d}}_{\frakq_{\ubar d}}\IC(\scrS)) \cong \Hom^{\bullet + n}_{Z_{\Gz}(z)}(\scrL_z, \rmR\Gamma_{c}(\Sp^Q_z, \dot\scrS))\quad \text{for certain } n \in \bbZ. 
	\]
	The vanishing of $\rmH^{\bullet}_{c}(\Sp^Q_z, \dot\scrS)$ in odd degrees \cite[14.10(c)]{LYII} implies that the spectral sequence
	\[
		E^{p,q}_2 = \Ext^p_{Z_{\Gz}(z)}(\scrL_z, \rmH^{q}_{c}(\Sp^Q_z, \dot\scrS)) \Rightarrow \Hom^{p+q}_{Z_{\Gz}(z)}(\scrL_z, \rmR\Gamma_{c}(\Sp^Q_z, \dot\scrS))
	\]
	degenerates at $E_2$; in particular, the non-vanishing of $E^{0,q}_2$ implies that of $\Hom^{q}_{Z_{\Gz}(z)}(\scrL_z, \rmR\Gamma_{c}(\Sp^Q_z, \dot\scrS))$ for $q\in \bbZ$. 
	\par
	The implication \ref{prop:primitive-induction-iii} $\Rightarrow$ \ref{prop:primitive-induction-i} follows from \autoref{prop:parab-spiral}, which is applicable because $\mu$-positive $\epsilon$-spirals exist for every regular cocharacter $\mu\in \Lambda_{L}$. \par
	We prove \ref{prop:primitive-induction-i} $\Rightarrow$ \ref{prop:primitive-induction-ii} following the strategy of \cite[\S 7.2]{lusztig95b} and \cite[\S 5.2]{LYIV}. Assume that $\pi\in \Pi_{\xi}$. Choose an $\fraksl_2$-triplet $(z, h, f)$ with $z\in \rmC$, $h\in \frakg_{\ubar 0}$ and $f\in \frakg_{-\ubar d}$ and let $\iota\in \bfX_*(\Gz)$ be such that $\rmd\iota(1) = h$. Let $\frakp^z\in \frakP^{\epsilon}$ (with $\epsilon = d / |d|$) and $\frakh_*\in \frakG$ be the spiral and splitting factor attached to $(z, h, f)$ defined in \autoref{subsec:spiral-sl2}:
	\[
		\frakp^z_{n} = {^{\iota}_{\ge 2n/d}}\frakg_{\ubar{n}},\quad \frakh_{n} = {^{\iota}_{2n/d}}\frakg_{\ubar{n}},
	\]
	so that $z\in \frakh_d$ holds in particular. Let $\frakT_H$ and $\Pi_H$ denote the counterparts of $\frakT$ and $\Pi$ for $H$, see \autoref{rema:Z-graded}. \par
	 By \autoref{lemm:spiral-stab}, the inclusion $Z_{H_0}(z)\subseteq Z_{G_{\ubar 0}}(z)$ induces an isomorphism on component groups: $\pi_0(Z_{H_0}(z))\cong \pi_0(Z_{G_{\ubar 0}}(z))$, the restriction $\scrL\mid_{\rmC\cap \frakh_d}$ is an irreducible $H_0$-equivariant local system and the pair $\pi':= (\rmC \cap \frakh_d, \scrL\mid_{\rmO \cap \frakh_d})\in \Pi_H$ satisfies 
	 \[
		 \Ind_{\frakp^z_d}^{\frakg_{\ubar d}}\IC_{\pi'} \cong \IC_{\pi} \oplus \scrK
	 \]
	 for some $\scrK\in \Db(\frakg^{\nil}_{\ubar d})$ (see the arguments of \cite[\S 7.1]{LYI}). As $\IC_\pi$ lies in $\Db_{\Gz}(\frakg_{\ubar d}^{\nil})_{\xi}$, the orthogonal decomposition \eqref{eq:dco} from \autoref{subsec:dco} implies that $\pi'$ lies in the $\xi'$-block for some $\xi' = (L', \frakl'_*, \scrC')\in \frakT_{H}$ which is $\Gz$-conjugate to $\xi$ and satisfies $\frakl'_*\subseteq \frakh_*$; therefore, up to replacing $(z, h, f)$ and $\iota$ with a $\Gz$-conjugate, we may assume that $\frakl_*\subseteq \frakh_*$ holds and $\pi'$ lies in the $\xi$-block. Hence, there exists a $\bbZ$-graded parabolic subalgebra $\frakp_*\subseteq \frakh_*$ which contains $\frakl_*$ as $\bbZ$-graded Levi subalgebra such that the $\scrL_z$-isotypic component of the cohomology $\rmH^{\bullet}_c(\Sp^P_z, \dot\scrC)$ is non-zero, where
	\[
		\Sp^P_z = \left\{ gP_{0}\in H_{0} / P_{0}\;\vert\;\Ad_{g^{-1}}z\in (q^1)^{-1}(\rmO)\right\}\xrightarrow{p^1}\ci\frakl_d, \quad gP_{0}\mapsto q^1(\Ad_{g^{-1}}z),
	\]
	$q^1: \frakp_d\to \frakl_d$ is the projection and $\dot\scrC$ is the inverse image $(p^1)^*\scrC$. Up to replacing $P$ with $aPa^{-1}$ for some $a\in H_0$, we may assume that $z\in \frakp_d$. \par

	The $\theta$-isotropic Levi subgroup $M = Z_G(Z_L^{\circ})$ satisfies $M\cap H = L$. Pick any $\theta$-stable parabolic group $Q\subseteq G$ which contains $M$ as Levi factor and satisfies $Q \cap H = P$. Taking $\iota$-fixed points of $\Sp^Q_z$ yields an identity of virtual representations:
	\begin{equation}\label{eq:groth}
		\sum_{n\in \bbZ} (-1)^n[\rmH_c^n(\Sp^{Q}_z, \dot\scrS)] = \sum_{n\in \bbZ} (-1)^n[\rmH_c^n((\Sp^{Q}_z)^{\iota}, \dot\scrS )].
	\end{equation}
	\par
	The connected components of $(\Sp^{Q}_z)^{\iota}$ can be described as follows: let $g_1Q_{\ubar 0}, \cdots, g_rQ_{\ubar 0}$ be a complete list of representatives of $H_{\ubar 0}$-orbits in $(\Gz / Q_{\ubar 0})^{\iota}$, where $1 = g_1, g_2, \cdots, g_r\in \Gz$. For each $i\in \left\{ 1, \cdots, r \right\}$, set $P^i = g_iQg_i^{-1}\cap H$, $L^i = g_iMg_i^{-1} \cap H$. The Lie algebras $\frakp^i = \Lie P^i$ and $\frakl^i = \Lie L^i$ are equipped with the induced $\bbZ$-gradings from $\frakh_*$. Let $q^i:\frakp^i_d\to \frakl^i_d$ denote the projection. The following diagram is commutative:
	\[
		\begin{tikzcd}[column sep=2cm]
			\frakq_{\ubar d}\arrow{d}{q}\arrow{r}{\Ad_{g_i}} & \Ad_{g_i}\frakq_{\ubar d} \arrow{d}{\Ad_{g_i}\circ q\circ \Ad_{g_i}^{-1}}\arrow[hookleftarrow]{r} & \frakp^i_d \arrow{d}{q^i} \\
			\frakm_{\ubar d}\arrow{r}{\Ad_{g_i}} & \Ad_{g_i}\frakm_{\ubar d}\arrow[hookleftarrow]{r} & \frakl^i_d
		\end{tikzcd}
	\]
	We have $(\Sp^Q_z)^{\iota} = \bigsqcup_i \Sp^{P^i}_z$, where
	\[
		\Sp^{P^i}_z = \left\{ aP^i_{0}\in H_{0} / P^i_{0}\;\vert\;\Ad^{-1}_{a}z\in \Ad_{g_i}q^{-1}(\rmO)\right\}\hookrightarrow \Sp^Q_{e},\quad aP^i_{0}\mapsto ag_iP_0.
	\]
	Up to renumbering, we may assume that $\Sp^{P^i}_z$ is non-empty for $1 \le i \le r'$ and empty for $r'< i \le r$. For each $i\in \left\{ 1, \cdots, r' \right\}$, up to replacing $g_i$ with $ag_i$ for some $a\in H_0$ satisfying $aP^i_0\in \Sp^{P^i}_z$, we may assume furthermore that $P^i\in \Sp^{P^i}_z$.
	\par
We show that $\Ad_{g_i}\rmO\cap \frakh_d$ coincides with the open $L^i_0$-orbit $\ci\frakl^i_d\subseteq \frakl^i_d$ for $i\in \left\{ 1, \cdots, r' \right\}$.
	The intersection is non-empty --- indeed, $\Ad_{g_i}^{-1}z$ lies in $q^{-1}(\rmO)$ and hence $z' :=\Ad_{g_i}q(\Ad_{g_i}^{-1}z) = q^i(z)$ lies in $\Ad_{g_i}\rmO\cap \frakh_d$. As $\rmO$ is a distinguished $M_{\ubar 0}$-orbit, so is the conjugate $\Ad_{g_i}\rmO$ a distinguished $g_iM_{\ubar 0}g_i^{-1}$-orbit and thus each element of $\Ad_{g_i}\rmO\cap \frakh_d$ is in a distinguished $L^i_0$-orbit. The openness of distinguished orbits in the $\bbZ$-graded setting \cite[4.3]{lusztig95b} implies that $\Ad_{g_i}\rmO\cap \frakh_d$ coincides with the unique open $L^i_0$-orbit $\ci\frakl^i_d$.  \par
	We show that the adjoint action of $g_i$ induces a $\bbZ$-graded isomorphism $\frakl_* \xrightarrow{\sim} \frakl^i_*$ for $i\in \left\{ 1, \cdots, r' \right\}$. Let $(z', h', f')$ be an $\fraksl_2$-triplet in $\frakl$ satisfying $z'\in \ci\frakl_d$, $h'\in \frakl_0$ and $f'\in \frakl_{-d}$ and $\iota'\in \bfX_*(L_0)$ such that $\rmd\iota'(1) = h'$. Up to replacing $(z', h', f', \iota')$ with an $L_0$-conjugate, we may assume that $\iota'\in \bfX_*(S)$. Then, $\Ad_{g_i}z'\in \Ad_{g_i}\rmO\cap \frakh_d =\ci\frakl^i_d$ holds. On the other hand, we have $\Ad_{g_i}z'\in \Ad_{g_i}{^{\iota'}_2}\frakl = {^{\Ad_{g_i}\iota'}_2}\frakl^i$. We may write $g_i = ab$ with $a\in H_0$ and $b\in N_{\Gz}(S)$, so that the image of $\Ad_{b}\iota'\in \bfX_*(S)$ commutes with that of $\Ad^{-1}_{a}\iota = \iota$. Hence, we can define a cocharacter $\iota''\in \bfX_*(L^i_0)$ by $\iota''(t) = \iota(t)\Ad_{g_i}\iota'(t^{-1})$. Now, $\iota''$ centralises the element $z'$, which lies in the distinguished $L^i_0$-orbit $\ci\frakl^i_d$; the distinguishedness implies that the image of $\iota''$ lies in $Z_{L^i}$ and hence $\frakl^i_n = {^{\Ad_{g_i}\iota'}_{2n}}\frakl^i_n$ holds for $n\in \bbZ$, which implies the claim. \par

	For $i\in \left\{ 1, \cdots, r' \right\}$, the adjoint action of $g_i$ transports the cuspidal pair $(\ci\frakl_d, \scrC)$ on $\frakl_d$ to a cuspidal pair $(\ci\frakl^i_d, \scrC^i)$ on $\frakl^i_d$ and $(\rmO, \scrS)$ to the distinguished primitive pair $(\rmO^i, \scrS^i)$ on $\Ad_{g_i}\frakm_{\ubar d}$ attached to the LY-cuspidal support $(L^i, \frakl^i_*, \scrC^i)$.
Consider
\[
	p^i:\Sp^{P^i}_z\to \ci\frakl^i_d,\quad aP^i_0\mapsto q^i(\Ad_{a^{-1}}z),\quad \text{where $q^i:\frakp^i_d\to \frakl^i_d$ is the projection}.
\]
The restriction $\dot\scrS\mid_{\Sp^{P^i}_z}$ is isomorphic to $\dot\scrC^i:=(p^{i})^*\scrC^i$. It follows that we have a decomposition:
\[
	\rmH^{\bullet}_c((\Sp^{Q}_z)^{\iota}, \dot\scrS ) \cong \bigoplus_{i=1}^{r'}\rmH^{\bullet}_c(\Sp^{P_i}_z, \dot\scrC^i) = \rmH^{\bullet}_c(\Sp^{P}_z, \dot\scrC)\oplus \bigoplus_{i=2}^{r'}\rmH^{\bullet}_c(\Sp^{P_i}_z, \dot\scrC^i).
\]
It is known by \cite[7.6]{lusztig95b} and \cite[14.10(c)]{LYII} that the cohomology groups $\rmH^{\bullet}_c(\Sp^{Q}_z, \dot\scrS )$ and $\rmH^{\bullet}_c(\Sp^{P_i}_z, \dot\scrC^i)$ for each $i\in \left\{ 1, \cdots, r' \right\}$ vanish in odd degrees; this parity condition and the equation \eqref{eq:groth} imply the existence of an isomorphism in $\bbC\pi_0(Z_{\Gz}(z))\mof$:
\[
	\bigoplus_{n\in \bbZ}\rmH_c^n(\Sp^{Q}_z, \dot\scrS)\cong \bigoplus_{n\in \bbZ}\rmH^{n}_c(\Sp^{P}_z, \dot\scrC)\oplus \bigoplus_{i=2}^{r'}\bigoplus_{n\in \bbZ}\rmH^{n}_c(\Sp^{P_i}_z, \dot\scrC^i).
\]
Since the right-hand side contains $\scrL_z$ as direct summand, so does the left-hand side. 
\end{proof}
\begin{rema}
	The spiral version of \autoref{prop:primitive-induction}, namely that every pair $\pi\in \Pi$ appears in the fibres of the cohomology of the induced complex $\Ind^{\frakg_{\ubar d}}_{\frakp_d}\IC(\scrC)$ for some $\epsilon$-spiral $\frakp_*\in \frakP^{\epsilon}_{\xi}$ with $\epsilon = d / |d|$, is an immediate consequence of \cite[13.8(a)]{LYII} and the orthogonal decomposition \eqref{eq:dco}. The conjunction of \autoref{prop:primitive-induction} and \autoref{prop:parab-spiral} shows that $\frakp_*$ can be chosen to be $\mu$-positive for some regular cocharacter $\mu\in \Lambda_{L}$ which depends on $\pi$. This statement can be expected in view of the representation theory of degenerate double affine Hecke algebras --- namely, non-zero submodules of \emph{standard modules} are not annihilated by the Knizhnik--Zamolodchikov functors.
\end{rema}

\section{Supercuspidal pairs and supercuspidal supports}\label{sec:supercuspidal}

We keep the setup of \autoref{subsec:G}. In this section, we introduce the notions of supercuspidal pairs in $\Pi$ and supercuspidal supports on $\frakg_{\ubar d}$. The objective is to prove \autoref{theo:supercuspidal}, which provides useful characterisations of supercuspidal pairs. It turns out that supercuspidal pairs are precisely the distinguished primitive pairs studied in \autoref{subsec:dist-prim}-\autoref{subsec:parab-prim}.

\subsection{Supercuspidal pairs}\label{subsec:supercuspidal}
\begin{defi}\label{defi:supercuspidal}
	A pair  $\left(\rmO, \scrS \right)\in \Pi$ is called \emph{supercuspidal} if for every $\theta$-stable proper parabolic subalgebra $\frakq\subsetneq \frakg$ with nil-radical $\frakv$, the vanishing condition
	\[
		\rmH^\bullet_c\left( \rmO\cap \left(z + \frakv_{\ubar d} \right), \scrS\mid_{ \rmO\cap \left(z + \frakv_{\ubar d} \right)} \right) = 0
	\]
	holds for every $z\in \frakq_{\ubar d}$.
\end{defi}
Equivalently, $(\rmO, \scrS)$ is a supercuspidal pair if the extension $j_{\rmO !}\scrS\in \Db_{\Gz}(\frakg^{\nil}_{\ubar d})$ is annihilated by the parabolic restriction (\ref{subsec:parab-ind}) along every proper $\theta$-stable parabolic algebra of $\frakg$, where $j_{\rmO}:\rmO\hookrightarrow \frakg^{\nil}_{\ubar d}$ is the inclusion of orbit.

\begin{theo}\label{theo:supercuspidal} 
	Given a pair $\pi = \left( \rmC, \scrL \right)\in \Pi$, the following conditions are equivalent:
	\begin{enumerate}
		\item\label{theo:supercuspidal-i}
			The pair $\pi$ is supercuspidal.
		\item\label{theo:supercuspidal-iii}
			The $\Gz$-orbit $\rmC$ is distinguished (\autoref{subsec:dist-prim}) and the local system $\scrL$ is clean.
		\item\label{theo:supercuspidal-ii}
			There exists an LY-cuspidal support $\xi = (L,  \frakl_*, \scrC)\in \frakT$ such that $Z_L /(Z_G\cap \Gz)$ is a finite group and $\pi$ is isomorphic to the primitive pair attached to $\xi$. 
	\end{enumerate}
\end{theo}
\begin{proof}[Proof of \autoref{theo:supercuspidal}]
	Let $\xi = (L, \frakl_*, \scrC)\in \frakT$ be an LY-cuspidal support such that $\pi\in \Pi_{\xi}$. 
	\par
	We prove \ref{theo:supercuspidal-ii} $\Rightarrow$ \ref{theo:supercuspidal-i}. Let $Q\subset G$ be any $\theta$-stable proper parabolic subgroup. There exists by \autoref{coro:theta-isotropic} a cocharacter $\lambda\in \bfX_*(\Gz)$ such that $\frakq = {^\lambda_{\ge 0}}\frakg$ and $M = Z_G(\lambda)$ is a $\theta$-stable Levi factor of $Q$. Let $j_{\rmC}: \rmC\to \frakg^{\nil}_{\ubar d}$ be the inclusion of orbit. It follows that $j_{\rmC !}\scrL$ lies in $\Db_{\Gz}\left( \frakg_{\ubar d}^{\nil} \right)_{\xi}$ by \autoref{lemm:ext-block} and thus $\Res^{\frakg_{\ubar d}}_{\frakq_{\ubar d}} j_{\rmC !}\scrL = 0$ by adjunction since no $\Gz$-conjugate of $L$ is contained in $M$ (otherwise, $L$ would be centralised by some $\Gz$-conjugate of $\lambda$, contradicting the finiteness of $Z_L/(Z_G\cap \Gz)$). It follows that $\pi$ is supercuspidal. \par
	We prove \ref{theo:supercuspidal-iii} $\Rightarrow$ \ref{theo:supercuspidal-ii}. Let $(\rmO, \scrS)\in \Pi_{\xi}$ be the primitive pair on $\frakg_{\ubar d}$ attached to $\xi$, which is clean by \autoref{prop:cleanness}. If $\rmC \neq \rmO$, then the cleanness of $\scrS$ and $\scrL$ implies that 
	\[
		\Hom^{\bullet}_{\Gz}(\IC_{\pi}, \IC(\scrS)) \cong \Hom^{\bullet + \dim \rmO - \dim \rmC}_{\Gz}(j_{\rmC!}\scrL, j_{\rmO*}\scrS)  = 0,
	\]
	which contradicts \autoref{lemm:HomFG}; hence $\rmC = \rmO$ holds and $\pi \cong (\rmO, \scrS)$ by \autoref{lemm:supp}. The distinguishness of the orbit $\rmC$ implies that $L$ is semisimple. 
		\par
		We prove \ref{theo:supercuspidal-i} $\Rightarrow$ \ref{theo:supercuspidal-iii}. Set $M = Z_G(Z_L^{\circ})$ and let $(\rmO, \scrS)\in \Pi_M$ be the primitive pair on $\frakm_{\ubar d}$ attached to $\xi$. By \autoref{prop:primitive-induction}\ref{prop:primitive-induction-i}$\Rightarrow$\ref{prop:primitive-induction-iii}, there exists a $\theta$-stable parabolic subgroup $Q\subseteq G$ which contains $M$ as Levi factor such that 
		\begin{equation}\label{equa:HomLS}
			\Hom^{\bullet}_{\Gz}(j_{\rmC!}\scrL, \Ind^{\frakg_{\ubar d}}_{\frakq_{\ubar d}}\IC(\scrS))\neq 0.
		\end{equation}
		The adjunction (\autoref{subsec:parab-ind}) implies that $\Res^{\frakg_{\ubar d}}_{\frakq_{\ubar d}}j_{\rmC!}\scrL \neq 0$. Then, the supercuspidality of $\pi$ implies that $M = Q = G$ holds. In particular, $\rmO$ is a distinguished $\Gz$-orbit. Consequently, \eqref{equa:HomLS} implies $j^!_{\rmC}\IC(\scrS)\neq 0$ and the cleanness of $(\rmO, \scrS)$ shown in \autoref{prop:cleanness} implies that $\rmC = \rmO$ implies. Then, \autoref{lemm:supp} implies that $\pi\cong (\rmO, \scrS)$ holds. Hence, $\rmC = \rmO$ is distinguished and $\pi$ is clean. 
 %Let $(P^M_0, \frakp^M_*)$ be any spiral on $\frakm_*$ (\autoref{subsec:spiral-levi}) with splitting $(L, L_0, \frakl_*)$ and let $(P_0, \frakp_*)$ be a spiral on $\frakg_*$ with splitting $(L, L_0, \frakl_*)$ which satisfy the condition of \autoref{subsec:trans}. Then $\Ind^{\frakg_{\ubar d}}_{\frakq_{\ubar d}}\Ind^{\frakm_{\ubar d}}_{\frakp^M_d}\IC(\scrC)\cong \Ind^{\frakg_{\ubar d}}_{\frakp_{d}}\IC(\scrC)$ is a non-zero semisimple complex on $\frakg_{\ubar d}$. Suppose that $M\neq G$. The adjunction of parabolic induction-restriction (\autoref{subsec:parab-ind}) and the supercuspidality of $(\rmO, \scrS)$ yields $\Hom^{\bullet}_{\Gz}(\Res^{\frakg_{\ubar d}}_{\frakq_{\ubar d}}\IC(\scrS), \Ind^{\frakm_{\ubar d}}_{\frakp^M_d}\IC(\scrC)) = 0$, which contradicts \autoref{lemm:HomFG}. Therefore, $M = G$ holds. Let $(\rmO', \scrS')\in \Pi_{\Gz}(\frakg_{\ubar d})_{\xi}$ be the primitive pair, which is clean by \autoref{prop:cleanness}. Then, the implication \ref{theo:supercuspidal-ii} $\Rightarrow$ \ref{theo:supercuspidal-i} proven above implies that $(\rmO', \scrS')$ is supercuspidal and \autoref{prop:distingue2} shows that $\rmO'$ is distinguished. This implies that $\IC(\scrS')$ is the only simple object in $\Perv_{\Gz}(\frakg^{\nil}_{\ubar d})_{\xi}$ by \autoref{prop:dist-induction}. It follows that $(\rmO, \scrS)\cong (\rmO', \scrS')$, so $\rmO$ is distinguished and $\scrS$ is clean. 
\end{proof}
%As a consequence to~\autoref{theo:supercuspidal}, supercuspidal pairs are primitive pairs on distinguished nilpotent orbits, so \autoref{prop:dist-induction} apply to supercuspidal pairs.  

\subsection{Supercuspidal supports}\label{subsec:supp-supercuspidal}
\begin{defi}\label{defi:super-support}
	A \emph{supercuspidal support} on $\frakg_{\ubar d}$ is a triplet $(M, \rmO, \scrS)$, where $M\subseteq G$ is a $\theta$-isotropic Levi subgroup and $(\rmO, \scrS)\in \Pi_{M}$ is a supercuspidal pair on $\frakm_{\ubar d}$. 
\end{defi}
Two supports $(M, \rmO, \scrS)\cong (M', \rmO', \scrS')$ are said to be $\Gz$-conjugate if there exists a pair $(g, \eta)$ consisting of an element $g\in \Gz$ such that $gMg^{-1} = M'$, $\Ad_g\rmO = \rmO'$ and an isomorphism $\eta: g^* \scrS'\cong \scrS$. \par

Given an LY-cuspidal support $\xi = (L, \frakl_*, \scrC)$, put $M = Z_G(Z^{\circ}_{L})$ and let $(\rmO, \scrS)$ be the primitive pair on $\frakm_{\ubar d}$ attached to $\xi$. Then, \autoref{theo:supercuspidal} implies that $(\rmO, \scrS)$ is a supercuspidal pair on $\frakm_{\ubar d}$ and, conversely, every supercuspidal support arises this way. We call $\zeta = (M, \rmO, \scrS)$ the \emph{supercuspidal support attached to $\xi$}. \par

Let $\zeta = (M, \rmO, \scrS)$ be a supercuspidal support on $\frakg_{\ubar d}$.  In view of the cleanness (\autoref{theo:supercuspidal}) of $\scrS$, we have $j_!\scrS[\dim \rmO] \cong \IC(\scrC) \cong j_*\scrS[\dim \rmO]$, where $j:\rmO\hookrightarrow \frakm_d$ is the inclusion. 
\begin{defi}
	Given a supercuspidal support $\zeta = (M, \rmO, \scrS)$, we say that $\zeta$ is a \emph{supercuspidal support} of a pair $(\rmC, \scrL)\in \Pi$ if
	 there exists a $\theta$-stable parabolic subgroup $Q\subseteq G$ with Levi decomposition $Q = MV$ such that the following holds:
\begin{eq}
	\Hom_{\pi_0(Z_{\Gz}(z))}\left( \scrL_z, \rmH^{\bullet}_c(\Sp^Q_z,  \dot\scrS) \right)\neq 0\quad  \forall z\in \rmC,
\end{eq} 
where
\[
	\Sp^Q_z = \left\{ gQ_{\ubar 0}\in G_{\ubar 0} / Q_{\ubar 0}\;\vert\;\Ad_{g^{-1}}z\in q^{-1}(\rmO)\right\}\xrightarrow{p}\rmO, \quad gQ_{\ubar 0}\mapsto q(\Ad_{g^{-1}}z),
\]
$q:\frakq_{\ubar d}\to \frakm_{\ubar d}$ is the projection and $\dot\scrS = p^*\scrS$.
%\begin{eq}
	%\Hom_{Z_{\Gz}(z)}\left( \scrL_z, \rmH^{\bullet}(\Ind_{\frakq_{\ubar d}}^{\frakg_{\ubar d}}\IC\left( \scrS \right))_z \right)\neq 0\quad, \forall z\in \rmC.
%\end{eq} 
\end{defi}

\begin{theo}\label{theo:comparaison-bloc}
	Let $\xi\in \frakT$ be an LY-cuspidal support and let $\zeta$ be the supercuspidal support attached to $\xi$. Then, $\zeta$ is a supercuspidal support a pair $\pi\in \Pi$ if and only if $\pi$ lies in the $\xi$-block $\Pi_{\xi}$.
\end{theo}
\begin{proof}
	The statements follow immediately from \autoref{theo:supercuspidal} and \autoref{prop:primitive-induction}.
\end{proof}
In particular, the supercuspidal support of a pair $(\rmC, \scrL)\in \Pi$ is unique up to $\Gz$-conjugacy.

\begin{exam}[Principal block]\label{exam:principal1}
	Choose a pinning $E = (B, T, a)$ as in \autoref{subsec:grading} and set $S = (T^{\sigma_E})^{\circ}$. Consider the LY-cuspidal support $\xi_0 = (S, \fraks_*, \delta_0)$, where $\fraks_*$ is the trivial $\bbZ$-grading on $\fraks = \Lie S$ and $\delta_0$ is the constant sheaf with coefficients in $\bbC$ on $\fraks_d = 0$. The supercuspidal support attached to $\xi_0$ is $\zeta_0 = (T, 0, \delta_0)$. Given $z\in \frakg_{\ubar d}^{\nil}$, let 
	\[
		(G / B)^{\theta}_{z} = \{g B\in G/B\;\vert\; \Ad_{g^{-1}} z \in \frakb;\; \forall\omega\in \mu_m,\;g^{-1}\theta(\omega)(g)\in B\}.
	\]
	There is an isomorphism
	\[
		\bigsqcup_{gB\in (G / B)^{\theta}}\Sp^{gBg^{-1}}_z\xrightarrow{\sim} (G / B)^{\theta}_{z},\quad h(gBg^{-1})\mapsto hgB.
	\]
	Therefore, the simple constituents of $\rmH^{\bullet}((G / B)^{\theta}_{z}, \bbC)$ for various $z$ form the principal series $\Pi_{\xi}$. They recover the set of Deligne--Langlands--Lusztig parameters studied in \cite{KL}. 
\end{exam}

\section{Characterisation of bi-orbital sheaves}\label{sec:biorbital} 
In this section, we study bi-orbital sheaves on a graded Lie algebra. The main theorem, which characterises the bi-orbital sheaves as parabolically induced from bi-orbital supercuspidal sheaves, is stated in \autoref{subsec:statement} and proven in \autoref{subsec:proof}. We shall apply the generalised Springer correspondence to relate sheaves with modules of degenerate double affine Hecke algebras and then characterise bi-orbital sheaves via the algebraic Knizhnik--Zamolodchikov functor. 
\subsection{Main theorem}\label{subsec:statement}
Keep the setup of \autoref{subsec:G}. We shall make use of the Fourier--Sato transform on $\frakg_{\ubar d}$, which is an involutive equivalence on $\Cq$-monodromic complexes (\autoref{subsec:hyperbolic}):
\[
	\FS:\Db_{\Gz}(\frakg_{\ubar n})^{\Cq\operatorname{-mon}}\xrightarrow{\sim}\Db_{\Gz}(\frakg_{-\ubar n})^{\Cq\operatorname{-mon}}
\]
for each $\ubar n\in \bbZ / m$, see \cite[\S 3.7]{KS90}.
\begin{defi}
	A simple perverse sheaf $\scrF\in \Irr\Perv_{\Gz}(\frakg_{\ubar d})$ is called \emph{orbital} if $\scrF$ is concentrated in $\frakg_{\ubar d}^{\nil}$; it is called \emph{anti-orbital} if its Fourier--Sato transform $\FS(\scrF)$ is concentrated in $\frakg_{-\ubar d}^{\nil}$; it is called \emph{bi-orbital} if it is both orbital and anti-orbital. 
\end{defi}
Since $\supp \FS(\scrF)$ for a monodromic sheaf $\scrF$ equals the image of the singular support $\operatorname{SS}(\scrF)$ under the projection $\frakg_{\ubar d}\times \frakg_{-\ubar d}\to \frakg_{-\ubar d}$ the bi-orbital condition for $\scrF$ is equivalent to that $\operatorname{SS}(\scrF)$ lies in $\frakg^{\nil}_{\ubar d}\times \frakg^{\nil}_{-\ubar d}$. Let $\Pi^\nil \subseteq \Pi$ denote the subset of pairs $\pi\in \Pi$ such that $\IC_{\pi}$ is bi-orbital.
The main theorem is the following:
\begin{theo}\label{theo:biorbital}
	Let $\pi\in \Pi$. Then, $\pi\in \Pi^\nil$ holds if and only if there exist a supercuspidal support $(M, \rmO, \scrS)$ and a $\theta$-stable parabolic subgroup $Q\subseteq G$ containing $M$ as Levi factor such that $(\rmO, \scrS)\in \Pi^\nil_M$ holds and $\IC_\pi$ is a direct summand of the perverse cohomology of $\Ind^{\frakg_{\ubar d}}_{\frakq_{\ubar d}}\IC(\scrS)$. 
\end{theo}
The proof of \autoref{theo:biorbital} will occupy the rest of the section. 
\begin{rema}
	Let $\xi\in \frakT$ with supercuspidal pair $(\rmO, \scrS)\in\Pi_{M, \xi}$ attached to $\xi$. Note that even if $\IC(\scrS)$ is bi-orbital, there can be non-bi-orbital sheaves in the same block, namely that $\Pi^{\nil}_{\xi}\subset \Pi_{\xi}$ may not hold. \autoref{theo:biorbital} means that bi-orbital sheaves are exactly the simple constituents of parabolic induction from bi-orbital supercuspidal sheaves. 
\end{rema}

\subsection{Relative finite root system}\label{subsec:rel-rs}
Fix $\xi\in \frakT$. We describe the finite root system which underlies the relative affine root system $(\fraka_{\xi}, R^{\aff}_{\xi})$ introduced in \autoref{subsec:rars}. Let $\zeta = (M, \rmO, \scrS)$ be the supercuspidal support attached to $\xi$. In particular, $M = Z_G(Z^{\circ}_L)$ is the smallest Levi subgroup of $G$ containing $L$ and $(\rmO, \scrS)$ is the primitive pair attached to $\xi$. %\hl{Modification n\'ecessaire pour le cas tordu!}
Set $\frake_{\zeta,\bbQ} = \bfX_*(Z_L^{\circ})_{\bbQ}$ and $\frake_{\zeta} = \frake_{\zeta, \bbQ}\otimes \bbR$. Then, $\fraka_{\xi}$ is a principal homogeneous $\frake_{\zeta}$-space under the action of translations. For each relative affine root $\alpha\in R^{\aff}_{\xi}$, let $\dot\alpha$ be its derivative, considered as element of the dual vector space $\frake_{\zeta}^*$ of $\frake_{\zeta}$. Define
\[
	R_{\zeta} = \left\{\dot\alpha \in \frake_{\zeta}^*\;\vert\; \alpha\in R^{\aff}_{\xi} \right\}.
\]
Then, $(\frake_{\zeta}, R_{\zeta})$ is the underlying finite root system of $(\fraka_{\xi}, R^{\aff}_{\xi})$. Let $W_{\zeta}$ denote its Weyl group. 
The elements of $W^{\aff}_{\xi}$ which act as translations on $\fraka_{\xi}$ form a free abelian group $\Lambda_{\xi}$. It can be shown that $\Lambda_{\xi}$ is a subgroup of $\bfX_*(Z^{\circ}_{L})$. There is a canonical isomorphism $W_{\xi}^{\aff} / \Lambda_{\xi} \cong  W_{\zeta}$.
Under the quotient map $p:W_{\xi}^{\aff}\to W_{\zeta}$, the stabiliser $W^{\aff}_{\xi, \bfx}$ is mapped injectively into $W_{\zeta}$. Then, \autoref{lemm:weyl-stab} implies that 
\begin{equation}\label{equa:WG0}
	p(W^{\aff}_{\xi, \bfx})\cong N_{\Gz}(Z^{\circ}_{L}) / M_{\ubar 0}. 
\end{equation}\par
The following is immediate from \autoref{lemm:weyl-stab}:
\begin{lemm}\label{lemm:Wzeta}
	The conjugacy action of $N_{\Gz}(Z^{\circ}_L)$ on $\frake_{\zeta}$ induces an isomorphism $N_{\Gz}(Z^{\circ}_L) / M_{\ubar 0}\cong p(W^{\aff}_{\xi, \bfx})$. \hfill\qedsymbol
\end{lemm}

Let $\frakQ_{\zeta}$ denote the set of $\theta$-stable parabolic subalgebra $\frakq$ of $\frakg$ which contains $\frakm$ as Levi factor, indexed by the set of Weyl chambers of the root system $(\fraka_{\xi}, R_{\xi})$. The group $p(W^{\aff}_{\xi, \bfx})$ acts $\frake_{\zeta}$ and induces an action on $\frakQ_{\zeta}$ via \eqref{equa:WG0}. Let $\ubar \frakQ_{\zeta}$ denote the set of $p(W^{\aff}_{\xi, \bfx})$-orbits, or equivalently (due to \autoref{lemm:Wzeta}), the set of $\Gz$-conjugacy classes in $\frakQ_{\zeta}$.  %Then, $W_{\xi}$ acts on $\frakQ_{\xi}$. 

\subsection{Sheaf-theoretic realisation of \texorpdfstring{$\bfV$}{V}}\label{subsec:Vfaisc}
Assume for the moment that $G$ is simply connected and almost simple. Set $\epsilon = d / |d|$. Recall the dDAHA $\bfH_{\xi}$ attached to $\xi$ from \autoref{subsec:springer}. The sheaf-theoretic realisation of the spectral completion $\calH_{\xi} = \calC^{\bfx}\bfH_{\xi}$ provided by \autoref{theo:springer} allows us to express the algebraic KZ functor $\bfV$ (see \autoref{subsec:AKZ}) in sheaf-theoretic terms. Set
\[
	P_{\bfV} = \bigoplus_{\frakq\in \ubar \frakQ_{\zeta}} \;\bigoplus_{\frakp_*\in \ubar \frakP^{\epsilon}_{\xi}}\Hom^{\bullet}_{\Gz}(\Ind^{\frakg_{\ubar d}}_{\frakq_{\ubar d}}\IC(\scrS), \Ind^{\frakg_{\ubar d}}_{\frakp_{d}} \IC(\scrC))^{\wedge}_0.
\]
Via \autoref{theo:springer}, $P_{\bfV}$ acquires an $\calH_{\xi}$-module structure. Recall the idempotent $1_{\bfV}\in \calH_{\xi}$ introduced in \autoref{subsec:AKZ}.
\begin{prop}\label{prop:KZ-sheaf}
	Assume that $G$ is simply connected and almost simple. Then, there is a canonical isomorphism of $\calH_{\xi}$-modules $P_{\bfV}\cong \calH_{\xi} 1_{\bfV}$. In other words, $P_{\bfV}$ is a compact projective $\calH_{\xi}$-module and the algebraic KZ functor $\bfV$ is isomorphic to $\Hom_{\calH_{\xi}}(P_{\bfV}, \relbar)$. 
\end{prop}
\begin{proof}
	We view $\Lambda_{\xi}$ as subgroup of $\bfX_*(Z_L^{\circ})$. Fix $\alpha_0\in \Delta^{\aff}$ such that $\Delta = \{\dot\alpha\;\vert\; \alpha\in \Delta^{\aff}\setminus\left\{ \alpha_0 \right\}\}$ is a basis for $(\frake_{\zeta}, R_{\zeta})$. Let $\gamma\in \Lambda_{\xi}$ be an element such that $\langle \gamma, \alpha\rangle \gg 0$ for every $\alpha\in \Delta$. Write $t_{-\gamma}\in W^{\aff}_{\xi}$ for the image of $-\gamma$ under the inclusion $\Lambda_{\xi}\subseteq  W_{\xi}^{\aff}$ described in \autoref{subsec:rel-rs}. Given $w\in W_{\zeta}$, the $\epsilon$-spiral ${^{\tilde w}}\frakp_*$ with $\tilde w = w t_{-\gamma}$ is $(-w\gamma)$-positive in the sense of \autoref{defi:positive}. Consequently, setting ${^w}\frakq = {^{w\gamma}_{\le 0}}\frakg$, we have by \autoref{prop:parab-spiral}
\[
	\Ind^{\frakg_{\ubar d}}_{{^{\tilde w}}\frakp_{d}} \IC(\scrC) \cong \Ind^{\frakg_{\ubar d}}_{{^w}\frakq_{\ubar d}}\IC(\scrS).
\]
It follows that
\begin{align*}
	P_{\bfV} &= \bigoplus_{\frakq\in \ubar \frakQ_{\zeta}} \;\bigoplus_{\frakp_*\in \ubar \frakP^{\epsilon}_{\xi}}\Hom^{\bullet}_{\Gz}(\Ind^{\frakg_{\ubar d}}_{\frakq_{\ubar d}}\IC(\scrS), \Ind^{\frakg_{\ubar d}}_{\frakp_{d}} \IC(\scrC))^{\wedge}_0 \\
	&\cong\bigoplus_{w\in p(W^{\aff}_{\xi, \bfx})\backslash W_{\zeta} }\;\bigoplus_{y\in   W^{\aff}_{\xi, \bfx}\backslash W^{\aff}_{\xi}}\Hom^{\bullet}_{\Gz}(\Ind^{\frakg_{\ubar d}}_{{^{\tilde w}}\frakp_{d}}\IC(\scrC), \Ind^{\frakg_{\ubar d}}_{{^{y}}\frakp_{d}} \IC(\scrC))^{\wedge}_0  \cong \calH_{\xi} 1_{\bfV}.
\end{align*}

\end{proof}

We obtain an algebraic characterisation for \emph{parabolically induced sheaves}:
\begin{coro}\label{coro:KZ-sheaf}
	Let $\pi\in\Pi_{\xi}$. Then, $\IC_{\pi}$ is a direct summand of $\pH^{\bullet}\Ind^{\frakg_{\ubar d}}_{\frakq_{\ubar d}}\IC(\scrS)$ for some $\frakq\in \frakQ_{\zeta}$ if and only if the simple module $\calL_{\pi}\in \calH_{\xi}\mof$ attached to $\pi$ under \autoref{coro:springer} satisfies $\bfV(\calL_{\pi})\neq 0$.
\end{coro}
\subsection{Projective and relatively injective modules}\label{subsec:proj-inj}

Set $\epsilon = d / |d|$. Consider the subcategories $\calH_{\xi}\proj\subset \calH_{\xi}\Mod$ and $\calH_{\xi}\inj_{\calZ_{\xi}}\subset \calH_{\xi}\Mod$ from \autoref{subsec:AKZ}. The same definitions apply to the opposite $\calH_{\xi}^{\op}$. \par

Given $\scrF\in \Db_{\Gz}(\frakg_{\ubar d})$, we set
\begin{align*}
	\calM_{\scrF} = \bigoplus_{\frakp_*\in \ubar\frakP^{\epsilon}_{\xi}}\Hom^{\bullet}_{\Gz}(\scrF, \Ind^{\frakg_{\ubar d}}_{\frakp_d}\IC(\scrC))^{\wedge}\in \calH_{\xi}\Mod, \\
	\calN_{\scrF} = \bigoplus_{\frakp_*\in \ubar\frakP^{\epsilon}_{\xi}}\Hom^{\bullet}_{\Gz}(\Ind^{\frakg_{\ubar d}}_{\frakp_d}\IC(\scrC), \scrF)^{\wedge}\in \calH_{\xi}^{\op}\Mod,
\end{align*}
where the $\calH_{\xi}$-module structure is given by the Yoneda product.

\begin{prop}\label{prop:proj-inj}
	For every simple orbital sheaf $\scrF\in \Irr\Perv_{\Gz}(\frakg^{\nil}_{\ubar d})_{\xi}$, we have $\calM_{\scrF}\in \calH_{\xi}\proj$ and $\calN_{\scrF}\in \calH_{\xi}^{\op}\proj$. For every simple anti-orbital sheaf $\scrF\in \Irr\Perv_{\Gz}(\frakg_{\ubar d})^{\nil}_{\xi}$, we have $\calM_{\scrF}\in \calH_{\xi}\inj_{\calZ_{\xi}}$ and $\calN_{\scrF}\in \calH_{\xi}^{\op}\inj_{\calZ_{\xi}}$. 
\end{prop}
\begin{proof}
	Let $\scrF\in \Irr\Perv_{\Gz}(\frakg^{\nil}_{\ubar d})_{\xi}$. Then, there exists $\frakp_*\in \frakP^{\epsilon}_{\xi}$ such that $\scrF$ is a direct summand of the perverse cohomology of $\Ind^{\frakg_{\ubar d}}_{\frakp_{d}} \IC(\scrC)$. It follows that $\calM_{\scrF}$ is a direct summand of $\calH_{\xi}$. Therefore $\calM_{\scrF}\in \calH_{\xi}\proj$ holds. Similar arguments show that $\calN_{\scrF}\in\calH_{\xi}^\op\proj$ holds.\par
	Let $\scrF\in \Irr\Perv_{\Gz}(\frakg_{\ubar d})^{\nil}_{\xi}$. Then, there exists $\frakp'_*\in \frakP^{-\epsilon}_{\xi}$ such that $\scrF$ is a direct summand of the perverse cohomology of the $(-\epsilon)$-spirally induced complex $\scrG := \Ind^{\frakg_{\ubar d}}_{\frakp'_{d}} \IC(\scrC)$. \autoref{lemm:dualM} below implies that $\calM_{\scrG}\in \calH_{\xi}\inj_{\calZ_{\xi}}$ holds.
\end{proof}

\begin{lemm}\label{lemm:dualM}
	Let $\frakp'_*\in \frakP_{\xi}^{\epsilon}$ and $\frakp''_*\in \frakP_{\xi}^{-\epsilon}$ be opposite spirals. Put $\scrG = \Ind^{\frakg_{\ubar d}}_{\frakp'_d}\IC(\scrC)$ and $\scrG' = \Ind^{\frakg_{\ubar d}}_{\frakp''_d}\IC(\scrC)$. Then, are canonical isomorphisms in $\calH_{\xi}^{\op}\Mod$:
	\[
		\calM^{\vee}_{\scrG} := \bigoplus_{x\in W^{\aff}_{\xi}\bfx} \Hom_{\calZ_{\xi}}(\calM_{\scrG}1_x, \calZ_{\xi})\cong \calN_{\scrG'},\quad \calM^{\vee}_{\scrG'}\cong \calN_{\scrG}.
	\]
\end{lemm}
\begin{proof}
	For any $\scrF\in \Db_{\Gz}(\frakg^{\nil}_{\ubar d})_{\xi}$, we have isomorphisms:
	\begin{align*}
		\Hom^{\bullet}_{\Gz}(\scrG, \scrF)&\cong \Hom^{\bullet}_{L_0}(\IC(\scrC), \cores^{\frakg_{\ubar d}}_{\frakp'_d}\scrF) \cong \Hom^{\bullet}_{L_0}(\cores^{\frakg_{\ubar d}}_{\frakp'_d}\scrF, \IC(\scrC))^{\vee} \\
		&\cong \Hom^{\bullet}_{L_0}(\Res^{\frakg_{\ubar d}}_{\frakp''_d}\scrF, \IC(\scrC))^{\vee} \cong \Hom^{\bullet}_{\Gz}(\scrF, \scrG')^{\vee},
	\end{align*}
	where the first one and the last one are due to adjunction, the second one \autoref{lemm:duality} below, the third one \autoref{prop:hyperbolic-spiral}. 
	Since these isomorphisms are functorial in $\scrF$, we obtain the isomorphism $\calM^{\vee}_{\scrG}\cong \calN_{\scrG'}$ by summing over $\scrF = \Ind^{\frakg_{\ubar d}}_{\frakp_d}\IC(\scrC)$ for $\frakp_*\in \ubar\frakP^{\epsilon}_\xi$. The proof of the other isomorphism is similar.
\end{proof}
\begin{lemm}\label{lemm:duality}
	Set $\calS = (\rmH^{\bullet}_{Z^{\circ}_L})^{W^{\aff}_{\xi, \bfx}}$, regarded as commutative formal differential graded algebra. There is a bi-functorial isomorphism in $\DG(\calS)$:
	\[
		\RHom_{L_0}(\scrF, \scrG) \cong \RHom_{\calS}(\RHom_{L_0}(\scrG, \scrF), \calS[-2D])
	\]
	some $D\in \bbN$.
\end{lemm}
\begin{proof}
	Let $\calB = \rmH^{\bullet}_{Z^{\circ}_L}$. The arguments of \cite[Proposition 2.4]{RR21} apply for $\bbZ$-graded Lie algebras $\frakl_*$ and show that there is an equivalence of triangulated categories:
	\begin{equation}\label{equa:bloc-cusp}
		\Db_{L_0}(\frakl_\eta)_{\xi} \cong \DG(\calB),%\quad \scrF\mapsto \RHom_{L_0}(\IC(\scrC), \scrF)
	\end{equation}
	where $\DG(\calB)$ is the derived category of finitely generated dg-modules over $\calB$. As $Z^{\circ}_L$ is a torus, $\calB$ is a graded polynomial ring with generators in degree $2$. As $W^{\aff}_{\xi,\bfx}$ acts on $\fraka_{\xi}$ by reflections, the invariant subalgebra $\calS$ is a graded polynomial ring and $\calB$ is a symmetric algebra over $\calS$ with trace map $\tau:\calB\to \calS$, which is homogeneous of degree $-2D$, where $D = \#\{ \alpha\in R^{\aff}_{\xi}\;\vert\; \alpha(\bfx) = 0 \}/2$. \par

	Suppose that $M, N$ are semi-free dg-modules over $\calB$. The trace pairing
	\[
		\Hom^{\bullet}_{\calB}(M, N) \times \Hom^{\bullet}_{\calB}(N, M)\to \Hom^{\bullet}_{\calB}(N, N)\xrightarrow{\tr} \calB\xrightarrow{\tau}\calS[-2D]
	\]
	induces a quasi-isomorphism $\Hom^{\bullet}_{\calB}(M, N)\to \Hom^{\bullet}_{\calS}(\Hom^{\bullet}_{\calB}(N, M), \calS[ -2D])$. Via the equivalence \eqref{equa:bloc-cusp}, it produces the required isomorphism.
\end{proof}

\subsection{Proof of \autoref{theo:biorbital}}\label{subsec:proof}
The ``if'' part is due to the fact that the parabolic induction preserves complexes on the nilpotent cone	and commutes with the Fourier--Sato transform, see \autoref{subsec:parab-ind}. We prove the ``only if'' part. Assume that $\IC_{\pi}$ is bi-orbital.  \par
Similar arguments as in the proof of \autoref{lemm:HomFG} reduce the problem to the case where $d = 1$ and $G$ is simply connected and almost simple. Under this assumption, we may apply the results of \autoref{subsec:springer} and \autoref{subsec:Vfaisc}. Let $\xi\in \frakT$ be an LY-cuspidal support such that $\pi\in \Pi_{\xi}$ and let $\calH_{\xi}$ be the spectral completion of degenerate double affine Hecke algebra attached to $\xi$ (\autoref{subsec:springer}). \autoref{prop:proj-inj} implies that $\calM_{\pi} := \calM_{\IC_{\pi}}$ lies in $\calH_{\xi}\proj\cap \calH_{\xi}\inj_{\calZ_{\xi}}$. By \autoref{theo:KZ}, this condition implies that the cosocle $\calL_{\pi} = \operatorname{cosoc}_{\calH_{\xi}}(\calM_{\pi})\in \calH_{\xi}\mof$ satisfies $\bfV(\calL_{\pi})\neq 0$. Hence, \autoref{coro:KZ-sheaf} implies that $\IC_{\pi}$ is a direct summand of the perverse cohomology of $\Ind^{\frakg_{\ubar 1}}_{\frakq_{\ubar 1}}\IC(\scrS)$ for some $\frakq\in \frakQ_{\zeta}$. \par
It remains to show that $\IC(\scrS)$ is bi-orbital. Indeed, the adjunction yields
\[
	\Hom^{\bullet}_{M_{\ubar 0}}(\Res^{\frakg_{\ubar 1}}_{\frakq_{\ubar 1}} \IC_{\pi}, \IC(\scrS))\cong \Hom^{\bullet}_{M_{\ubar 0}}(\IC_{\pi}, \Ind^{\frakg_{\ubar 1}}_{\frakq_{\ubar 1}} \IC(\scrS)) \neq 0,
\]
so that $\IC(\scrS)$ is a direct summand of $\pH^{\bullet}\Res^{\frakg_{\ubar 1}}_{\frakq_{\ubar 1}} \IC_{\pi}$ by \autoref{prop:dist-induction}. Since 
			\[
				\FS\left(\Res^{\frakg_{\ubar 1}}_{\frakq_{\ubar 1}} \IC_{\pi} \right)\cong \Res^{\frakg_{-\ubar 1}}_{\frakq_{-\ubar 1}} \FS(\IC_{\pi})\in \Db_{M_{\ubar 0}}(\frakm^{\nil}_{-\ubar 1}),
			\]
			it follows that $\IC(\scrS)$ is bi-orbital.  \hfill \qed

%With this result at hand, we can derive analogues of \autoref{theo:biorbital}:
	\section{Applications}\label{sec:app}
We derive some consequences from the main results \autoref{theo:supercuspidal} and \autoref{theo:biorbital}. 

\subsection{Bi-orbital cuspidal pairs}\label{subsec:cusp}
We keep the setup of \autoref{subsec:G}. Let $\Db_{\Gz}(\frakg^{\nil}_{\ubar d})^{\mathrm{old}}$ be the thick subcategory of $\Db_{\Gz}(\frakg^{\nil}_{\ubar d})$ generated by the image of the parabolic induction functor (\autoref{subsec:parab-ind})
\[
	\Ind^{\frakg_{\ubar d}}_{\frakq_{\ubar d}}: \Db_{M_{\ubar 0}}(\frakm^{\nil}_{\ubar d})\to \Db_{\Gz}(\frakg^{\nil}_{\ubar d})
\]
along every $\theta$-stable proper parabolic subalgebra $\frakq\subset \frakg$ and from any $\theta$-stable Levi factor $\frakm\subseteq \frakq$. We call a pair $\pi\in \Pi$ \emph{cuspidal} if $\IC_{\pi}$ does not lie in $\Db_{\Gz}(\frakg_{\ubar d})^{\mathrm{old}}$.  \par
From \autoref{theo:supercuspidal}, it is known that supercuspidal pairs are cuspidal. The converse is true under the bi-orbital condition (\autoref{subsec:statement}):
\begin{theo}\label{theo:cusp-biorbital}
	Every cuspidal pair in $\Pi^{\nil}$ is supercuspidal.
\end{theo}
\begin{proof}
	Suppose that $\pi\in \Pi^{\nil}$ is cuspidal. Then, \autoref{theo:biorbital} shows that there exists a supercuspidal support $\zeta = (M, \rmO, \scrS)$ on $\frakg_{\ubar d}$ satisfying $(\rmO, \scrS)\in \Pi^{\nil}_M$ and a parabolic subalgebra $\frakq\in \frakQ_{\zeta}$ such that $\IC_\pi$ is a direct summand of $\pH^{\bullet}\Ind^{\frakg_{\ubar d}}_{\frakq_{\ubar d}}\IC(\scrS)$, so it lies in the thick subcategory generated by $\Ind^{\frakg_{\ubar d}}_{\frakq_{\ubar d}}\IC(\scrS)$. The cuspidality of $\pi$ implies that $\frakq = \frakg$ holds and it follows that $\IC_{\pi}\cong \IC(\scrS)$. 
\end{proof}
\begin{rema}
	Cuspidal pairs may not be supercuspidal in general. The former has been studied in \cite{VX21, VX23}. In a forthcoming joint work with Tsai C.-C., Xue T. and K. Vilonen, we will show that every cuspidal sheaf (namely, IC-extension of a cuspidal pair) on a $\bbZ/m$-graded Lie algebra arises from a nearby-cycle construction which takes bi-orbital supercuspidal sheaves as input. The nearby-cycle functor increases the semisimple part of the singular support of a cuspidal sheaf as result. \autoref{theo:cusp-biorbital} can therefore be regarded as the base case in our classification of cuspidal sheaves.
\end{rema}

\subsection{Affine Hecke algebras at roots of unity}\label{subsec:groj}

In the rest of this section, let $G$ be a connected complex reductive group. Choose a Borel pair $(B, T)$ for $G$. Let $\frakg = \Lie G$, $\frakb = \Lie B$ and $\frakt = \Lie T$ denote the corresponding Lie algebras and set 
\begin{align*}
	\dot\frakg &= \left\{ (gB, z)\in G/B\times \frakg\;\vert\; \Ad_{g}^{-1} z\in [\frakb, \frakb] \right\},\\
	\ddot\frakg &= \left\{ (gB, g'B, z)\in G/B\times G/B\times \frakg\;\vert\; \Ad_{g}^{-1} z\in [\frakb,\frakb],\; \Ad_{g'}^{-1} z\in  [\frakb, \frakb] \right\} = \dot\frakg\times_{\frakg}\dot\frakg.
\end{align*}
There is an action of $G\times \bbC^{\times}$ on $\dot\frakg$ and $\ddot\frakg$ by:
\[
	((h,r), (gB, z)) \mapsto  (hgB, r^{-1}\Ad_h z),\quad ((h,r), (gB, gB', z)) \mapsto  (hgB, hg'B, r^{-1}\Ad_h z).
\]
The fixed point subvarieties of a pair $(s, q)\in T\times \bbC^{\times}$ in $\frakg$, $\dot\frakg$ and $\ddot\frakg$ under this action are given by
\begin{align*}
	\frakg^s_q &= \left\{z\in \frakg\;\vert\; \Ad_s z = qz\right\},\quad
	(G/B)^{s} = \left\{ gB\in G/B\;\vert\; g^{-1}sg\in B\right\}, \\
	\dot\frakg^s_q &= \dot\frakg\cap ((G/B)^{s}\cap \frakg^s_q),\quad
	\ddot\frakg^s_q = \dot\frakg^s_q\times_{\frakg^s_q}\dot\frakg^s_q.
\end{align*}
They are stable under the action of the subgroup $G^{s}\times \bbC^{\times} = Z_G(s)\times \bbC^{\times}\subseteq G\times \bbC^{\times}$. Given any $z\in \frakg^{s}_q\cap \frakg^{\nil}$, the centraliser $G^{s}_z = Z_G(s, z)$ acts by left multiplication on the Springer fibre 
\[
	(G/B)^{s}_z = \left\{ gB\in (G/B)^s\;\vert\; \Ad_{g^{-1}}z\in \frakb  \right\},
\]
inducing an action of the component group $\pi_0(G^s_z)$ on the cohomology group $\rmH^{\bullet}((G/B)^s_z)$.

\begin{lemm}\label{lemm:principal}
	Let $(s, q)\in T\times \bbC^{\times}$ and assume that $q\neq 1$ and $G^s$ is connected. Then, given an element $z\in \frakg^{s}_q\cap \frakg^{\nil}$ and a character $\chi\in \pi_0(G^s_z)^{\wedge}$, the following conditions for $(s, \chi)$ are equivalent:
	\begin{enumerate}
		\item\label{lemm:principal-i}
			The $\chi$-isotypic component of the $\pi_0(G^s_z)$-action on $\rmH^{\bullet}((G/B)^s_z)$ is non-zero and the intersection complex $\IC_{\chi}$ with coefficients in $\chi$ has nilpotent singular support:
			\[
				\operatorname{SS}(\IC_{\chi})\subseteq (\frakg^s_q \times \frakg^s_{q^{-1}})\cap (\frakg^{\nil}\times\frakg^{\nil}).
			\]
		\item\label{lemm:principal-ii}
			$\IC_{\chi}$ is a direct summand of $\pH^{\bullet}\pi^s_*\bbC$, where $\pi^s: \dot\frakg^{s}_q\to \frakg^{s}_q$ is the projection.

	\end{enumerate}
	%\begin{enumerate}
		%\item\label{coro:principal-i}
		%$(\rmC, \scrL)\in \Pi_{\xi}$ if and only if $\rmH^{\bullet}((G/B)^{\theta}_z, \bbC)^\scrL\neq 0$ for $z\in \rmC$.
		%\item\label{coro:principal-ii}
			%$(\rmC, \scrL)\in \Pi_{\xi}^\nil$ if and only if there exists a $\theta$-stable Borel subgroup $B'$ containing $T$ such that $\IC(\scrL)$ is a simple constituent of $\pH^{\bullet}\Ind^{\frakg_{\ubar d}}_{\frakb'_{\ubar d}}\IC(0)$.
	%\end{enumerate}
\end{lemm}
\begin{proof}
	If $q$ is not a root of unity, then $\frakg^s_{q^{\pm 1}} \subseteq \frakg^{\nil}$ holds, so the statement is reduced to the main theorem of \cite{KL}. We assume therefore that $q$ is a root of unity of order $m\in \bbZ_{\ge 1}$. \par
	We treat first the case where $s^m\in Z_G$ holds. In this case, we may find a cocharacter $\iota:\bbC^{\times}\to G^{\ad} = G/Z_G$ such that $\iota(q) = \Ad_s$. We pick a $\iota$-stable Borel pair $(B, T)$ of $G$. Let $\theta = \iota\mid_{\mu_m}$ be the restriction. Then, $\theta$ yields a $\bbZ/m$-grading on $G$ as in \autoref{subsec:grading} and $\xi_0 = (T,\frakt_*, \delta_0)\in \frakT$ is an LY-cuspidal support. We see readily that $\frakg^s_q = \frakg_{\ubar 1}$ and $G^s = G_{\ubar 0}$ hold. We are in the situation of \autoref{exam:principal1} and $(G/B)^{\theta} = (G/B)^s$. Let $\scrL_{\chi}$ denote the $\Gz$-equivariant local system given by $(\scrL_{\chi})_z = \chi$. Then, \autoref{theo:comparaison-bloc} implies that $(\Ad_{\Gz} z, \scrL_{\chi})$ lies in $\Pi_{\xi_0}$ if and only if the $\chi$-isotypic component of the $G^s_z$-action on $\rmH^{\bullet}((G/B)^s_z)$ is non-zero. The statement follows from \autoref{theo:biorbital}. \par
	In the general case where $s\in T$ may not satisfy $s^m\in Z_G$, we consider 
	\[
		H = Z_G(s^m), \quad B_H = B\cap H, \quad \frakh := \Lie H = \bigoplus_{\omega\in \mu_m}\frakg^s_\omega.
	\]
	We have, $\frakg^s_{q^{\pm 1}} = \frakh^s_{q^{\pm 1}}$ and $G^s = H^s$ hold. We define the analogue of $\dot\frakg$ for $\frakh$:
	\[
		\dot\frakh = \left\{ (hB_H, z)\in (H/B_H)\times \frakh\;\vert\; \Ad_h^{-1} z\in [\frakb,\frakb]\cap \frakh\right\}.
	\]
	The varieties of fixed points $(H/B_H)^s$ and $\dot\frakh^s_q$ are defined in a similar manner.
	As $s^m\in Z_H$, the previous case applies to $H$, so that the conditions (i) and (ii) are equivalent for the datum $(H, z, \chi)$. \par
	Choose a lifting $\dot w\in N_G(T)$ for every $w\in W(G, T) = N_G(T) / T$ and set
	\[
		\Sigma = \{\dot w\;\vert\; w\in W(G,T),\quad \dot w B\dot w^{-1} \cap H = B_H \}\subseteq N_G(T).
	\]
	Then, there are $G^s$-equivariant morphisms \\
	\begin{minipage}[r]{.5\linewidth}
		\[
			\begin{tikzcd}
				\dot\frakh^s_q\times \Sigma \arrow{r}{\rho}[swap]{\sim}\arrow{d}{\pi_H^s\circ \pr_1}&\dot\frakg^s_q \arrow{d}{\pi^s} \\
				\frakh^s_q \arrow{r}{=} & \frakg^s_q, \\
			\end{tikzcd}
		\]
	\end{minipage}
	\begin{minipage}[l]{.7\linewidth}
		$\rho(hB_H, z, \dot w) = (h\dot wB, z)$ \\
		$\pi^s(gB, z) = z$ \\
		$\pi_H^s(hB_H, z) = z$,
	\end{minipage}
	where $\rho$ is an isomorphism. This diagram induces isomorphisms
	\begin{align*}
		\Hom_{G^s_z}(\chi, \rmH^{\bullet}((G/B)^s_z))&\cong \Hom_{H^s_z}(\chi, \rmH^{\bullet}((H/B_H)^s_z))\otimes \bbC\Sigma \\
		\pH^{\bullet}\pi^s_*\bbC&\cong \pH^{\bullet}\pi^s_{H*}\bbC\otimes \bbC\Sigma.
	\end{align*}
	From these isomorphisms we deduce that the condition \ref{lemm:principal-i} for $(G, z, \chi)$ is equivalent to \ref{lemm:principal-i} for $(H, z, \chi)$ and \ref{lemm:principal-ii} for $(G, z, \chi)$ equivalent to \ref{lemm:principal-ii} for $(H, z, \chi)$. This concludes the proof.
\end{proof}

In the rest of this section, we assume furthermore that the derived subgroup $[G, G]$ is simply connected and semisimple. Let $R = R(G, T)$ be the set of roots, $\Delta\subseteq R$ the basis attached to $B$ and $W = N_G(T)/T$ the Weyl group generated by the simple reflections $\left\{ s_{\alpha} \right\}_{\alpha\in \Delta}$. The affine Hecke algebra $\bfK$ attached to $G$ with parameter $q := \exp(2\pi i(\eta / m))$ is the unital associative algebra over $\bfC$ generated by $\left\{ X^{\mu} \right\}_{\mu\in \bfX^*(T)}\cup \left\{ T_w \right\}_{w\in W}$ subject to the following relations:
\begin{align*}
	&X^{0} = 1;\quad X^{\mu}X^{\nu} = X^{\mu + \nu}, \quad \forall\mu,\nu\in \bfX^*(T);\\
	&T_wT_y = T_{wy} \quad \forall w,y\in W\; \text{such that $\ell(wy) = \ell(w) + \ell(y)$};\\
	&(T_{\alpha} - q)(T + 1) = 0,\quad \forall\alpha\in \Delta\text{, where $T_{\alpha} = T_{s_{\alpha}}$ };\\
	& T_{\alpha} X^{\mu} - X^{s_{\alpha}(\mu)}T_{\alpha} = (q - 1)(X^{\mu} - X^{s_{\alpha}(\mu)})/(1 - X^{-\alpha}),\quad \forall\alpha\in \Delta,\; \forall\mu\in \bfX^*(T).
\end{align*}
The Deligne--Langlands correspondence for affine Hecke algebras proven by Kazhdan--Lusztig \cite{KL} makes use of a ring isomorphism
\begin{equation}\label{equa:isoK1}
	\calK^{s} = \bbC[T/W]^{\wedge}_{[s]}\otimes_{\bbC[T/W]} \bfK \xrightarrow{\sim} K_0^{G^s}(\ddot\frakg^s_q)^{\wedge}_1,\quad [s] = Ws\in T/W.
\end{equation}
The bivariant Riemann--Roch theorem yields an isomorphism 
\begin{equation}\label{equa:isoK2}
	K_0^{G^s}(\ddot\frakg^s_q)^{\wedge}_1\cong \rmH_*^{G^s}(\ddot\frakg^s_q)^{\wedge}_0.
\end{equation}
Via the homomorphisms \eqref{equa:isoK1} and \eqref{equa:isoK2} one obtains an action of the affine Hecke algebra $\bfK$ on the cohomology group $\rmH^{\bullet}((G/B)^s_z, \bbC)$ for every $z\in \frakg^{\nil}\cap \frakg^s_q$, called {\itshape Springer action}, which commutes with the $G^s_z$-action. For $(s, q)\in T \times \bbC^{\times}$, let $\Pi^{\nil}_{s,q}$ denote the set of $G^s$-conjugacy classes of pairs $(z, \chi)$ satisfying condition \ref{lemm:principal-i} of \autoref{lemm:principal} and for each $(z, \chi)\in \Pi^{\nil}_{s,q}$, we set
\[
	\ba\Delta_{z, \chi} = \Hom_{G^s_z}(\chi, \rmH^{\bullet}((G/B)^s_z, \bbC))\in \bfK\mof_s.
\]
Let $\operatorname{cosoc}_{\bfK}\ba\Delta_{z, \chi}$ denote the cosocle of $\ba\Delta_{z, \chi}$ as $\bfK$-module.
\begin{theo}
	For $q\neq 1$ and $s\in T$, the Springer action of $\bfK$ with parameter $q$ induces a bijection
	\[
		\Pi^\nil_{s,q}\xrightarrow{\sim}\Irr\bfK\mof_{s},\quad (z,\chi) \mapsto \operatorname{cosoc}_{\bfK}\ba\Delta_{z, \chi}.
	\]
\end{theo}
\begin{proof}
	Let $\calE = \End^{\bullet}_{G^s}(\pi^s_*\bbC)^{\wedge}_0$ denote the completed extension algebra. 
	The general theory of extension algebras \cite[\S 8]{CG} implies that there is a ring isomorphism $\calE\cong \rmH^{G^s}_{\bullet}(\ddot\frakg^s_q)^{\wedge}_0$ and the assignment
	\[
		\scrF \mapsto \operatorname{cosoc}_{\calE}\Hom^{\bullet}_{G^s}(\scrF, \pi^s_*\bbC)^{\wedge}_0 = \Hom_{G^s}(\scrF, \pH^{\bullet}\pi^s_*\bbC)
	\]
	is a bijection from the set of (isomorphism classes of) simple constituents of $\pH^{\bullet}\pi^s_*\bbC$ to the set of simple $\calE$-modules. The simple connectedness of $[G,G]$ implies that $G^s$ is connected. By \autoref{lemm:principal}, the simple constituents are simple perverse sheaves of the form $\IC_{\chi}$ with $(z, \chi)\in \Pi^{\nil}_{s, q}$. The arguments of \cite{kato17} via weight structure show that $\ba\Delta_{z, \chi}$ is a quotient of $\Hom^{\bullet}_{G^s}(\IC_{\chi}, \pi^s_*\bbC)^{\wedge}_0$, so they have the same simple cosocle as $\calE$-modules. Since $\ba\Delta_{z, \chi}$ is finite-dimensional, the isomorphisms \eqref{equa:isoK1} and \eqref{equa:isoK2} imply that
	\[
		\operatorname{cosoc}_{\bfK}\ba\Delta_{z, \chi} = \operatorname{cosoc}_{\calE}\ba\Delta_{z, \chi} = \Hom_{G^s}(\IC_{\chi}, \pH^{\bullet}\pi^s_*\bbC).
	\]
	This concludes the proof.
\end{proof}
\begin{rema}
	This result was originally announced by I. Grojnowski in \cite[Theorem 1]{groj94}. It generalises the Deligne--Langlands correspondence established in \cite{KL} to the case where $q$ is allowed to be a root of unity.
\end{rema}

\iffalse 
	\subsection{Lusztig's question}
	Lusztig has shown in \cite[\S 7]{lusztig89b} that given $q\in \bbC^{\times}$ the $q$-eigenspace $\frakg^{s}_q$ is contained in the nilpotent cone $\frakg^{\nil}$ for every semisimple element $s\in G$ if and only if $q$ is not a root of the Poincar\'e polynomial: 
	\begin{equation}\label{equa:poincare}
		\sum_{w\in W} q^{\ell(w)} \neq 0. 
	\end{equation}
	We obtain as corollary an affirmative answer to a question raised by Lusztig:
	\begin{coro}
		The Deligne--Langlands correspondence proven in \cite{KL} remains valid for $\bfK$ as long as the condition\eqref{equa:poincare} holds.
	\end{coro}
	\begin{proof}
		Consider first the case where $q$ and $s$ are of finite order.
	\end{proof}
	Another consequence is the following:
	\begin{coro}
		Assume that \eqref{equa:poincare} holds. Then, the projective dimension the regular $\bfK$-bimodule is at most $2\dim S$. 
	\end{coro}
\fi

\section{Equivariant sheaves under \texorpdfstring{$\theta$}{theta}-invariant subgroups}\label{sec:equiv}
We have been working in the preceding sections with $\Gz$-equivariant sheaves, where $\Gz = (G^{\theta})^{\circ}$ is the neutral component of the invariant subgroup $G^{\theta}$, in order to avoid technical complexity with disconnected groups. In this section, we explain the relation between $\Gz$-equivariant and $G^{\theta}$-equivariant sheaves. Analogues of the main results for $G^{\theta}$-equivariant sheaves can be proven via change of equivariance.

\subsection{Component group}
Let $(B, T)$ be a $\theta$-stable Borel pair for $G$, $\rho:\til G\to G$ the simply connected cover of the derived subgroup $[G,G]$ and $\til T = \til G\times_G T$. The component group of $G^{\theta}$, denoted by $\Gamma$, can be described in terms of group cohomology:
\[
	\Gamma := G^{\theta} / \Gz = \ker(\rmH^1(\mu_m, \pi_1(G))\to \rmH^1(\mu_m, \til G)),
\]
where $\mu_m$ acts on $\pi_1(G)$ and $\til G$ via $\theta$. In particular, $\Gamma$ is a finite abelian group. 

\subsection{Change of equivariance for sheaves}
We will write $\hat\Pi$ and $\hat\Pi^{\nil}$ for the $G^{\theta}$-equivariant counterparts of $\Pi$ and $\Pi^{\nil}$. Similarly, given a $\theta$-isotropic Levi subgroup $M\subseteq G$, we will write $\hat\Pi_M$ and $\hat\Pi^{\nil}_M$ for the $M^{\theta}$-equivariant counterparts of $\Pi_M$ and $\Pi^{\nil}_M$. The notion of supercuspidal pairs (\autoref{defi:supercuspidal}), supercuspidal support (\autoref{defi:super-support}) and cuspidal pairs (\autoref{subsec:cusp}) also have an obvious $G^{\theta}$-equivariant counterpart. \par

Given subgroups $\Gz\subseteq K\subseteq H\subseteq G^{\theta}$, we define an induction functor
\[
	\Ind^{H}_{K} := f_*f^* = f_!f^!: \Db_{K}(\frakg_{\ubar d})\to \Db_{H}(\frakg_{\ubar d}),
\]
where $f$ is the finite \'etale morphism:
\[
	f: H\times^{K} \frakg_{\ubar d}\to \frakg_{\ubar d},\quad [g : z]\mapsto \Ad_gz.
\]
The finite \'etale property of $f$ implies that $\Ind^{H}_{K}$ is perverse t-exact and it preserves the complexes concentrated in nilpotent cone. Moreover, the induction is transitive in the obvious sense and there is a natural bi-adjunction between $\Ind^{H}_K$ and the restriction of equivariance: $\Res^H_K:\Db_{H}(\frakg_{\ubar d})\to \Db_{K}(\frakg_{\ubar d})$. 
\par

The adjoint action induces an action of $G^{\theta}$ on $\Pi$ preserving $\Pi^{\nil}$ and this action factorises through the component group $G^{\theta}\to \Gamma$. Given $\pi = (\rmC , \scrL )\in \Pi$, we write $[\pi]\in \Pi/\Gamma$ for the $\Gamma$-orbit of $\pi$ and set
\[
	\quad H_\pi = \Stab_{G^{\theta}}(\pi),\quad \Gamma_{\pi} = H_{\pi} / \Gz,\quad  \calA_{\pi} = \End_{H_\pi}(\Ind^{H_\pi}_{\Gz}\IC_{\pi}). 
\]
It is a standard fact that $\calA_{\pi}$ is a twisted group algebra of $\Gamma_{\pi}$ whose isomorphism class is classified by an element of the group cohomology $\rmH^2(\Gamma_{\pi}, \bbC^{\times})$. The latter vanishes because $\Gamma_{\pi}\subseteq \Gamma$ is abelian. Therefore, $\calA_{\pi}$ is semisimple and commutative.\par

Let $\calA_{\pi}^{\wedge}$ denote the set of simple characters of $\calA_{\pi}$. Given $\chi\in \calA_{\pi}^{\wedge}$, the $\chi$-isotypic component of $\Ind^{H_\pi}_{\Gz}\IC_{\pi}$ is a simple $H_{\pi}$-equivariant perverse sheaf, denoted by $\scrF_{\chi, \pi}$.
%\par
%The relation between $\Pi$ and $\ha\Pi$ is the following:
\begin{prop}
	The induction $\Ind^{G^{\theta}}_{H_{\pi}} \scrF_{\chi, \pi}$ is a simple $G^{\theta}$-equivariant perverse sheaf. This yields a bijection:
	\[
		\bigsqcup_{[\pi]\in \Pi / \Gamma} \calA_{\pi}^{\wedge} \xrightarrow{\sim} \hat\Pi 
	\]
	which sends $(\pi, \chi)$ with $\pi\in \Pi$ and $\chi\in \calA^{\wedge}_{\pi}$ to the pair $\hat\pi\in \hat\Pi$ such that $\Ind^{G^{\theta}}_{H_{\pi}}\scrF_{\chi, \pi} = \IC_{\hat\pi}$. Moreover, $\pi\in \Pi^{\nil}$ (resp. $\pi$ is cuspidal / supercuspidal) if and only if $\hat\pi\in \hat\Pi^{\nil}$ (resp. $\hat\pi$ is cuspidal / supercuspidal).\hfill\qedsymbol
\end{prop}

\subsection{Parabolic induction and restriction}\label{subsec:parab'}
In the setup of \autoref{subsec:parab-ind}, we have the following induction diagram:
\[
	\frakm_{\ubar d}\xleftarrow{\alpha} G^{\theta}\times^{U_{\ubar 0}}\frakq_{\ubar d}\xrightarrow{\beta} G^{\theta}\times^{Q^{\theta}}\frakq_{\ubar d}\xrightarrow{\gamma} \frakg_{\ubar d}, \qquad
	\frakg_{\ubar d}\xleftarrow{\delta} \frakq_{\ubar d}\xrightarrow{\varepsilon} \frakm_{\ubar d}.
\]
The parabolic induction extends to $G^{\theta}$-equivariant and $M^{\theta}$-equivariant derived categories:
\[
	\ha\Ind^{\frakg_{\ubar d}}_{\frakq_{\ubar d}} = \gamma_*(\beta^*)^{-1}\alpha^*[\dim \frakv_{\ubar 0} + \dim \frakv_{\ubar d}]:\Db_{M^{\theta}}(\frakm_{\ubar d})\to \Db_{G^{\theta}}(\frakg_{\ubar d})
\]
It commutes with the transitivity with the induction of equivariance:
\[
	\ha\Ind^{\frakg_{\ubar d}}_{\frakq_{\ubar d}}\Ind^{M^{\theta}}_{M_{\ubar d}} \cong \Ind^{G^{\theta}}_{\Gz}\Ind^{\frakg_{\ubar d}}_{\frakq_{\ubar d}}
\]
We define similarly $\ha\Res^{\frakg_{\ubar d}}_{\frakq_{\ubar d}}$ and $\ha\cores^{\frakg_{\ubar d}}_{\frakq_{\ubar d}}$, so that there is are adjunctions: $\ha\Res^{\frakg_{\ubar d}}_{\frakq_{\ubar d}}\dashv\ha\Ind^{\frakg_{\ubar d}}_{\frakq_{\ubar d}}\dashv\ha\res^{\frakg_{\ubar d}}_{\frakq_{\ubar d}}$. 
\subsection{Variant of main results}
We state the counterparts of \autoref{theo:supercuspidal}, \autoref{theo:biorbital} and \autoref{theo:cusp-biorbital} for $G^{\theta}$-equivariant sheaves, the proof of which is left to the reader.
\begin{coro}
	Given a pair $\hat\pi = \left( \hat\rmC, \hat\scrL \right)\in \hat\Pi$, the following conditions are equivalent:
	\begin{enumerate}
		\item
			The pair $\hat\pi$ is supercuspidal.
		\item
			The orbit $\hat\rmC$ is distinguished and the local system $\hat\scrL$ is clean.
	\end{enumerate}\hfill\qedsymbol
\end{coro}

\begin{coro}
	Let $\hat\pi\in \hat\Pi$. Then, $\hat\pi\in \hat\Pi^\nil$ holds if and only if there exist a supercuspidal support $\hat\zeta = (M, \hat\rmO, \hat\scrS)$ and a parabolic subalgebra $\frakq\in \frakQ_{\hat\zeta}$ such that $(\hat\rmO, \hat\scrS)\in \hat\Pi^\nil_M$ holds and $\IC_{\hat\pi}$ is a simple constituent of the perverse cohomology of $\ha\Ind^{\frakg_{\ubar d}}_{\frakq_{\ubar d}}\IC(\hat\scrS)$. \hfill\qedsymbol
\end{coro}
\begin{coro}
	Every cuspidal pair in $\hat\Pi^{\nil}$ is supercuspidal.\hfill\qedsymbol
\end{coro}

\appendix
\section{Reminder on dDAHA}\label{sec:AHA-DAHA}
For the convenience of the reader, we provide a brief review of some basic results about degenerate double affine Hecke algebras, spectral completion and algebraic KZ functor. These materials can be found in \cite{liu21,liu22,liu23}.

\subsection{Degenerate double affine Hecke algebras}\label{subsec:dDAHA}
Let $\left( \fraka, R^{\aff}, \Delta^{\aff} \right)$ be a based irreducible reduced affine root system, where $\fraka$ is a euclidean affine space, $R^{\aff}$ is a set of $\bbR$-valued affine functions on $\fraka$ and $\Delta^{\aff}\subseteq R^{\aff}$ is a basis. The affine Weyl group is denoted by $W^{\aff}$ and the extended affine Weyl group $W^{\ext}$.  \par

Let $c = ( c_\alpha )_{\alpha\in R^{\aff}}$ be a $W^{\ext}$-invariant family of complex numbers. Let $\bbC[\fraka]$ denote the ring of regular functions on $\fraka_{\bbC} = \fraka\otimes_{\bbR}\bbC$. The degenerate double affine Hecke algebra with parameters $c$ attached to $\left( \fraka, R^{\aff}, \Delta^{\aff} \right)$ is the unital associative $\bbC$-algebra on the vector space $\bfH = \bbC W^{\aff} \otimes \bbC[\fraka]$ subject to the following relations:
\begin{itemize}
	\item
		Each of the subspaces $\bbC W^{\aff}\otimes 1$ and $1\otimes \bbC[\fraka]$ is given the usual ring structure;
	\item
		$a\in \bbC W^{\aff}$ and $f\in \bbC[\fraka]$ multiply by juxtaposition: $(a\otimes 1)(1\otimes f) = a\otimes f$; 
	\item
		$\alpha\in \Delta^{\aff}$ and $f\in \bbC[\fraka]$ satisfy the following equation:
		\begin{eq}
			\left(s_{\alpha} \otimes 1 \right)(1 \otimes f) - \left(1 \otimes \prescript{s_{\alpha}}{}f \right)(s_{\alpha} \otimes 1) = 1\otimes c_{\alpha} \frac{f - \prescript{s_{\alpha}}{}f}{\alpha}.
		\end{eq}
		%where $d_{\alpha}\in\bbR_{>0}$ is the minimal number such that $\alpha + d_{\alpha}\delta\bbZ\subseteq  R^{\aff}$.
\end{itemize}

For each $x\in\fraka_{\bbC}$, let $\frako_{x}\subset \bbC[\fraka]$ be the defining ideal of $x$, which is the maximal ideal generated by $\{f - f(x)\;\vert\; f\in \bbC[\fraka]\}$.
Given any module $M\in \bfH\Mod$, consider for each $x\in \fraka_{\bbC}$ the generalised $x$-weight space in $M$:
\begin{eq}
	M_{x} = \bigcup_{N\ge 0}\left\{ a\in M\;\vert\;\frako_{x}^N a = 0\right\}.
\end{eq}
For any $\bfx\in \fraka_{\bbC}$, we define $\rmO_{\bfx}\left( \bfH \right)$ to be the full subcategory of $\bfH\Mod$ consisting of finitely generated $\bfH$-modules $M$ satisfying
\begin{eq}
	M = \bigoplus_{x\in W^{\aff} \bfx} M_{x}.
\end{eq}
In other words, the polynomial subalgebra $\bbC[\fraka]$ acts locally finitely on $M$ with eigenvalues in the $W^{\aff}$-orbit of $\bfx\in \fraka_{\bbC}$. 

\subsection{Spectral completion}\label{subsec:spectral}
We review the notion of spectral completion. Given any $W^{\aff}$-orbit $[\bfx] := W^{\aff}\bfx$ in $\fraka_{\bbC}$, the spectral completion of $\bfH$ at $[\bfx]$ is defined to be
\[
	\calH^{\bfx} = \bbC[\fraka]^{\wedge}_{[\bfx]}\otimes_{\bbC[\fraka]}\bfH,\quad \bbC[\fraka]^{\wedge}_{[\bfx]} = \bigoplus_{x\in [\bfx]} \bbC[\fraka]^{\wedge}_x.
\]
We will write $\calH = \calH^{\bfx}$ when $\bfx\in \fraka$ is fixed. It can be shown that the ring structure on $\bfH$ induces a natural non-unital ring structure on $\calH^{\bfx}$. For each $x\in [\bfx]$, let $1_x\in \bbC[\fraka]^{\wedge}_{[\bfx]}$ denote the idempotent corresponding to the projection onto the direct factor $\bbC[\fraka]^{\wedge}_{x}\subset\bbC[\fraka]^{\wedge}_{[\bfx]}$. We will also denote $1_x = 1_x\otimes 1\in \calH^{\bfx}$.  \par
Let $\calH^{\bfx}\Mod$ denote the category of left $\calH^{\bfx}$-modules $N$ satisfying $\calH^{\bfx}N = N$. The assignment $M\mapsto \bbC[\fraka]^{\wedge}_{[\bfx]}\otimes_{\bbC[\fraka]}M$ yields an exact functor $\calC^{\bfx}: \bfH\Mod\to \calH^{\bfx}\Mod$, called spectral completion.  
\begin{theo}\label{theo:compl}
	The spectral completion $\calC^{\bfx}$ restricts to an equivalence of categories $\rmO_{\bfx}(\bfH)\xrightarrow{\cong} \calH^{\bfx}\mof_{\fl}$. 
\end{theo}
\subsection{Generic fibre}
The categorical centre of $\calH$
\[
	\calZ^{\bfx} = \End(\id_{\calH\Mod})
\]
can be identified with the completion of the ring $\bbC[\fraka]^{W^{\aff}_{\bfx}}$ at the defining ideal of $\bfx\in \fraka_{\bbC}$. By the Chevalley--Shephard--Todd theorem, $\calZ^{\bfx}$ is a complete regular local ring and $\bbC[\fraka]^{\wedge}_x$ is free over $\calZ^{\bfx}$ for every $x\in [\bfx]$. Let $K$ be the field of fractions of $\calZ^{\bfx}$. 
\begin{prop}\label{prop:HK}
	The spectral completion $\calH^{\bfx}$ is free over $\calZ^{\bfx}$. Moreover, $\calH^{\bfx}_K = K\otimes_{\calZ^{\bfx}}\calH^{\bfx}$ is isomorphic to the subalgebra of endomorphisms of finite rank of some countable-dimensional $K$-vector space.
\end{prop}
\begin{proof}
	The first statement follows from the decomposition $\bfH = \bbC[a]\otimes\bbC W^{\aff}$. Consider the polynomial representation $P = \bfH / \bfH \left\{ s_{\alpha} - 1\;\vert\; \alpha\in \Delta^{\aff} \right\}$ of $\bfH$. The spectral completion $\calP^{\bfx} = \calC^{\bfx}P\in \calH^{\bfx}\Mod$ is a faithful $\calH^{\bfx}$-module free over $\calZ^{\bfx}$. One can verify that the action map $\calH^{\bfx}\to \operatorname{end}_{\calZ^{\bfx}}(\calP^{\bfx})$ becomes an isomorphism after base change to $K$, where $\operatorname{end}_{\calZ^{\bfx}}(\calP^{\bfx})\subset \End_{\calZ^{\bfx}}(\calP^{\bfx})$ is the non-unital subalgebra of endomorphisms of finite rank.
\end{proof}

\subsection{Algebraic KZ functor}\label{subsec:AKZ}
Let $\frake$ be the vector space of translations on $\fraka$. For each affine root $\alpha\in R^{\aff}$, let $\dot \alpha\in \frake^*$ denote its differential. The finite root system underlying $(\fraka, R^{\aff})$ is defined to be $(\frake, R)$, where $R = \left\{ \dot \alpha\;\vert\; \alpha\in R^{\aff} \right\}$. There exists $\alpha_0\in \Delta^{\aff}$ such that $\Delta := \left\{ \dot \alpha\;\vert\; \alpha\in \Delta^{\aff} \setminus \left\{ \alpha_0 \right\} \right\}$ is a basis for $(\frake, R)$. We fix the choice of $\alpha_0$. 
Let $\Lambda\subseteq W^{\aff}$ denote the subgroup of translations, it is a subgroup of finite index of the coroot lattice of the underlying finite root system of $(\fraka, R^{\aff})$. The quotient group $W^{\aff} / \Lambda$ is canonically isomorphic to the Weyl group $W$ of $(\frake, R)$. Moreover, the stabiliser $W^{\aff}_{\bfx} = \Stab_{W^{\aff}}(\bfx)$ is mapped injectively into $W$ under the quotient map $p:W^{\aff}\to W^{\aff} / \Lambda\cong W$. 

Fix an orbit $W^{\aff}\bfx\subseteq \fraka_{\bbC}$. For each $[w]\in W^{\aff} / W^{\aff}_{\bfx}\Lambda $, choose an element $\bfx_w\in w\bfx + \Lambda$ such that $\Re\langle \bfx_w, \alpha\rangle \ll 0$ holds for every $\alpha\in \Delta$. Consider the idempotent $1_{\bfV} := \sum_{[w]\in W^{\aff} / W^{\aff}_{\bfx}\Lambda }1_{\bfx_w}\in \calH^{\bfx}$. The projective $\calH$-module $\calH 1_{\bfV}$ is independent of the choice of $\left\{ \bfx_w \right\}_{[w]\in W^{\aff} /  W^{\aff}_{\bfx}\Lambda}$ up to canonical isomorphisms. The algebraic KZ functor introduced in \cite{liu22} is the following: 
\[
	\bfV: \calH^{\bfx}\Mod\to 1_{\bfV}\calH^{\bfx}1_{\bfV}\Mod,\quad  M\mapsto \Hom_{\calH}(\calH 1_{\bfV}, M) =  1_{\bfV}M.
\]
\par

Let $(\calH^{\bfx})^{\op}$ denote the opposite ring of $\calH^{\bfx}$. Let $\calH^{\bfx}\proj$ (resp. $(\calH^{\bfx})^{\op}\proj$) denotes the category of compact projective left $\calH^{\bfx}$-modules (resp. compact projective right $\calH^{\bfx}$-modules) and let $\calH^{\bfx}\inj_{\calZ^{\bfx}}$ denote the category of compact left $\calH^{\bfx}$-modules $M$ such that the dual
\[
	M^{\vee} = \bigoplus_{x\in W^{\aff}\bfx}\Hom_{\calZ^{\bfx}}(1_{x}M, \calZ^{\bfx}) \in (\calH^{\bfx})^{\op}\Mod
\]
lies in $(\calH^{\bfx})^{\op}\proj$. Such objects are called \emph{relatively injective} over $\calZ^{\bfx}$.
\begin{theo}\label{theo:KZ}
	Let $L\in \calH^{\bfx}\mof$ be a simple $\calH$-module. Then, $\bfV(\calC^{\bfx}L)\neq 0$ holds if and only if the projective cover of $L$ lies in $\calH^{\bfx}\proj\cap \calH^{\bfx}\inj_{\calZ^{\bfx}}$. 
\end{theo}
\begin{rema}
	It has been shown in \cite{liu22} that the idempotent subalgebra $1_{\bfV}\calH^{\bfx}1_{\bfV}$ is isomorphic to the completion of a certain affine Hecke algebra along a central character. However, we do not need this identification. 
\end{rema}

\printbibliography

\end{document}